\documentclass[11pt, letterpaper]{amsart}

\usepackage{tikz-cd}

\usepackage[utf8]{inputenc}
\usepackage[OT2]{fontenc}
\usepackage[T1]{fontenc}

\usepackage{palatino}

\usepackage[english]{babel}
\usepackage{amsmath}
\usepackage{amsfonts}
\usepackage{mathrsfs}
\usepackage{amssymb}
\usepackage{stmaryrd}
\usepackage{chemfig}
\usepackage{mathtools}
\usepackage{amsthm}
\usepackage{hyperref}
\usepackage{wasysym}
\usepackage{xcolor}

\theoremstyle{plain}
\newtheorem{theo}{Theorem}[section]
\newtheorem{prop}[theo]{Proposition}
\newtheorem{coro}[theo]{Corollary}
\newtheorem{lemma}[theo]{Lemma}
\newtheorem{sublemma}[theo]{Sublemma}

\theoremstyle{definition}
\newtheorem{madeflemma}[theo]{Definition/Lemma}
\newtheorem{madef}[theo]{Definition}
\newtheorem{nota}[theo]{Notation}
\newtheorem{computation}[theo]{Computation}
\newtheorem*{nota*}{Notation}
\newtheorem{rmq}[theo]{Remark}
\newtheorem{ex}[theo]{Example}

\theoremstyle{remark}
\newcommand\cyrillic[1]{{\fontencoding{OT2}\fontfamily{wncyr}\selectfont #1}}
\newcommand\mathcyr[1]{\text{\cyrillic{#1}}}
\newcommand{\sh}{\textnormal{\mathcyr{Sh}}}

\renewcommand{\r}[1]{\textcolor{black}{#1}}
\newcommand{\red}[1]{\textcolor{black}{#1}}

\newcommand{\aff}{\mathrm{aff}}
\newcommand{\proj}{\mathrm{proj}}
\newcommand{\wond}{\mathrm{wond}}

\newcommand{\zero}{\hat{0}}
\newcommand{\un}{\hat{1}}

\newcommand\numberthis{\addtocounter{equation}{1}\tag{\theequation}}

\renewcommand{\Bar}{\mathrm{B}}
\newcommand{\Basis}{\mathrm{B}}

\newcommand{\MM}{\mathsf{MM}}

\newcommand{\Ch}{\mathrm{Ch}\,}
\renewcommand{\d}{\mathrm{d}}

\newcommand{\Q}{\mathbb{Q}}
\renewcommand{\L}{\mathcal{L}}
\newcommand{\N}{\mathbb{N}}
\newcommand{\Z}{\mathbb{Z}}
\newcommand{\G}{\mathcal{G}}
\renewcommand{\S}{\mathcal{S}}
\newcommand{\A}{\mathcal{A}}
\newcommand{\I}{\mathcal{I}}
\renewcommand{\H}{\mathcal{H}}
\newcommand{\At}{\mathrm{At}}
\newcommand{\Fact}{\mathrm{Fact}}

\newcommand{\iso}{\mathrm{iso}}
\newcommand{\Symcat}{\mathrm{Sym}}
\newcommand{\nucat}{\mathcal{V}}
\newcommand{\fcat}{\mathcal{F}}
\newcommand{\F}{\mathfrak{F}}
\renewcommand{\C}{\mathcal{C}}
\newcommand{\D}{\mathcal{D}}
\newcommand{\E}{\mathcal{E}}
\newcommand{\LBSi}{\textbf{LBS}_{\mathrm{irr}}}
\newcommand{\LBS}{\mathfrak{LBS}}
\newcommand{\LBSm}{\mathfrak{LBSmod}}
\newcommand{\LBScat}{\textbf{LBS}}
\renewcommand{\deg}{\mathrm{deg}}
\newcommand{\Hom}{\mathrm{Hom}}
\renewcommand{\ker}{\mathrm{ker}\,}
\newcommand{\Mor}{\mathrm{Mor}}
\newcommand{\im}{\mathrm{Im}\,}
\newcommand{\Id}{\mathrm{Id}}
\newcommand{\rk}{\mathrm{rk}}

\newcommand{\fr}{\mathrm{fr}}
\newcommand{\Op}{\mathfrak{Op}}
\newcommand{\op}{\mathrm{op}}

\newcommand{\Gen}{\mathrm{Gen}}

\newcommand{\lt}{\mathsf{Lt}}

\newcommand{\gCR}{\mathrm{\textbf{grComRing}}}

\newcommand{\Ind}{\mathrm{Ind}}
\newcommand{\ind}{\mathrm{ind}}

\newcommand{\Comp}{\mathrm{Comp}}

\newcommand{\Sym}{\mathbb{S}}

\newcommand{\FY}{\mathbb{FY}}
\newcommand{\OS}{\mathbb{OS}}
\newcommand{\OSbar}{\overline{\mathbb{OS}}}
\newcommand{\OSbara}{\overline{\mathrm{OS}}}
\newcommand{\OSa}{\mathrm{OS}}

\renewcommand{\P}{\mathbb{P}}

\renewcommand{\O}{\mathcal{O}}

\newcommand{\B}{\mathcal{B}}

\newcommand{\Cov}{\mathrm{Cov}}

\newcommand{\indtotor}{ \vartriangleleft^* }

\newcommand{\EL}{\mathrm{EL}}

\newcommand{\M}{\mathfrak{M}}

\newcommand{\PD}{\mathrm{PD}}

\newcommand{\odd}{\mathrm{odd}}

\newcommand{\Cx}{\mathbb{C}}

\newcommand{\cd}{\mathrm{cd}}

\newcommand{\Fil}{\mathrm{F}}

\newcommand{\n}{\mathrm{n}}
\newcommand{\Adm}{\mathrm{Adm}}
\newcommand{\K}{\mathbb{K}}

\newcommand{\FYa}{\mathrm{FY}}

\renewcommand{\int}{\mathrm{int}}
\newcommand{\ext}{\mathrm{ext}}
\newcommand{\Tree}{\mathcal{T}}
\newcommand{\Mop}{\mathbb{M}}
\newcommand{\tot}{\mathrm{tot}}

\title{Matroids, Feynman categories, and Koszul duality}
\author{Basile Coron}
\date{}

\begin{document}
\maketitle
\begin{abstract}
We show that various combinatorial invariants of matroids such as Chow rings and Orlik--Solomon algebras may be assembled into ``operad-like'' structures. Specifically, one obtains several operads over a certain Feynman category which we introduce and study in detail. In addition, we establish a Koszul-type duality between Chow rings and Orlik--Solomon algebras, vastly generalizing a celebrated result of Getzler. This provides a new interpretation of combinatorial Leray models of Orlik--Solomon algebras. \end{abstract}
\tableofcontents
\newpage
\section{Introduction}
This story finds its origin in the celebrated work of De Concini and Procesi \cite{de_concini_wonderful_1995}. In that article the authors construct special compactifications for every projective complement $\P(\Cx^n) \setminus \bigcup_{H \in \H} \P H$ of some hyperplane arrangement $\H$. Those compactifications are called ``wonderful'' because the complement of $\P(\Cx^n) \setminus \bigcup_{H \in \H} \P H$ in the compactification is a divisor with normal crossings. Wonderful compactifications are obtained by successively blowing up some of the intersections of hyperplanes of $\H$. \textcolor{black}{The choice of the set of intersections to blow up, denoted by $\G$,} matters and the resulting blown up space $\overline{Y}_{\H, \G}$ is a wonderful compactification when $\G$ is a \emph{building set} of the lattice $\L_{\H} = \{\bigcap_{H \in I} H, \, \, I \subset \H\}$, ordered by reverse inclusion. The building set condition on $\G$ ensures that all the non-transversal intersections will be blown up. \\

Wonderful compactifications are naturally stratified by the exceptional divisors obtained in the process of blowing up. There is one exceptional divisor $\D_{G}$ for each element $G$ of~$\G$. Those divisors have normal crossings and an intersection of divisors $\D_{G_1},\textcolor{black}{\ldots}, \D_{G_n}$ is nonempty exactly when $\{G_1,\textcolor{black}{\ldots} , G_n\}$ forms a \emph{nested set} of~$\G$. The non-empty intersections of divisors are the closed strata of the stratification. De Concini and Procesi discovered that the closed strata are in fact isomorphic to products of ``smaller'' wonderful compactifications. With this in mind, simply considering inclusions of the strata in the compactifications leads to interesting additional structures. For instance, when looking at the family of braid arrangements $\mathsf{Braid_n}$ consisting of the diagonal hyperplanes $\{z_i = z_j\} \subset \Cx^n$, together with the unique minimal building set, this gives us the well-known structure of an operad on the Deligne--Mumford compactification\textcolor{black}{s} of moduli spaces of pointed curves of genus zero \textcolor{black}{(see Getzler \cite{Getzler_1994})}. Passing to the homology gives the operad of hypercommutative algebras called $\mathsf{Hypercom}$ \textcolor{black}{(see Manin \cite{Manin_1999} for a textbook reference on this subject)}. When looking at the family of boolean arrangements $\mathsf{Bool_n}$ consisting of coordinate hyperplanes in $\Cx^n$, together with the unique maximal building set, we also get an operad-like structure, this time on the so-called \textcolor{black}{Losev--Manin} spaces \textcolor{black}{(see Losev and Manin \cite{LM_2000})}.\\

The aim of this article is to develop a formalism enabling us to see the family of all possible wonderful compactifications of all possible hyperplane arrangements, with structural morphisms given by inclusions of strata, as one big ``operad''. This way, we will be able to interpret commonalities between the above examples as mere consequences of the properties enjoyed by the unified structure. For instance, we shall be interested in finding a unified presentation by generators and quadratic relations of this operad, as well as a minimal model, which will be done by Koszul duality theory (see Loday and Vallette \cite{LV_2012} for an introduction). Fortunately, recent years have seen the advent of general theories that were successful in developing a language as well as methods to deal with ``operad-like'' structures, such as generalized operads (Borisov and Manin \cite{Borisov2008}), patterns (Getzler \cite{Getzler2009}), operadic categories (Batanin and Markl \cite{Batanin_2015}) and so on. One of those theories is that of Feynman categories (Kaufmann and Ward \cite{kaufmann_feynman_2017}), which we will use in this article to fulfil our goal. \\

Instead of dealing with hyperplane arrangements we shall be working at the more general level of matroids (see Welsh \cite{welsh_matroid_1976} for an introduction), which form a combinatorial abstraction of hyperplane arrangements. For our purposes it will be enough to consider only simple loopless matroids, which can be axiomatized by geometric lattices. This axiomatization is the most convenient to us and is the one we will be working with in this paper. \\

The wonderful compactifications do depend on the hyperplane arrangement itself, together with the choice of the building set, but their cohomology algebra\textcolor{black}{s} only depend on the intersection lattice and its building set. In  \cite{feichtner_chow_2003}, Feichtner and Yuzvinsky introduced a generalization of those cohomology rings for every pair \textcolor{black}{consisting} of a geometric lattice $\L$ and a building set $\G \subset \L$, which we will denote by $\FYa(\L, \G)$. In the case where $\G$ is equal to $\L \setminus \{\zero \}$, the ring \textcolor{black}{$\FYa(\L, \L \setminus \{\zero \})$} is \textcolor{black}{commonly referred to} as the combinatorial Chow ring of $\L$. \textcolor{black}{In general, the rings $\FYa(\L, \G)$} have been extensively studied and are known to satisfy very strong properties. For instance even though these rings are not necessarily cohomology rings of projective complex varieties, they all have a Hodge theory, meaning that they satisfy Poincaré duality as well as the Kähler package (see Adiprasito, Huh, and Katz \cite{Huh_2018}, and, Pagaria and Pezzoli \cite{Pagaria_2021} for general building sets). In \cite{Bibby_2021}, Bibby, Denham and Feichtner show that the morphisms between cohomology algebras of wonderful compactification\textcolor{black}{s} induced by the inclusions of strata can also be generalized to the purely combinatorial setting. This means that we have morphisms between the rings $\FYa(\L, \G)$, indexed by nested sets. \\

In Section \ref{secfeycat} we construct a Feynman category $\LBS$ \textcolor{black}{(those letters standing for lattice, building set respectively)} such that the family of algebras $\{\FYa(\L, \G) \}_{(\L, \G)}$ together with the above morphisms forms a \textcolor{black}{co}operad of type $\LBS$, that is simply a monoidal functor from \textcolor{black}{$\LBS^{\op}$} to some symmetric monoidal category (in the present case the category of graded commutative algebras). In short, the objects of $\LBS$ will be the pairs $(\L, \G)$ with $\G$ a building set of $\L$, and the morphisms will be given by the nested sets. The heart of the problem is to find a suitable ``composition'' of nested sets, which is given in Subsection \ref{secFeycons}. \red{In Subsection \ref{presentation} we give a presentation of the Feynman category $\LBS$. This presentation allows one to describe an operad of type $\LBS$ concretely as a collection of objects $\O(\L,\G)$ (in some symmetric monoidal category) indexed by pairs $(\L, \G)$ with $\G$ a building set of $\L$, together with symmetries, and structural morphisms 
$$ \O([\zero, G], \G\cap [\zero,G])\otimes \O([G,\un], (G\vee\G)\setminus\{G\}) \rightarrow \O(\L, \G),$$ (with $\zero$ and $\un$ the minimum and maximum of $\L$ respectively) which must satisfy some ``associativity'' axioms and compatibility relations (see Corollary \ref{corodescriptionoperad} for the exact statement).\\}

Going back to braid arrangements, in \cite{Getzler_1994} Getzler has shown that $\mathsf{Hypercom}$ is Koszul with Koszul dual the operad $\mathsf{Grav}$ consisting of the cohomology algebras of the projective complements of the braid arrangements, with operations given by residue morphisms. \textcolor{black}{To sum up, Getzler's argument is that the bar construction of $\mathsf{Hypercom}$ can be identified with the second page of the Leray spectral sequence associated to the inclusion of the projective arrangement complement inside the wonderful compactification, and that this spectral sequence converges at the third page (necessarily to $\mathsf{Grav}$) by a mixed Hodge theoretic argument}. Cohomology algebras of projective complements of hyperplane arrangements are called projective Orlik--Solomon algebras and only depend on the intersection lattice (see Orlik-Solomon \cite{OS_1980}). They can be generalized to arbitrary geometric lattices. We will denote those algebras by $\OSbara(\L)$, for $\L$ any geometric lattice. \\

In this article we show that the family of projective Orlik--Solomon algebras $\{ \OSbara(\L) \}$ also has a \textcolor{black}{co}operadic structure over $\LBS$, with morphisms given by combinatorial generalizations of residue morphisms. When restricted to partition lattices (intersection lattices of braid arrangements) this gives back the linear dual of the operad $\mathsf{Grav}$. Additionally, in \cite{Bibby_2021} the authors show that the combinatorial Chow rings can be assembled to form a combinatorial model $B(\L, \G)$ of $\OSbara(\L)$:
$$\OSbara(\L) \xrightarrow{\sim} B(\L, \G). $$
\textcolor{black}{This model is a combinatorial generalization of the second page of the spectral sequence mentioned in the preceding paragraph. In Section \ref{seckoszul} we explain that $\{B(\L, \G)\}$ is also an operad over $\LBS$ and can in fact be identified (just as in the classical operadic case) with the bar construction of the operad $\{\FYa^{\vee}(\L, \G)\}$ (where $(-)^{\vee}$ denotes the linear dual). This immediately implies the following corollary.} 
\begin{coro}
The operad $\{\FYa^{\vee}(\L,\G)\}$ is Koszul with Koszul dual $\{\OSbara(\L)\}$.
\end{coro}
We give an alternative proof of this fact by implementing a Gröbner bases machinary for operads over $\LBS$ and proving the following results.
\begin{prop}
If an operad over $\LBS$ admits a quadratic Gröbner basis then it is Koszul.
\end{prop}
\begin{theo}
The operad $\{\FYa^{\vee}(\L, \G) \}$ admits a quadratic Gröbner basis.
\end{theo}
\textcolor{black}{This is also an extension of a known proof of the Koszulness of $\mathsf{Hypercom}$ (see Dotsenko and Khoroshkin \cite{Dotsenko_2012})}. Finally, let us mention that evidence of ``operad-like'' structures related to combinatorics of building sets/nested sets were already highlighted by Forcey and Ronco \cite{Forcey_Ronco_2022}, and Rains \cite{Rains_2010}, in specific settings and using different formalisms. \\

Here is the general layout of the article.  \\

Section \ref{secprelim} is devoted to introducing the combinatorial ingredients that will be used in the construction of the Feynman category $\LBS$.\\

Section \ref{secfeycat} is the core of this paper. This is where we define the Feynman category $\LBS$. We also prove the important fact that $\LBS$ admits a graded presentation.\\

In Section \ref{secoperads} we show that the families $\{\FYa(\L, \G) \}$, $\{\OSa(\L)\}$, $\{\OSbara(\L)\}$ form operads over the Feynman category $\LBS$. \\

In Section \ref{secgrobner} we develop a theory of Gröbner bases for operads over $\LBS$ and prove that $\{\FYa(\L, \G)\}$ admits a quadratic Gröbner basis. \\

In Section \ref{seckoszul} we show that $\LBS$ is a cubical Feynman category which implies that there is a Koszul duality theory for operads over $\LBS$. We then show that $\{\FYa(\L, \G)\}$ is Koszul via two different methods. \\

Finally, in Section \ref{secdirec} we give some last remarks toward possible generalizations and modifications of $\LBS$. \\

\textbf{Acknowledgements.} The author would like to thank Vladimir Dotsenko for his invaluable support and countless crucial discussions. We also thank Clément Dupont for his generous availability and many interesting conversations, as well as Bérénice Delcroix-Oger. 

This research is part of the author's PhD. The author would like to thank Eva Maria Feichtner and Karim Adiprasito for agreeing to be his thesis referees. 

The author wishes to express gratitude to the anonymous referee for their exceptional efforts, which have tremendously improved this article. 

This research was supported by the University of Strasbourg Institute for Advanced Study through the French national program ``Investment for the future'' [IdEx-Unistra, fellowship USIAS-2021-061 of V.~Dotsenko] and by the French national research agency [grant ANR-20-CE40-0016].

\newpage

\section{Combinatorial preliminaries}\label{secprelim}
In this section we introduce the main combinatorial objects which will be used throughout this paper. \\

It is important to note that in this article we work at a strictly combinatorial level. However, many of the main protagonists in this story have a geometric origin and although the geometric picture is not formally required, it is our main source of inspiration and therefore we will try to draw this picture whenever possible.
\subsection{Lattices, building sets and nested sets}

\begin{madef}[Lattice]
A finite poset $\L$ is called a \textit{lattice} if every pair of elements in $\L$ admits a supremum and an infimum.
\end{madef}
The supremum of two elements $G_1, G_2$ is denoted by $G_1 \vee G_2$ and called their \textit{join}, while their infimum is denoted by $G_1 \wedge G_2$ and called their \textit{meet}.
\begin{rmq}
Since $\L$ is \textcolor{black}{assumed} to be finite, having supremums and infimums for pairs of elements implies having supremums and infimums for any subset $S$ of $\L$, which will be denoted by $\bigvee S$ and $\bigwedge S$ respectively. As a consequence, every lattice admits an upper bound (the supremum of $S = \L$) and a lower bound (the infimum of $S = \L$) which will be denoted by $\un$ and $\zero$ respectively. 
\end{rmq}
\begin{madef}[Geometric lattice]
A finite lattice $(\L, \leq)$ is said to be \textit{geometric} if it satisfies the following properties:
\begin{itemize}
\item For every pair of elements $G_1 \leq G_2$, all the maximal chains of elements between $G_1$ and $G_2$ have the same \textcolor{black}{length}. \textit{(Jordan--Hölder property)}
\item The rank function $\rho: \L \rightarrow \N$, which assigns to any element $G$ of $\L$ the \textcolor{black}{length} of any maximal chain of elements \textcolor{black}{between} $\zero$ and $G$, satisfies the inequality
\begin{equation*}
\rho(G_1 \wedge G_2) + \rho(G_1 \vee G_2) \leq \rho(G_1) + \rho(G_2)
\end{equation*}
for every $G_1$, $G_2$ in $\L$. \textit{(Sub-modularity)}
\item Every element in $\L$ can be obtained as the supremum of some set of atoms (i.e. elements of rank~$1$). \textit{(Atomicity)}
\end{itemize}
\end{madef}
\textcolor{black}{One of the  reasons to focus on this particular class of posets comes from the study of geometric objects called hyperplane arrangements. A hyperplane arrangement is a finite collection of hyperplanes in a shared finite dimensional vector space. To each hyperplane arrangement $\H$ one can associate a lattice called the intersection lattice of $\H$, defined as the set of all possible intersections of hyperplanes in $\H$ ordered by reverse inclusion (the join in this lattice being given by the intersection of subspaces). Essentially, this lattice captures the linear dependencies between the hyperplanes of the hyperplane arrangement, but without keeping track of the exact linear relations. The key point is that the intersection lattice of a hyperplane arrangement is always a geometric lattice. The Jordan--Hölder property arises from the existence of dimension in linear algebra, submodularity comes from Grassmann's identity and the atomicity is inherent in the definition of the intersection lattice.} In fact, one may think of geometric lattices as a combinatorial abstraction of hyperplane arrangements. In addition, \textcolor{black}{the datum of a geometric lattice} is equivalent to the datum of a loopless simple matroid via the lattice of flats construction (see \cite{welsh_matroid_1976} for a reference on matroid theory) and therefore it has connections to many other areas in mathematics (graph theory for instance). \\

Here is a list of some important well-known geometric lattices.
\begin{ex}\label{exgeolatt}
\leavevmode
\begin{itemize}
\item If $X$ is any finite set, the set $\mathcal{P}(X)$ of subsets of $X$ ordered by inclusion is a geometric lattice with join the union and meet the intersection. It is the intersection lattice of the hyperplane arrangement of coordinate hyperplanes in $\Cx^{X}$. Those geometric lattices are called boolean lattices and denoted by $\mathcal{B}_{X}$.
\item If $X$ is any finite set, the set $\Pi_{X}$ of partitions of $X$, \textcolor{black}{viewed as equivalence relations on $X$, is ordered by inclusion of the corresponding subsets of $X\times X$}. It is the intersection lattice of the so-called \textit{braid arrangement} which consists of the diagonal hyperplanes $\{z_i = z_j\}$ in $\Cx^{X}$. Those geometric lattices are called partition lattices.
\item If $\textcolor{black}{\Gamma} = (V, E)$ is any graph one can construct the graphical matroid $M_{\textcolor{black}{\Gamma}}$ associated to $\textcolor{black}{\Gamma}$ and then consider $\L_{\textcolor{black}{\Gamma}}$ the lattice of flats associated to $M_{\textcolor{black}{\Gamma}}$ (see \cite{welsh_matroid_1976} for the details of \textcolor{black}{that} construction). Those \textcolor{black}{geometric} lattices are said to be \emph{graphical}. This family of geometric lattices contains the \textcolor{black}{previous two} because $\mathcal{B}_{X}$ is the lattice associated to any tree with \textcolor{black}{edge set} $X$ and $\Pi_{X}$ is the lattice associated to the complete graph with \textcolor{black}{vertex set} $X$. For any graph $\textcolor{black}{\Gamma} = (V, E)$ the geometric lattice $\L_{\textcolor{black}{\Gamma}}$ is the intersection lattice of the hyperplane arrangement $\{\{z_{u} = z_{v}\}, (u,v) \in E\}$ in $\mathbb{C}^{V}$.
\end{itemize}
\end{ex}

We have the following important fact about geometric lattices.
\begin{prop}[\cite{welsh_matroid_1976}]
Let $(\L, \leq)$ be a geometric lattice. For every $G_1 \leq G_2 \in \L$, the interval $[G_1, G_2] =\{G \in \L \, | \, G_1\leq G \leq G_2\}$ ordered by the restriction of $\leq$ is a geometric lattice.
\end{prop}
In the rest of this article every lattice will be assumed to be geometric unless stated otherwise.
\begin{madef}[Building set]\label{defbuilding}
Let $\L$ be a geometric lattice. A \textit{building set} $\G$ of $\L$ is a subset of $\L \setminus \{\zero\}$ such that for every element $X$ of $\L$ the morphism of posets
\begin{equation}\label{isobuilding}
\prod_{G \in \textcolor{black}{\Fact_{\G}(X)}}[\zero, G] \xrightarrow{\vee} [\zero, X]
\end{equation}
is an isomorphism (where \textcolor{black}{$\Fact_{\G}(X)$} denotes the set of maximal elements of $\G \cap [\zero,\textcolor{black}{X}])$).
\end{madef}
The elements of \textcolor{black}{$\Fact_{\G}(X)$} will be called the \emph{factors} of $X$ in $\G$. \textcolor{black}{In plain English the above definition means that for every element $X$, the elements below $X$ can be uniquely written as a join of elements below the factors of $X$. Heuristically, one may think of building sets as an axiomatization of the notion of connectedness (in graph theory, topology, etc). The building set prescribes which elements are ``connected'' and for every $X$ the factors of $X$ are the ``connected components'' of $X$.}
\begin{madef}[Built lattice]
The datum of a lattice $\L$ and a building set $\G$ of $\L$ will be called a \textit{built lattice}. If $\G$ contains $\un$ we will say that $(\L, \G)$ is \textit{irreducible}.
\end{madef}
The definition of a building set makes sense for a larger class of posets, as shown in \cite{feichtner_chow_2003}, but in this paper we will restrict ourselves to the case of geometric lattices. In this particular context, building sets are geometrically motivated by the construction of wonderful compactifications for hyperplane arrangement complements. In a nutshell, building sets are sets of intersections of a hyperplane arrangement that one can successively blow up in order to obtain a wonderful compactification of its complement (see \cite{de_concini_wonderful_1995} for more details). Each blowup creates a new exceptional divisor, so the wonderful compactification is equipped with a family of irreducible divisors indexed by $\G$. This family of divisors forms a normal crossing divisor when $\G$ is a building set. \\

There are a few key examples to keep in mind throughout this story.
\begin{ex}\label{exbs}
\leavevmode
\begin{itemize}
\item Every lattice $\L$ admits $\L\setminus\{\zero\}$ as a building set. \textcolor{black}{We will refer to that building set as the maximal building set of $\L$}. 
\item Every lattice $\L$ also admits a unique minimal building set which consists of all the elements $G$ of $\L$ such that $[\zero,G]$ is not a product of proper subposets.
\item From the definition one can see that a building set of some lattice $\L$ must contain all the atoms of $\L$. If $\L$ is a boolean lattice (see Example \ref{exgeolatt}) then its set of atoms is in fact a building set (the minimal one). This fact characterizes boolean lattices.
\item If $\L$ is the lattice of partitions of some finite set (see Example \ref{exgeolatt}) then the subset of partitions with only one block having strictly more than one element is a building set of $\L$. This is the minimal building set of $\L$.
\item If $\L$ is a graphical lattice (see Example \ref{exgeolatt}) then the set of elements of $\L$ corresponding to sets of edges which are connected \textcolor{black}{(in the sense of graph theory)} forms a building set of $\L$, called the graphical building set of $\L$. This family of examples contains \textcolor{black}{the previous two} (by considering totally disconnected graphs for the former and complete graphs for the latter). \textcolor{black}{In general the graphical building set may be different from the minimal building set. For instance for star-shaped graphs every subset of edges is connected and therefore in that case the graphical building set coincides with the maximal building set.}
\item Alternatively, if $\textcolor{black}{\Gamma} = (V,E)$ is a graph one can consider the boolean lattice $\B_{V}$. This lattice has a building set made up of the ``tubes'' of $\Gamma$, that is, sets of vertices of $\Gamma$ such that the induced subgraph on those vertices is connected. This leads to the notion of graph associahedra introduced in \cite{Carr_Devadoss_2004}.
\end{itemize}
\end{ex}
A key fact about building sets is that any interval \textcolor{black}{$[X_1, X_2]$} in some built lattice $(\L, \G)$ admits an ``induced'' building set \textcolor{black}{which we now describe}. We start by introducing a useful notation.
\begin{nota}
For any element \textcolor{black}{$X$} of some lattice $\L$ and any subset \textcolor{black}{$S$} of $\L$, we denote by \textcolor{black}{$X\vee S$} the set of elements of $\L$ which can be obtained as the join of \textcolor{black}{$X$} and some element of \textcolor{black}{$S$}. We also denote by $S_{\leq X}$ the set of elements in $S$ below $X$. 
\end{nota}
{\color{black}\begin{madeflemma}[Induced building set]
Let \textcolor{black}{$X_1 < X_2$} be two elements in some lattice $\L$ and let $\G$ be a building set of $\L$. The set $$(X_1\vee \G) \cap [X_1, X_2] \setminus \{X_1\},$$ which we will denote by $\Ind_{[X_1, X_2]}(\G)$, is a building set of $[X_1, X_2]$. We will call this building set the \emph{induced building set} on $[X_1, X_2]$.
\end{madeflemma}
\begin{proof}
That result can be extracted from \cite{Bibby_2021} Lemma 2.8.5. For completeness we provide an explicit proof. Let $X$ be some element in $[X_1, X_2]$. By isomorphism \eqref{isobuilding} applied to $X$ viewed in $\L$ we have the isomorphisms
\begin{align*}
    [X_1, X] & \simeq \prod_{G \in \Fact_{\G}(X)\setminus \Fact_{\G}(X_1)}[\bigvee \Fact_{\G}(X_1)_{\leq G} ,G] \\
             & \simeq \prod_{G \in \Fact_{\G}(X)\setminus \Fact_{\G}(X_1)}[X_1, G\vee X_1].
\end{align*}
What is left is to prove the equality
\begin{equation}\label{eqindproof}
    \Fact_{\Ind_{[X_1, X_2]}(\G)}(X) = X_1 \vee (\Fact_{\G}(X)\setminus \Fact_{\G}(X_1)).
\end{equation}
It is clear that the right-hand side of \eqref{eqindproof} is included in $\Ind_{[X_1, X_2]}(\G)\cap [X_1, X]$. Let $G$ be some element of $\Fact_{\G}(X)\setminus \Fact_{\G}(X_1)$ such that $X_1 \vee G$ is less than or equal to some element $X_1 \vee G' \in \Ind_{[X_1, X_2]}(\G)\cap [X_1, X]$ with $G'$ some element in $\G$. The element $G'$ must be below one of the factors of $X$. By isomorphism \eqref{isobuilding} this factor must be $G$, which implies that we have in fact $X_1 \vee G = X_2\vee G$. In other words $X_1 \vee G$ is maximal in $\Ind_{[X_1, X_2]}(\G)\cap [X_1, X]$, and we have just proved that the right-hand side of \eqref{eqindproof} is included in the left-hand side. On the other hand if $G\vee X_1$ is an element of $\Fact_{\Ind_{[X_1, X_2]}(\G)}(X)$, with $G$ some element of $\G$, then by maximality of $G\vee X_1$ we get that $G$ is maximal in $\G\cap [\zero, X]$ which proves that $G\vee X_1$ belongs to $X_1 \vee (\Fact_{\G}(X)\setminus \Fact_{\G}(X_1))$. 
\end{proof}}
We will often write $\Ind(\G)$ instead of  $\Ind_{[X_1, X_2]}(\G)$ if the interval can be deduced from the context. \textcolor{black}{Note that the above operation of restricting a building set to an interval $[X_1, X_2]$ can be broken down into first restricting to $[\zero, X_2]$, which is just taking the intersection with $[\zero, X_2]$, and then restricting to $[X_1, \un (=X_2)]$, which is much less innocent as it actually modifies the elements of the building set.} We have the following key lemma.
\begin{lemma}\label{lemmadoubleinduction}
For any elements $X_1 \leq X_2 \leq X_3 \leq X_4$ in some lattice $\L$ with building set $\G$, we have the equality of building set\textcolor{black}{s}
\begin{equation*}
\Ind_{[X_2, X_3]}(\Ind_{[X_1, X_4]}(\G)) = \Ind_{[X_2, X_3]}(\G).
\end{equation*}
\end{lemma}
{\color{black}
\begin{proof}
\begin{align*}
\Ind_{[X_2, X_3]}(\G) &= (X_2\vee \G) \cap [X_2, X_3] \setminus \{X_2\} \\
&= ((X_2\vee X_1)\vee \G) \cap [X_2, X_3] \setminus \{X_2\} \\ 
&= (X_2\vee (X_1 \vee \G)) \cap [X_2, X_3] \setminus \{X_2\} \\
&= (X_2\vee ((X_1 \vee \G) \cap [X_1, X_4] \setminus \{X_1\})) \cap [X_2, X_3] \setminus \{X_2\} \\
&= (X_2\vee \Ind_{[X_1, X_4]}(\G))\cap [X_2, X_3] \setminus \{X_2\} \\
&= \Ind_{[X_2, X_3]}(\Ind_{[X_1, X_4]}(\G)).
\end{align*}
\end{proof}}
\textcolor{black}{For any subset $I$ in some poset we will denote by $\max I$ the set of maximal elements of $I$. If we have $I = \max I$, that is, if the elements of $I$ are pairwise incomparable, the subset $I$ will be called an \textit{antichain}.}

\begin{madef}[Nested set]
Let $(\L, \G)$ be a built lattice. A subset $\S$ of $\G$ is called a \textit{nested set} if for every \textcolor{black}{antichain} $\A$ included in $\S$, the join of the elements of $\A$ does not belong to $\G$ whenever $\A$ contains at least two elements. A nested set $\S \subset \G$ is said to be \red{\textit{spanning}} if it contains $\max \G$.
\end{madef}
\begin{ex}
A chain of elements in some building set $\G$ is always nested and those are the only nested sets of the maximal building set ($\G = \L \setminus \{\zero \}).$
\end{ex}
{\color{black}\begin{ex}\label{exnestedsetshuffletree}
For any integer $n\geq 2$, consider the partition lattice $\Pi_{\{1,\ldots,n\}}$ with its minimal building set (see Example \ref{exbs}). The elements of $\G_{\min}$ can be described by their unique non trivial equivalence class, which is a subset of $\{1,\ldots, n\}$ having at least two elements. To a spanning nested set $\S$ of $(\Pi_{\{1,\ldots,n\}}, \G_{\min})$ one can associate a planar rooted tree $\mathcal{T}$ with leaves indexed by $\{1,\ldots, n \}$ as follow. First, set the inner vertices of $\mathcal{T}$ to be the set $\S$ and declare that some element $G\in \S$ is a parent of some other element $G'\in \S$ if $G$ is covered by $G'$ in $\S$. Declare also that a leaf $i \in \{1,\ldots,n\}$ is the parent of some element $G \in \S$ if $G$ is the unique minimal element of $\S$ containing $i$ (when viewed as a set). Finally, embed $\mathcal{T}$ in the plane in the unique way such that at every vertex $v$ of $\mathcal{T}$, the minima of the sets of leaves of the parent of $v$ are increasing when read from left to right. A rooted planar tree satisfying this property is called a shuffle tree. One can check that this construction establishes a bijection between the spanning nested sets of $(\Pi_{\{1,\ldots,n\}}, \G_{\min})$ and the shuffle trees with leaves labelled from $1$ to $n$ and such that every vertex has at least two parents. 
\end{ex}}
Geometrically, nested sets correspond to sets of divisors in the wonderful compactification which have a nontrivial intersection. There are two crucial lemmas regarding nested sets. \textcolor{black}{In the rest of this article, by nested antichain we mean an antichain which is nested.} 
\begin{lemma}[\cite{FK_2004} Proposition 2.8]\label{factors}
Let $\G$ be a building set of a geometric lattice $\L$ and let $X$ be any element of $\L$. The subset $\Fact_{\G}(X)$ is a nested antichain in $\G$ and furthermore it is the only nested antichain in $\G$ having join $X$.
\end{lemma}
{\color{black}
\begin{proof}
If the join $G$ of a subset $I$ of $\Fact_{\G}(X)$ belongs to $\G$ then $G$ must be below one of the factors of $X$. By isomorphism \eqref{isobuilding} this immediately implies $I = \{G\}$, which proves that $\Fact_{\G}(X)$ is a nested set. It is clear by isomorphism \eqref{isobuilding} that $\Fact_{\G}(X)$ has join $X$. If a nested antichain $\A$ has join $X$ then by \cite{FK_2004} Proposition 2.8 (1) $\implies$ (2), we have that $\A$ must be equal to $\Fact_{\G}(X)$, which proves the uniqueness statement. 
\end{proof}}
\begin{lemma}\label{lemmaforest}
A nested set of $\L$ is a forest in the Hasse diagram of $\L$. More precisely for every $K$ in \textcolor{black}{$\L\setminus\{\zero\}$}, $\S_{>K}$ is either empty or has a unique minimal element.
\end{lemma}
{\color{black}
\begin{proof}
For every $K$ different from $\zero$, if $\S_{>K}$ is non-empty then by the previous proposition $\min \S_{>K}$ is the set of factors of $\bigvee \min \S_{>K}$. By isomorphism \eqref{isobuilding}, two intervals $[\zero, G_1]$, $[\zero, G_2]$ with $G_1, G_2$ two different factors of $\bigvee \min \S_{>K}$ can only have intersection $\{\zero\}$, which in our case is impossible since those two intervals must contain $K$. As a consequence $\min \S_{>K}$ must be a singleton. 
\end{proof}}
{\color{black}
We next introduce a map $\Comp^{X_2}_{X_1}: \Ind_{[X_1,X_2]}(\G) \rightarrow \G$ for any elements $X_1 \leq X_2$ in $\L$, which we will use later to define the composition of nested sets. One has the following lemma. 
\begin{madeflemma}\label{defcomp}
Let $X_1 \leq X_2$ be two elements of $\L$ and let $Y$ be an element of $\Ind_{[X_1,X_2]}(\G)$. The set of elements $G \in \G$ such that $Y = X_1 \vee G$ has unique maximal element. In addition, this element is also the unique factor of $Y$ in $\G$ which is not a factor of $X_1$. We will denote this element by $\Comp^{X_2}_{X_1}(Y).$
\end{madeflemma}
\begin{proof}
Let $G \in \G$ be such that $Y = X_1 \vee G$ and which is maximal for that property. Denote by $\{G_i, 1 \leq i\leq n\}$ the factors of $X_1$ in $\G$. We have the equalities
\begin{equation*}
Y = X_1\vee G = \left(\bigvee_{1 \leq i \leq n} G_i \right)\vee G = \left(\bigvee_{\substack{1 \leq i \leq n \\ G_i \nleq G}} G_i \right)\vee G.
\end{equation*}
The antichain $\A \coloneqq \{G_i, G_i \nleq G\}\cup \{G\}$ is nested because if an antichain $\A'\subset \A$ of cardinality greater than $1$ has a join in $\G$ then the nestedness of $\Fact_{\G}(X_1)$ implies that $\A'$ must contain $G$, but then $\bigvee \A'$ contradicts the maximality of $G$. By Lemma \ref{factors} the set $\A$ is exactly the set of factors of $Y$ in $\G$, which concludes the proof. 
\end{proof}
\begin{ex}\label{excompmax}
If $\G$ is the maximal building set of $\L$ then for all $X_1\leq X_2$ in $\L$ the induced building set on $[X_1, X_2]$ is the maximal building set of $[X_1, X_2]$, and the map $\Comp_{X_1}^{X_2}$ is simply the inclusion of $[X_1, X_2]\setminus\{X_1\}$ in $\G$. 
\end{ex}
\begin{ex}
Consider the partitions $G_1 = 13|2|4|5|6$, $G_2 = 25|1|3|4|6$ in the minimal building set of $\Pi_{\{1,\ldots,6\}}$ and define $X \coloneqq G_1 \vee G_2 = 13|25|4|6$. By definition the element $Y \coloneqq X \vee 56|1|2|3|4$ belongs to $\Ind_{[X, \un]}(\G_{\min})$. One has $Y = 13|4|256$ and thus one can see that the only factor of $Y$ which is not a factor of $X$ is $256$. In other words, we have $\Comp_{X_1}^{\un}(Y) = 256|1|3|4$.
\end{ex}
If $X_1,X_2$ can be deduced from the context we will omit them. We have the two lemmas. 
\begin{lemma}\label{rmkinj}
For all $X_1 \leq X_2 \in \L$, the map $\Comp_{X_1}^{X_2}$ is injective.
\end{lemma}
\begin{proof} 
The map $\Comp_{X_1}^{X_2}$ has a left inverse given by taking the join with $X_1$. 
\end{proof}
\begin{lemma}\label{lemmacompincr}
For all $X_1 \leq X_2 \in \L$, the map $\Comp_{X_1}^{X_2}$ is increasing.
\end{lemma}
\begin{proof}
Let $G_1\leq G_2$ be two elements in some induced building set $\Ind_{[X_1, X_2]}(\G)$. The element $\Comp(G_1)$ belongs to $\G$ and is below $G_2$, which implies that it must be below one of the factors of $G_2$. Those factors are $\Comp(G_2)$ and some factors of $X_1$. Since $\Comp(G_1)$ is not below any of the factors of $X_1$, it must be below $\Comp(G_2)$. 
\end{proof}}
\subsection{Combinatorial invariants}
To the objects introduced in the previous section, one can associate various rings that generalize cohomology rings in the realizable case, \textcolor{black}{that is, when the geometric lattice is the intersection lattice of a complex hyperplane arrangement.}
\subsubsection{The Feichtner--Yuzvinsky rings}

\begin{madef}
For every built lattice $(\L, \G)$ we define the Feichtner--Yuzvinsky graded commutative ring $\FYa(\L, \G)$ by
\begin{equation*}
\FYa(\L, \G) = \Z[x_G, \,G \in \G]/ \I_{\aff},
\end{equation*}
with all the generators in degree $2$, and $\I_{\aff}$ the ideal generated by elements
\begin{equation*}
\sum_{G \geq H}x_G
\end{equation*}
for every atom $H$, and elements
\begin{equation*}
\prod_{G \in X}x_G
\end{equation*}
for every set $X \subset \G$ which is not nested.
\end{madef}
In the realizable case, the ring $\FYa(\L, \G)$ is the cohomology ring of the wonderful compactification \textcolor{black}{of the projectivized arrangement complement, with choice of} building set $\G$ (see \cite{de_concini_wonderful_1995} for the computation of the cohomology ring). Those rings \textcolor{black}{were defined} for arbitrary built lattices by Feichtner and Yuzvinsky in \cite{feichtner_chow_2003}. \\

The Feichtner--Yuzvinsky rings admit two other useful presentations.
\begin{prop}
For every built lattice $(\L, \G)$ we have the other classical presentation
\begin{equation*}
\FYa(\L, \G) \simeq \Z[x_G, \, G \in \G \setminus \{\un \}] / \I_{\proj}
\end{equation*}
where $\I_{\proj}$ is the ideal generated by elements
\begin{equation*}
\sum_{\un > G \geq H_1}x_G - \sum_{\un > G \geq H_2}x_G
\end{equation*}
for every pair of atoms $H_1$ and $H_2$, and elements
\begin{equation*}
\prod_{G \in X}x_G
\end{equation*}
for every set $X \subset \G\setminus \{\un\}$ which is not nested.

Additionally, we have the presentation
\begin{equation*}
\FYa(\L, \G) \simeq \Z[h_G, \, G \in \G] / \I_{\wond}
\end{equation*}
where $\I_{\wond}$ is the ideal generated by relations
\begin{equation*}
h_H
\end{equation*}
for every atom $H$ and
\begin{equation*}
\prod_{G' \in \A} (h_G - h_{G'})
\end{equation*}
for every $G \in \G$ and $\A$ an antichain in $\G$ such that $\bigvee \A$ is equal to $G$. The change of variable between the last presentation and the defining presentation is given by $$h_G = \sum_{G' \geq G}x_{G'}.$$
\end{prop}
The first (defining) presentation will be called the \textit{affine} presentation, the second the \textit{projective} presentation and the last one the \textit{wonderful} presentation. The first two presentations appear in \cite{feichtner_chow_2003} (as a definition) while the second appeared first in \cite{EHKR_2010} for the braid arrangement and in \cite{Backman_Spencer_Eur_2020} for general maximal building sets. It is widely used in \cite{Pagaria_2021}.
\begin{proof}
The proof can be found in \cite{Pagaria_2021} (Theorem 2.9).
\end{proof}
In \cite{feichtner_chow_2003}, the authors address the issue of finding a Gröbner basis \textcolor{black}{of $\I_{\aff}$} (see \cite{BW_1993} for a reference on Gröbner bases) and they show that when considering any linear order on generators refining the reversed order on $\G$, \textcolor{black}{and the induced degree lexicographic ordering on monomials}, although the elements  defining $\I_{\aff}$ do not form a Gröbner basis in general, one can still describe a fairly manageable Gröbner basis.
\begin{theo}[\cite{feichtner_chow_2003} Theorem 2, \textcolor{black}{\cite{feichtner_chow_2003} Corollary 1}]\label{theogrobnerFY}
Elements of the form $$(\prod_{G \in \S}x_G) h_{G'}^{\rho(G') - \rho(\bigvee S)}$$ with $\S$ any nested set and $G'$ any element of $\G$ satisfying $G' > \bigvee \S$ ,  together with the usual $\prod_{G \in X} x_G$ for every non-nested set $X$, form a Gröbner basis of \textcolor{black}{$\I_{\aff}$} for \textcolor{black}{the degree lexicographic ordering on monomials associated to} any linear order on generators refining the reversed order of $\L$. The normal monomials with respect to this Gröbner basis are monomials of the form
\begin{equation*}
x_{G_1}^{\alpha_1}... x_{G_n}^{\alpha_n}
\end{equation*}
where the $G_i$'s form a nested set $\S$ and for every $i \leq n$ we have $\alpha_i < \rk [\bigvee \S_{< G_i}, G_i]$.
\end{theo}
\red{Sets of indices satisfying the condition above with respect to some nested set $\S$ will be called \textit{$\S$-admissible}.} Using that the above monomials form a linear basis of $\FYa(\L, \G)$ one can see that every Feichtner--Yuzvinsky algebra is in fact of finite dimension and that the part of maximal grading, which is $2(\rk(\L) - 1)$, has dimension one (generated by $x_{\un}^{\rk (\L)-1}$). \\

\textcolor{black}{Essentially using the fact that the change of variable from the generators $x_G$ to the generators $h_G$ has a triangular shape,} one gets the following corollary.
\begin{coro}\label{CoroGrobnerWonderful}
Elements of the form $(\prod_{G' \in \A} (h_G - h_{G'}))h_G^{\rho(G) - \rho(\A)}$ for every antichain $\A \subset \G$ without atoms and every $G \geq \bigvee \A$ form a Gröbner basis of \textcolor{black}{$\I_{\wond}$} for \textcolor{black}{the degree lexicographic ordering on monomials associated to} any linear order on generators refining the reversed order of $\L$. \red{The normal monomials with respect to this Gröbner basis are monomials of the form 
\begin{equation*}
h_{G_1}^{\alpha_1}... h_{G_n}^{\alpha_n}
\end{equation*}
where the $G_i$'s form a nested set $\S$ and $(\alpha_i)_i$ is $\S$-admissible.} 
\end{coro}
\begin{proof}
One can see that the leading terms of those elements are terms of the form $(\prod_{G' \in \A}  h_{G'}) h_G^{\rho(G) - \rho(\bigvee \A)}$ and therefore the normal monomials with respect to those relations are elements of the form $\prod_{G \in \S} h_G^{\alpha_G}$ for any nested set $\S$ and any $\S$-admissible indices $(\alpha_G)_G$, which are in obvious bijection with the normal monomials for the affine Gröbner basis. This proves that those monomials form a linear basis of $\FYa_{\wond}(\L, \G)$ (linearly independent with the right cardinality), which implies that the elements above form a Gröbner basis of $\I_{\wond}$.
\end{proof}
This Gröbner basis will be of use in subsequent sections. In the next and last preliminary subsection we introduce another important combinatorial invariant.
\subsubsection{The Orlik--Solomon algebras}
In this document ``graded commutative'' means with Koszul signs. {\color{black}A \textit{circuit} in a geometric lattice is a set of atoms $C = \{H_i, \, 1 \leq i \leq n \}$ such that $\rho ( \bigvee C)$ is $n-1$ and $\rho(\bigvee X) = | X | $ for all proper subsets $X \subset C$.} 
\begin{madef}[Orlik--Solomon algebra]
Let $\L$ be a geometric lattice. We define the Orlik--Solomon graded commutative algebra $\OSa(\L)$ by
\begin{equation*}
\OSa(\L) = \Lambda [e_H, \, H \textrm{ atom of } \L  ] / \I
\end{equation*}
where $\I$ is the ideal generated by elements of the form $\delta (e_{H_1} \wedge ... \wedge e_{H_n}) $ for any circuit $\{H_1, ... , H_n \}$ and $\delta$ is the unique derivation of degree $-1$ satisfying $\delta(e_H) = 1$. All the generators $e_H$ have degree $1$.
\end{madef}

In the complex realizable case \textcolor{black}{that} algebra is \textcolor{black}{isomorphic to} the cohomology ring of the complement of the hyperplane arrangement (see Orlik--Solomon \cite{OS_1980}).  \\

We denote by $\OSbara(\L)$ the subalgebra of $\OSa(\L)$ generated by elements of the form $e_H - e_{H'}$ for every pair of atoms $H, H'$. In the complex realizable case the algebra $\OSbara(\L)$ is the cohomology ring of the projective complement. We have the following important lemma.
\begin{lemma}[\cite{Yuzvinsky_2001} Section 2.4]\label{lemmaosbar}
For every geometric lattice $\L$ we have the equality
\begin{equation*}
\OSbara(\L) = \ker \delta = \im \delta.
\end{equation*}
\end{lemma}

\section{The Feynman category}\label{secfeycat}
In this section we show that the combinatorial objects introduced in the previous section (geometric lattices, building sets and nested sets) can be bundled up into a Feynman category.

\subsection{A short introduction to Feynman categories}
The notion of Feynman categories was introduced by Kaufmann and Ward in \cite{kaufmann_feynman_2017}. Loosely speaking Feynman categories encode types of operadic structures.
\begin{nota*}
Let $\C$ be a category. We denote by $\C^{\iso}$ the subcategory of $\C$ having the same objects as $\C$ but only its isomorphisms as morphisms. We denote by $\Symcat(\C)$ the free symmetric monoidal category generated by $\C$. For any functor $F: \C\rightarrow \D$ with $\D$ a symmetric monoidal category, there is a unique induced strong monoidal functor $\Symcat(F): \Symcat(\C) \rightarrow \D$ which restricts to $F$ on $\C$. If we are given a diagram of categories $\C \overset{F}{\rightarrow} \D \overset{G}{\leftarrow} \E$, the comma category $(F \downarrow G)$ is the category having for objects triples $(c \in \C, e \in \E, \phi: F(c) \rightarrow G(e))$ and for morphisms suitable commutative diagrams. If the functors $F$ and $G$ are clear from the context we will write instead $(\C \downarrow \E$).
\end{nota*}
\begin{madef}[Feynman category]\label{deffeycat}
A triple $\F = (\nucat, \fcat, \imath)$ is a \emph{Feynman category} if $\nucat$ is a groupoid, $\fcat$ is a symmetric monoidal category and $\imath:\nucat \rightarrow \fcat$ is a functor such that:
\begin{enumerate}
\item The functor $\imath$ induces an equivalence of categories $\Symcat(\imath): \Symcat(\nucat) \overset{~}{\rightarrow} \fcat^{\iso}$.
\item The functor $\imath$ induces an equivalence of categories $\Symcat((\fcat \downarrow \nucat)^{\iso}) \rightarrow (\fcat \downarrow \fcat)^{\iso}$ \r{(the functors used to define those comma categories being the identity of $\fcat$ and $\imath$).}
\item For every object $\star \in \nucat$, the comma category $(\fcat \downarrow \star)$ is essentially small (i.e. is equivalent to a small category).
\end{enumerate}
\end{madef}
A morphism in $\fcat$ which is not an isomorphism will be called a structural morphism. \red{An object of $\nucat$ will be called an arity.} \r{On a first reading the reader is encouraged to replace axioms $(1)$ and $(2)$ by their strong counterpart, that is, imagining that the involved equivalences of categories are in fact isomorphisms, or even just the identity. In that case axioms $(1)$ and $(2)$ can be roughly summarized by saying that the objects of $\fcat$ are words in some alphabet (the objects of $\nucat$). One can ``freely'' concatenate morphisms, say $w_1 \rightarrow w_1'$ and $w_2 \rightarrow w_2'$, to obtain a morphism $w_1w_2 \rightarrow w_1'w_2'$ (this comes from the monoidal structure). Finally, modulo isomorphisms, every morphism in $\fcat$ can be uniquely expressed as the concatenation of morphisms of the form $w \rightarrow w'$ where $w'$ consists of only one letter. In even simpler terms this means that $\fcat$ boils down to: 1) the letters of the alphabet and their isomorphisms, 2) the morphisms $w \rightarrow w'$ where $w'$ has only one letter, and 3) the rules for composing such morphisms. }
\begin{madef}[Operad over a Feynman category]\label{defopfeycat}
Let $\F = (\nucat, \fcat, \imath)$ be a Feynman category and $\C$ a symmetric monoidal category. An \textit{operad over} $\F$ in $\C$ is a strong monoidal functor from $\fcat$ to $\C$, and a \textit{cooperad over} $\F$ in $\C$ is a strong monoidal functor from $\fcat$ to $\C^{op}$. A \textit{module over} $\F$ in $\C$ is a functor from $\nucat$ to $\C$.
\end{madef}
(Co)operads (resp. modules) over $\F$ will also be called $\F$-(co)operads (resp. $\F$-modules).
\begin{ex}\label{exfeycatoperad}
As described in \cite{kaufmann_feynman_2017}, there exists a Feynman category $\Op$ encoding classical operads i.e. such that operads over $\Op$ are classical operads and modules over $\Op$ are $\Sym$-modules. \red{In a few words, the underlying groupoid of $\Op$ is given by the category with objects finite rooted corollas (i.e. trees with one inner vertex) and morphisms the isomorphisms of graphs preserving the root. The objects of $\fcat$ are disjoint unions of such corollas, and the structural morphisms are given by some particular graph operations (roughly the operations which do not create any cycle). For instance, starting with two corollas $c_1, c_2$ with leaves indexed by, say $\{1,2\}$ and $\{3,4,5\}$ respectively, one has the operation of gluing the root of the second corolla with leave $2$ of the first corolla, and then contracting that new edge. This gives a morphism in $\fcat$ from $c_1 \sqcup c_2$ to the corolla $c_3$ having leaves indexed by $\{1,3,4,5\}$. This morphism encodes the operation of ``composing at input $2$'': $$ \mu(x_1, x_2), \lambda(x_3, x_4, x_5) \rightarrow \mu(x_1, \lambda(x_3, x_4, x_5)).$$ }
\end{ex}

\subsection{Construction of the Feynman category}\label{secFeycons}
Let us start by defining the underlying groupoid of our Feynman category.
\begin{madef}
We define $LBS$ to be the groupoid having as objects the built lattices and morphisms
\begin{multline*}
\Mor_{LBS}((\L, \G), (\L', \G')) = \\ \{ f: \L' \xrightarrow{\sim} \L \textrm{ isomorphism of poset satisfying } f(\G') = \G \}.
\end{multline*}
We denote by $\LBSi$ the full subcategory of $LBS$ having as objects the irreducible built lattices.
\end{madef}
The groupoid $\LBSi$ will play the role of $\nucat$ in Definition \ref{deffeycat}. We will add morphisms to $LBS$ in order to get the right category $\fcat$.
\begin{prop}\label{dec}
The category LBS admits a symmetric monoidal structure $\otimes$ given by $$(\L, \G)\otimes(\L', \G') = (\L \times \L', \G\times\{\zero\} \cup \{\zero\}\times\G').$$ Furthermore the inclusion $\imath: \LBSi \rightarrow LBS$ induces an equivalence of categories $$\Symcat(\imath): \Symcat(\LBSi) \overset{~}{\rightarrow} LBS.$$
\end{prop}
\begin{proof}
The fact that a product of two geometric lattices is again a geometric lattice is classical and the proof can be found in \cite{welsh_matroid_1976}. Additionally  $\G\times\{\zero\} \cup \{\zero\}\times\G'$ is indeed a building set of $\L\times \L'$ because for any $(X,X') \in \L\times\L'$ we have the isomorphisms
\begin{align*}
[(\zero, \zero), (X,X')] & \simeq [\zero, X]\times[\zero, X'] \\
&\simeq \prod_{G \in \Fact_{\G}(X)} [\zero,G]\times \prod_{G' \in \Fact_{\G'}(X')}[\zero, G']\\
&\simeq \prod_{G'' \in \Fact_{\G\times\{\zero\} \cup \{\zero\}\times\G'}((X,X'))}[(\zero, \zero), G''].
\end{align*}

Besides, one can see that $\otimes$ is functorial and satisfies the associativity/symmetry axioms of a symmetric monoidal product, the unit being $(\{\zero\}, \emptyset)$.  \\

For the last claim we show that $\Symcat(\imath)$ is essentially surjective and fully faithful. Let $(\L, \G)$ be an object of $LBS$. If we denote by $\{G_i, i \leq n\}$ the factors of $\un$ in $\G$ then we have an isomorphism $[\zero, \un] \simeq \prod_{i}[\zero, G_i]$ and $\G$ is sent to $\G \cap [\zero, G_1] \cup ... \cup \G \cap [\zero, G_n]$ . In other words $(\L, \G)$ is isomorphic to $([\zero, G_1], \G \cap [\zero, G_1])\otimes ... \otimes ([\zero, G_n], \G \cap [\zero, G_n])$ and $\Symcat(\imath)$ is essentially surjective.\\

Finally, let $\bigotimes_{i \leq n}(\L_i, \G_i)$ and $\bigotimes_{j \leq n'}(\L'_j, \G'_j)$ be two elements of $\Symcat(\LBSi)$ (here $\otimes$ denotes the free symmetric monoidal product in $\Symcat(\LBSi)$) and let $\phi$ be an isomorphism in $LBS$ between $\bigotimes_{i \leq n}(\L_i, \G_i)$ and $\bigotimes_{j \leq n'}(\L'_j, \G'_j)$. Such an isomorphism is given by a bijection between the factors of both sides (the isomorphism induces a bijection between the maximal elements of the building set of the domain and the maximal elements of the building set of the target) together with isomorphisms between corresponding summands. This datum is exactly equivalent to an isomorphism in $\Symcat(\LBSi)$ between $\bigotimes_{i \leq n}(\L_i, \G_i)$ and $\bigotimes_{j \leq n'}(\L'_j, \G'_j)$ (via $\imath$), which proves that $\Symcat(\imath)$ is fully faithful.

\end{proof}

We now add structural morphisms to $LBS$ to get our Feynman category. Let $(\L ,\G)$ be an object of $\LBSi$ and $\S=\{G_i, i\leq n\}$ a spanning linearly ordered nested set of $\G$. For any $G \in \S$ we define $\tau_{\S}(G) \coloneqq \bigvee \S_{< G}$ \red{(with $<$ the lattice order)} and we set $$(\L_{\S}, \G_{\S}) \coloneqq \bigotimes_i ([\tau_{\S}(G_i), G_i], \Ind_{[\tau_{\S}(G_i), G_i]}(\G))$$ which is an object of $LBS$. We will view $\S$ as a new formal morphism $(\L_\S, \G_\S) \overset{\S}{\rightarrow} (\L, \G)$ that we  will add by hand to $LBS$. However, to this end one must specify how those new morphisms compose with each other and with the isomorphisms as well. This is the object of the next definition/lemma. \\

From now on, every nested set is assumed to be spanning unless stated otherwise. If $\S$ is a nested set, the intervals $[\tau_{\S}(G), G]$ for $G$ any element of $\S$ will be called the ``local intervals of $\S$''.\\

Let $\S= \{G_i, i \leq n\}$ be a nested set in $(\L, \G)$ and let there be given additional linearly ordered nested sets $\S_i$'s in each irreducible built lattice $([\tau_{\S}(G_i),G_i], \Ind_{[\tau_{\S}(G_i),G_i]}(\G))$. We define {\color{black}
\begin{equation}\label{defcompnested}
\S \circ (\S_i)_i \coloneqq \bigcup_i \,\{\Comp_{\tau_{\S}(G_i)}(K)\, , \, K \in \S_i \}
\end{equation}}
which comes naturally equipped with a linear order (by concatenating the linear orders of the $\S_i$'s). 
{\color{black} 
\begin{rmq}\label{rmkSinComp}
Note that since the $\S_i$'s are assumed to be spanning and we have $\Comp_{\tau_{S}(G_i)}(G_i) = G_i$ for all $i$, we must have $\S\circ(\S_i)_i \supset \S$.
\end{rmq}
\begin{ex}\label{excompnestedmax}
If $\G$ is the maximal building set of $\L$, then $\S$ is simply a chain in $\L$. Moreover, each induced building set on each local interval of $\S$ is the maximal building set of that local interval, which means that the $\S_i's$ are simply chains in each local interval. Finally, as we saw in Example \ref{excompmax}, for the maximal building set the Comp operation is simply the inclusion of the induced building set, which means that $\S\circ(\S_i)_i$ is simply the concatenation of the chains $\S_i, i\leq n$. In this case it is obvious that $\S\circ(\S_i)_i$ is again a nested set. Lemma \ref{lemmacompisnested} will show that this statement is true in general. Throughout this paper, the reader is encouraged to examine how each construction or statement simplifies when considering the maximal building set, in order to distinguish the difficulties arising from the combinatorics of building sets and nested sets from the rest. 
\end{ex}
\begin{computation}\label{computerelations}
Before jumping into technical lemmas, let us compute some compositions of small nested sets. First consider two elements $G_1\neq G_2 \in \G\setminus\{\un\}$ forming a nested antichain. The element $G_1\vee G_2$ belongs to both $\Ind_{[G_1, \un]}(\G)$ and $\Ind_{[G_2, \un]}(\G)$. We have 
\begin{equation*}
    \{G_1, \un\}\circ (\{G_1\}, \{G_1\vee G_2,\un\}) = \{\Comp_{\zero}^{G_1}(G_1), \Comp_{G_1}^{\un}(\un), \Comp_{G_1}^{\un}(G_1\vee G_2)\}. 
\end{equation*}
The first two elements in the right-hand term are obviously $G_1$ and $\un$. For the last element, notice that by Lemma \ref{factors} the factors of $G_1\vee G_2$ are $G_1$ and $G_2$, which means that we must have $\Comp_{G_1}^{\un}(G_1\vee G_2) = G_2$. To conclude we get 
\begin{equation*}
\{G_1, \un\}\circ (\{G_1\}, \{G_1\vee G_2,\un\}) = \{G_1, G_2, \un\} = \{G_2, \un\}\circ (\{G_2\}, \{G_1\vee G_2,\un\})
\end{equation*}
(the last equality by symmetry). On the other hand if we consider two comparable elements $G_1<G_2<\un$ in $\G$, then $G_2$ belongs to $\Ind_{[G_1, \un]}(\G)$, while $G_1$ belongs to $\Ind_{[\zero, G_2]}(\G)$, and we have 
\begin{equation*}
\{G_1,\un\}\circ (\{G_1\}, \{G_2, \un\}) = \{G_1, G_2, \un\} = \{G_2,\un\}\circ (\{G_1, G_2\}, \{\un\}).
\end{equation*}
\end{computation}}
We have the key lemma.
\begin{lemma}\label{lemmacompisnested}
For any built lattice $(\L, \G)$, the set $\S \circ (\S_i)_i$ defined above is a nested set of $(\L,\G)$.
\end{lemma}
\begin{proof}
\red{
By Lemma \ref{rmkinj}, Lemma \ref{lemmaforest} and Remark \ref{rmkSinComp} we have the equality 
\begin{equation}\label{eqpartition}
\S\circ(\S_i) = \S \sqcup \bigsqcup_i \Comp(\S_i \setminus\{G_i\}).
\end{equation}}
Let $\A = \{G_i \, |\, i \in I \}\sqcup \bigsqcup_{j\in J} \A_j$ be a non-empty antichain in $\S \circ (\S_i)_i$ partitioned according to partition \eqref{eqpartition}, with all the $\A_j$'s nonempty and such that the join of $\A$ belongs to $\G$. Our goal is to prove that $\A$ is a singleton. Since $\A$ is an antichain, $I$ and $J$ are disjoint. Let $M$ be the set of maximal elements of $\{G_i, i \in I \} \cup \{G_j, j \in J\}$. {\color{black}
\begin{sublemma}
The set $M$ is a singleton. 
\end{sublemma}}
\begin{proof}
The set $M$ is a nested antichain and so by Lemma \ref{factors} the set $M$ is the set of factors of $\bigvee M$. By definition of $M$ we have $\bigvee \A \leq \bigvee M$ and so if $\bigvee \A$ belongs to $\G$ then $\bigvee \A$ must be below one of the elements of $M$. By definition of $M$ this means that $M$ only contains that element. 
\end{proof}
If this singleton is contained in $\{G_i, i \in I\}$, then by the fact that $\A$ is an antichain we must have $\A = M$. If this singleton is contained in $\{G_j, j \in J\}$, let us denote this singleton $\{G_j\}$. {\color{black}
\begin{sublemma}
The set $\A_j$ is a singleton. 
\end{sublemma}}
\begin{proof}
By definition of $G_j$ we have the equality $$\left(\bigvee \A_j \right)\vee \tau_{S}(G_j) = \left(\bigvee \A_j \right)\vee \tau_{S}(G_j) \in \Ind_{[\tau_{S}(G_j), G_j]}(\G).$$
In addition, we have 
\begin{align*}
    \left(\bigvee \A_j\right) \vee \tau_{S}(G_j) & = \bigvee_{\substack{G \in \S_j \mathrm{ s.t. } \\ \Comp(G) \in \A_j}}(\Comp(G)\vee \tau_{S}(G_j)) \\
    &= \bigvee_{\substack{G \in \S_j \, \mathrm{ s.t. } \\ \Comp(G) \in \A_j}}G.
\end{align*}
By Lemma \ref{lemmacompincr} the set $\{G \in \S_j \, | \,  \Comp(G) \in \A_j\}$ is an antichain in the nested set $\S_j$ and so by the two previous equalities it must be a singleton, which implies the desired result. 
\end{proof}
Let us denote $\A_j = \{\Comp(K)\}$ with $K$ some element in~$\S_j$. We have
\begin{equation*}
\Comp(K) \vee \tau_{\S}(G_j) = \left(\bigvee \A \right)\vee  \tau_{\S}(G_j) = K.
\end{equation*}
Since $\bigvee \A$ belongs to $\G$, by the very definition of $\Comp(K)$ we must have $\bigvee \A \leq \Comp(K)$, or in other words every element of $\A$ must be below $\Comp(K)\in \A$. The set $\A$ being an antichain, $\A$ must be equal to the singleton~$\{\Comp(K)\}$.
\end{proof}

\begin{lemma}\label{caniso}
There is an isomorphism of built lattices $$\bigotimes_i (\L_{\S_i}, \G_{\S_i}) \xrightarrow{\Phi} (\L_{\S\circ (\S_i)_i}, \G_{\S \circ (\S_i)_i}).$$
\end{lemma}
\begin{proof}
\red{The irreducible summands of the built lattice on the left-hand side are indexed by $\bigsqcup \S_i$, which is also the case of the irreducible summands appearing on the right-hand side, by Lemma \ref{rmkinj} and Lemma \ref{lemmaforest}. We need to show that the irreducible built lattice indexed by some element $K$ in some $\S_i$ on the left-hand side is isomorphic to the irreducible built lattice indexed by the same element in the right-hand side. If we denote by $K_1,\ldots, K_p$ the maximal elements of $(\S_i)_{<K}$, the former irreducible built lattice is 
\begin{equation*}
([\tau_{\S_i}(K), K], \Ind_{[\tau_{\S_i}(K), K]}(\Ind_{[\tau_{\S}(G_i), G_i]}(\G)))
\end{equation*}
while the latter is
\begin{equation*}
\left(\left[\left(\bigvee_i \Comp(K_i)\right) \vee \left(\bigvee _{G_j< \Comp(K)} G_j\right), \Comp(K)\right], \Ind(\G)\right).
\end{equation*}}
There is an isomorphism between the two, given by taking the join with $\tau_{\S}(G_i)$. The fact that this is indeed an isomorphism of poset\red{s} comes from the building set isomorphism
\begin{equation*}
[\zero, K] \simeq [\zero, \Comp(K)]\times\prod_{\substack{G_j \in \max \S_{<G_i} \\ G_j \nleq \Comp(K)}}[\zero, G_j].
\end{equation*}
The fact that it sends building set to building set comes from Lemma \ref{lemmadoubleinduction}.
\end{proof}
Finally, we have to prove that the operation $\circ$ on nested sets is ``associative''.
\begin{lemma}
Let $\S$ be a nested set in $(\L, \G)$. Let $(\S_i)_i$ be nested sets in every local interval of $\S$ and for every $i$ let $(\S^i_j)_j$ be nested sets in every local interval of $\S_i$. The nested sets $\S \circ (\S_i \circ (\S^{i}_{j})_j)_j$ and $(\S \circ (\S_i)_i)\circ (\S^{i}_j)_{i,j}$ are equal (the last composition being performed via the isomorphism $\Phi$ of Lemma \ref{caniso}).
\end{lemma}
\begin{proof}
This is a statement about the $\Comp$ operation. When necessary we put the lattice in which we are doing the $\Comp$ operation in superscript. We denote $\S = \{G_i\, , \,  i\leq n\}$ and choose some $i_0 \leq n$. We then denote $\S_{i_0} = \{K_j\,,\, j \leq m \}$ and choose some $j_0 \leq m$.  We also define $I_0 \coloneqq \{i \,| \, G_i \in \max \S_{<G_{i_0}}\}$ and $J_0 \coloneqq \{j \,| \,K_j \in \max ((\S_{i_0})_{<K_{j_0}}) \}$. \\

We partition $I_0$ into the following subsets:
\begin{align*}
I_0^{\ext} &\coloneqq \{i \in I_0 \, | \, G_i \nleq \Comp(K_{j_0})\}, \\
I_0^{\int} &\coloneqq \{i \in I_0 \, | \, G_i \leq \Comp(K_{j_0}) \textrm{ and } \forall j \in J_0,\, G_i \nleq \Comp(K_j) \}, \\
I_0^{j} &\coloneqq \{i \in I_0 \, | \, G_i \leq \Comp(K_j) \}.
\end{align*}
It is a partition because the set $\{K_j, \, j \in J_0 \}$ is nested. The map $\Phi$ is taking the join with~$\bigvee_{i \in I_0^{\ext}} G_i$.\\

The lemma amounts to showing that for any $L$ in $\Ind_{[\tau_{\S_{i_0}}(K_{j_0}), K_{j_0}]}(\G)$ we have the equality
\begin{multline*}
\Comp^{\L}_{\tau_{\S}(G_{i_0})}(\Comp^{[\tau_{\S}(G_{i_0}), G_{i_0}]}_{\tau_{\S_{i_0}}(K_{j_0})}(L)) = \Comp_{\bigvee \{\Comp(K_j) ,\,  j \in J_0\} \vee \bigvee \{G_i,\, i \in I_0^{\int}\}}^{\L}(\Phi^{-1} (L)).
\end{multline*}
We have
\begin{multline*}
\Comp^{\L}_{\tau_{\S}(G_{i_0})}(\Comp^{[\tau_{\S}(G_{i_0}), G_{i_0}]}_{\tau_{\S_{i_0}}(K_{j_0})}(L)) \vee \bigvee_{j \in J_0} \Comp(K_j) \vee \bigvee_{i \in I_0^{\int}} G_i\vee \bigvee_{i \in I_0^{\ext}}G_i = \\
\Comp^{[\tau_{\S}(G_{i_0}), G_{i_0}]}_{\tau_{\S_{i_0}}(K_{j_0})}(L) \vee \bigvee_{j \in J_0} \Comp(K_j) \vee \bigvee_{i \in I_0^{\int}}G_i \vee \bigvee_{i \in I_0^{\ext}}G_i = L.
\end{multline*}
Applying $\Phi^{-1}$ to both sides we get
\begin{equation*}
\Comp^{\L}_{\tau_{\S}(G_{i_0})}(\Comp^{[\tau_{\S}(G_{i_0}), G_{i_0}]}_{\tau_{\S_{i_0}}(K_{j_0})}(L)) \vee \bigvee_{j \in J_0} \Comp(K_j) \vee \bigvee_{i \in I_0^{\int}} G_i  = \Phi^{-1}(L).
\end{equation*}
We need to prove that $\Comp^{\L}_{\tau_{\S}(G_{i_0})}(\Comp^{[\tau_{\S}(G_{i_0}), G_{i_0}]}_{\tau_{\S_{i_0}}(K_{j_0})}(L))$ is the biggest element in $\G$ which satisfies this equation. \\

Let $G' \in \G$ be such that we have $$G' \vee \bigvee_{j \in J_0} \Comp(K_j) \vee \bigvee_{i \in I_0^{\int}} G_i  = \Phi^{-1}(L).$$
Applying $\Phi$ on both sides we get
\begin{equation}\label{eqdefG'}
G' \vee \bigvee_{j \in J_0} \Comp(K_j) \vee \bigvee_{i \in I_0^{\int}} G_i \vee \bigvee_{i \in I_0^{\ext}} G_i = L.
\end{equation}
The element $G'\vee \bigvee_{i \in I}G_i$ is below $L$ and belongs to $\Ind_{[\tau_{\S}(G_{i_0}), G_{i_0}]}(\G)$ so it is below one of the factors of $L$ in $\Ind_{[\tau_{\S}(G_{i_0}), G_{i_0}]}(\G)$. Those factors are $\Comp^{[\tau_{\S}(G_{i_0}), G_{i_0}]}_{\tau_{\S_{i_0}}(K_{j_0})}(L)$ or $K_j$ for some $j$ in $J$. By equation \eqref{eqdefG'}, $G'\vee \bigvee_{i \in I}G_i$ cannot be below any $K_j$ so we have
\begin{equation*}
G'\vee \bigvee_{i \in I}G_i \leq \Comp^{[\tau_{\S}(G_{i_0}), G_{i_0}]}_{\tau_{\S_{i_0}}(K_{j_0})}(L),
\end{equation*}
which implies
\begin{equation*}
G' \leq \Comp^{\L}_{\tau_{\S}(G_{i_0})}(\Comp^{[\tau_{\S}(G_{i_0}), G_{i_0}]}_{\tau_{\S_{i_0}}(K_{j_0})}(L)).
\end{equation*}
\end{proof}
We are now in position to make the following definition.
\begin{madef}
$\LBScat$ is the monoidal category defined as follow\red{s}.
\begin{itemize}
\item The objects of $\LBScat$ are built lattices.
\item The morphisms of $\LBScat$ are generated (via composition and tensoring) by
    \begin{enumerate}
    \item Structural morphisms
\begin{equation*}
(\L_{\S}, \G_{\S}) \xrightarrow{\S} (\L, \G)
\end{equation*}
for every totally ordered nested set in some irreducible built lattice $(\L, \G)$. The composition of those morphisms is given by \eqref{defcompnested}.
    \item Isomorphisms between built lattices
\begin{equation*}
(\L, \G) \xrightarrow{\sim} (\L', \G')
\end{equation*}
for each isomorphism of poset\red{s} $f:\L' \rightarrow \L$ such that $f(\G') = \G$,
    \end{enumerate}
quotiented by relations
\begin{equation}\label{eqrelisoall}
(\L_{f(\S)}, \G_{f(\S)}) \overset{f(\S)}{\rightarrow}(\L, \G) \overset{f}{\rightarrow}(\L', \G')  \sim (\L_{f(\S)}, \G_{f(\S)}) \overset{\otimes f}{\rightarrow} (\L'_{\S}, \G'_{\S}) \overset{\S}{\rightarrow} (\L', \G')
\end{equation}
for any isomorphism of irreducible built lattice\red{s} $f: (\L, \G) \rightarrow (\L', \G')$, and relations
\begin{equation}\label{eqrelpermutation}
(\L_\S, \G_\S) \overset{\sigma}{\rightarrow} (\L_{\S_{\sigma}}, \G_{\S_{\sigma}}) \overset{\S_{\sigma}}{\rightarrow} (\L, \G) \sim (\L_\S, \G_\S) \overset{\S}{\rightarrow} (\L, \G)
\end{equation}
for every linearly ordered nested set $\S$ and $\sigma$ some permutation of the summands of $(\L_\S, \G_\S)$ (with $\S_{\sigma}$ the nested set equipped with the new linear order given by $\sigma$). Finally we also impose the relation
\begin{equation*}
\{\un \} \sim \Id_{(\L, \G)}.
\end{equation*}
\item The monoidal structure is the same as the one on $LBS$, which in addition acts on nested sets (which are now considered as morphisms) by disjoint union.
\end{itemize}

\end{madef}\begin{prop}
The triple $\LBS = (\LBSi, \LBScat, \imath)$ with $\imath$ the obvious inclusion is a Feynman category.
\end{prop}
{\color{black}
\begin{proof}
Proposition \ref{dec} shows that $\LBS$ satisfies Axiom (1) of Definition \ref{deffeycat}. For Axiom (2) notice that by relations \eqref{eqrelisoall} and \eqref{eqrelpermutation} every morphism in $\LBScat$, i.e. every object of $(\LBScat \downarrow \LBScat)$, can be factored as a composition $\S\circ \sigma \circ f$ with $\S$ a spanning nested set, $\sigma$ a permutation and $f$ an isomorphism of built lattices. The nested set $\S$ can be uniquely decomposed as a tensor of nested sets landing in irreducible built lattices. This decomposition gives an inverse of the functor $\Symcat((\LBScat \downarrow \LBSi)^{\iso}) \rightarrow (\LBScat \downarrow \LBScat)^{\iso}$. Axiom (3) is immediate. 
\end{proof}}

We conclude this subsection by proving a general lemma on the composition of nested sets which will be important later on.
\begin{lemma}\label{lemmarecons}
Let $\S$ be some spanning nested set in some irreducible built lattice and $\S' \subset \S$ a subset containing $\un$. For any $G'$ in $\S'$ we put $$\S'_{G'} \coloneqq (\S \vee \tau_{\S'}(G')) \cap (\tau_{\S'}(G'), G'].$$ For any $G'$ in $\S'$, $\S'_{G'}$ is a nested set in $([\tau_{\S'}(G'), G'], \Ind(\G))$ and we have the equality between nested sets:
\begin{equation}\label{eqdecrec}
\S = \S'\circ(\S'_{G'})_{G' \in \S'}.
\end{equation}
\end{lemma}
\begin{proof}
Let $G'$ be any element of $\S'$ and $G'_1, ..., G'_n$ the maximal elements of $\S'_{<G'}$. Let $G$ be some element in $\S$ which is below $G'$ and not below any of the $G'_i$'s. We denote $K \coloneqq \tau_{\S'}(G')\vee G$. By the fact that $\S$ is nested, $G$ is the maximal element of $\G$ satisfying the equality
\begin{equation*}
G\vee \tau_{\S'}(G') = K.
\end{equation*}
This proves the equality 
\begin{equation*}
G = \Comp_{\tau_{\S'}(G')}(\tau_{\S'}(G')\vee G),
\end{equation*}
which implies equality \eqref{eqdecrec}. \\

For the nestedness, assume we have $G_1, ... , G_k$ some elements of $\S$ such that $\tau_{\S'}(G')\vee G_1, ..., \tau_{\S'}(G')\vee G_k$ are elements of $(\tau_{\S'}(G'), G']$ and such that $\bigvee_i \tau_{\S'}(G')\vee G_i$ belongs to $\Ind_{[\tau_{\S'}(G'), \un]}(\G)$. The factors of $\bigvee_i \tau_{\S'}(G')\vee G_i$ in $\G$ are some of the $(G'_i)'s$ and the element $\Comp_{\tau_{\S'}(G')}(\bigvee_i \tau_{\S'}(G')\vee G_i)$. Since $\Comp_{\tau_{\S'}(G')}$ is increasing (it has a left inverse given by taking the join with $\tau_{\S'}(G')$), by the first part of this proof we have $$G_i = \Comp_{\tau_{\S'}(G')}(\tau_{\S'}(G')\vee G_i) \leq \Comp_{\tau_{\S'}(G')}(\bigvee_i \tau_{\S'}(G')\vee G_i)$$ for all $i \leq k$ and by the building set isomorphism we must have $$\bigvee_i G_i = \Comp_{\tau_{\S'}(G')}(\bigvee_i \tau_{\S'}(G')\vee G_i).$$ By nestedness of $\S$ the $G_i$'s do not form an antichain and therefore the $\tau_{S'}(G')\vee G_i$'s do not either. This proves that $(\tau_{\S'}(G')\vee \S) \cap (\tau_{\S'}(G'), \un]$ is a nested set.
\end{proof}

\subsection{Presentation of $\LBS$}\label{presentation}
In this section we give a presentation of the category $\LBScat$, that is, a set of structural morphisms that together with isomorphisms generate every other morphism in $\LBS$, via composition and tensoring. In general, having a presentation of some category $\fcat$, part of a Feynman category $\F = (\nucat, \fcat, \imath)$, is useful from a practical point of view because it makes defining operads over $\F$ a lot easier (one just needs to specify the operad on the generators and check that those morphisms satisfy the right relations).

\begin{ex}
In the case of classical operads, the Feynman category $\Op$ has a natural set of generating structural morphisms which can be represented by rooted trees with two inner vertices. This translates into the fact that operads can be defined by their partial compositions (compositions of two operations).
\end{ex}

\begin{prop}
Every morphism in $\LBScat$ can be obtained as a composition and tensoring of morphisms of the form $\{G, \un\}$ and isomorphisms.
\end{prop}
\begin{proof}
Iterate Lemma \ref{lemmarecons} with $\S'$ of the form $\{G, \un\}$.
\end{proof}
In the rest of the article, we will denote $[G] \coloneqq \{G, \un\}$ (viewed as a morphism in $\LBScat$).
\begin{prop}\label{presfey}
Relations between compositions of generators \red{$[G]$} are generated by the relations
\begin{equation}\label{eqrelchain}
\red{[G_1]} \circ (\red{[G_2]} \otimes \Id ) = \red{[G_2]} \circ (\Id \otimes \red{[G_1]})
\end{equation}
for every pair $G_1 < G_2$, relations
\begin{equation}\label{eqrelanti}
\red{[G_1]} \circ (\red{[G_1 \vee G_2]} \otimes \Id) = \red{[G_2]} \circ (\red{[G_1\vee G_2]} \otimes \Id) \circ \sigma^{2,3}
\end{equation}
for every pair $G_1$, $G_2$ of non comparable elements forming a nested set (with $\sigma^{2,3}$ the transposition swapping the two last summands) and relations
\begin{equation}\label{eqreliso}
f \circ \red{[f(G)]} = \red{[G]} \circ (f_{[G, \un]} \otimes f_{[\zero, G]})
\end{equation}
for every isomorphism $f: (\L, \G) \xrightarrow{\sim} (\L', \G')$ in $LBS$ and every element $G \in \G \setminus \{\un \}$.
\end{prop}
\begin{proof}
\red{The first two equations were proved by Computation \ref{computerelations} and the last equation is immediate.} In order to prove that those relations generate all the others we start with the following lemma.
\begin{lemma}\label{lemmadecmor}
Let $(\L, \G)$ be a built lattice and $\vartriangleleft$ a linear order on $\L$. Any morphism in $\LBS$ $$F: \bigotimes_i (\L_i , \G_i) \rightarrow (\L, \G)$$ can be uniquely written as a composition
\begin{equation}\label{eqdecmor}
\bigotimes_i (\L_i, \G_i) \xrightarrow{\bigotimes f_i} \bigotimes_i (\L'_i, \G'_i) \xrightarrow{\sigma} \bigotimes_i (\L'_{\sigma(i)}, \G'_{\sigma(i)}) \xrightarrow{\S} (\L, \G),
\end{equation}
where
\begin{itemize}
\item The $f_i$'s are isomorphisms in $\LBSi$.
\item $\sigma$ is a permutation of the summands.
\item $\S$ is a spanning nested set of $(\L, \G)$ with total order given by restriction of $\vartriangleleft$.
\end{itemize}
\end{lemma}
\begin{proof}
By iteration of relations \eqref{eqrelisoall} and \eqref{eqrelpermutation} one can see that every morphism can be written in the form \eqref{eqdecmor}. For the \red{uniqueness} we define an invariant $\iota(F): \bigsqcup_i (\L_i\setminus\{\zero\}) \rightarrow \L \setminus \{\zero\}$ for every morphism $F: \bigotimes_i (\L_i,\G_i) \rightarrow (\L, \G)$ in $\LBScat$ by setting
\begin{itemize}
\item For any nested set $\S$ in some built lattice $(\L, \G)$:
\begin{equation*}
\iota(\S): \bigsqcup_G ([\tau_{\S}(G), G]\setminus \{\tau_{\S}(G)\}) \hookrightarrow \L \setminus \{\zero\}
\end{equation*}
is the obvious inclusion.
\item For any isomorphism $(\L',\G') \xrightarrow{f} (\L, \G)$ in $\LBScat$, $\iota(f)$ is equal to $f^{-1}$.
\end{itemize}
and then extending to all morphisms by composition. One can see that $\iota$ preserves the relations \eqref{eqrelisoall} and $\eqref{eqrelpermutation}$ and therefore it passes to the quotient and gives a well-defined invariant for every morphism in $\LBS$. \\

Now if $F: \bigotimes_i (\L_i, \G_i) \rightarrow (\L, \G)$ can be written as a composition
\begin{equation*}
\bigotimes_i (\L_i, \G_i) \xrightarrow{\bigotimes f_i} (\L'_i, \G'_i) \xrightarrow{\sigma} (\L'_{\sigma(i)}, \G'_{\sigma(i)}) \xrightarrow{\S} (\L, \G),
\end{equation*}
then we see that $\S$ and the $f_i's$ can be extracted from $\iota(F)$ ($\S$ is just the image by $\iota(F)$ of $\bigsqcup_i \{\un_{\L_i} \}$ in $\L \setminus \{\zero\}$ and the $f_i's$ are just restrictions of $\iota(F)$ to the suitable subsets). Then, $\sigma$ is the only permutation that permutes the summands in the right order when $\S$ is given the order $\vartriangleleft$. This proves the \red{uniqueness}.
\end{proof}

Assume now that we have two sequences of morphisms in $\LBScat$
$$\varphi = A  \rightarrow X_1 \rightarrow \cdots \rightarrow X_n \rightarrow (\L, \G) \, \, \, \textrm{and} \, \, \, \psi = A \rightarrow X'_1 \rightarrow \cdots \rightarrow X'_{n'} \rightarrow (\L, \G),$$
such that every morphism in $\phi$ or $\psi$ is either a generator of the form $\Id \otimes \cdots \otimes \Id \otimes \red{[G]} \otimes \Id \otimes \cdots \otimes \Id$ or an isomorphism, and $(\L, \G)$ is some irreducible built lattice. We want to prove that if the composition of the morphisms of $\varphi$ is equal to the composition of the morphisms of $\psi$ then there is a chain of equivalences of the form \eqref{eqrelchain}, \eqref{eqrelanti} or \eqref{eqreliso} (possibly tensored and composed with other morphisms) between $\varphi$ and $\psi$. \\

By iterated use of relation \eqref{eqreliso} we can assume that the only morphisms of $\varphi$ and $\psi$ which are isomorphisms are at the beginning of $\varphi$ and $\psi$ respectively. Then, by Lemma \ref{lemmadecmor} we can assume that both $\varphi$ and $\psi$ only contain generators which are not isomorphisms. Let us denote by $\S$ the nested set obtained by composing the generators of $\varphi$, or equivalently $\psi$. We also denote by $\S_i$ (resp. $\S'_i$) the nested set obtained by composing the last $i$ morphisms of $\varphi$ (resp. $\psi$). By construction of the composition of nested sets (see Section \ref{secFeycons}) we have
\begin{equation*}
\S_1 \subsetneq ... \subsetneq \S_n = \S
\end{equation*}
and
\begin{equation*}
\S'_1 \subsetneq ... \subsetneq \S'_{n'} = \S,
\end{equation*}
and the \textcolor{black}{cardinality} of the nested sets increases exactly by one at each step. Let us denote $\S_1 = \{G_{\varphi}\}$. By the equations above there exist some $j \leq n'$ such that we have $$\S'_{j} \setminus \S'_{j-1} = \{G_{\varphi}\}.$$

One can find a chain of relations of the form \eqref{eqrelchain} and \eqref{eqrelanti} between $\psi$ and some $\psi'$ such that the first morphism of $\psi'$ is $\{G_{\varphi}\}$ (applying relations \eqref{eqrelchain} and \eqref{eqrelanti} allows one to swap successively the morphism corresponding to $G_{\varphi}$ in $\psi$ with the morphism right after, until it reaches the end). \\

This means that we can assume that the last morphism of $\psi$ is $\{G_{\varphi}\}$. We denote by $\S^{< G_{\varphi}}$ (resp. $\S^{>G_{\varphi}}$) the nested set obtained by composing the morphisms of $\varphi$ which correspond to generators in $[\zero, G]$ (resp. $[G, \un]$), and we denote similarly $\S'^{<G_{\psi}}$, $\S'^{>G_{\psi}}$ the same constructions but with $\psi$. We have
\begin{equation*}
\{G_{\varphi}\}\circ (\S^{< G_{\varphi}}, \S^{> G_{\varphi}}) = \S = \{G_{\varphi}\} \circ(\S'^{< G_{\varphi}}, \S'^{> G_{\varphi}})
\end{equation*}
which implies that we have
\begin{align*}
\S^{< G_{\varphi}} &= \S'^{< G_{\varphi}} \\
\S^{> G_{\varphi}} &= \S'^{> G_{\varphi}}
\end{align*}
and we conclude by induction.
\end{proof}
{\color{black}This leads to the following concrete description of an $\LBS$-operad.
\begin{coro}\label{corodescriptionoperad}
The datum of an $\LBS$-operad in a symmetric monoidal category $\C$ is equivalent to the datum of a collection of objects in $\C$ $$\O(\L, \G)$$ indexed by irreducible built lattices, together with isomorphisms in $\C$
$$\O(\L, \G) \xrightarrow{\sim} \O(\L', \G')$$ for every isomorphism of built lattices $f:(\L,\G) \xrightarrow{\sim} (\L', \G')$, together with morphisms in $\C$ 
\begin{equation*}
    \O([\zero, G], \Ind(\G))\otimes \O([G,\un], \Ind(\G)) \rightarrow \O(\L, \G)
\end{equation*}
for every element $G\neq \un$ in the building set of some irreducible built lattice $(\L, \G)$, satisfying relations \eqref{eqrelchain}, \eqref{eqrelanti},  and \eqref{eqreliso}. 
\end{coro}}
\begin{ex}\label{exshuffle}
Let $\O$ be an $\LBS$-operad. Let $m,n \geq 2$, $i\leq n$ be some integers and $\sigma:\{i+1,\ldots, m+n-1\} \rightarrow \{i+1,\ldots, m+n-1\}$ an $(m-1, n-i)$-shuffle, i.e. a permutation which is increasing on both $\{i+1,\ldots,i+m-1\}$ and $\{i+m,\ldots, m+n-1\}$. Consider the partition $G \in \Pi_{\{1,\ldots,m+n-1\}}$ with only one non-trivial equivalence class given by $\{i, \sigma(i+1),\ldots,\sigma(i+m-1)\}$. We have an isomorphism of built lattices 
\begin{equation*}
    (\Pi_{\{1,\ldots, m\}}, \G_{\min}) \xrightarrow[\sim]{f} ([\zero, G], \Ind(\G_{\min}))
\end{equation*}
given by identifying the partitions in $[\zero, G]$ with the partitions of $\{i, \sigma(i+1),\ldots,\sigma(i+m-1)\}$ and then identifying this latter set with $\{1,\ldots, m\}$ via~$\sigma$. We also have an isomorphism of built lattices 
\begin{equation*}
    (\Pi_{\{1,\ldots, n\}}, \G_{\min}) \xrightarrow[\sim]{g} ([G,\un], \Ind(\G_{\min}))
\end{equation*}
given by identifying the partitions in $[G,\un]$ with the partitions of the set of equivalence classes of $G$, and identifying this latter set with $\{1,\ldots, n\}$ via~$\sigma$. Precomposing $\O([G])$ with $\O(f)\otimes \O(g)$ gives a morphism 
\begin{equation*}
    \O(\Pi_{\{1,\ldots, m\}}, \G_{\min})\otimes \O(\Pi_{\{1,\ldots, n\}}, \G_{\min}) \xrightarrow{\circ_{i, \sigma}} \O(\Pi_{\{1,\ldots, m+n-1\}}, \G_{\min}).
\end{equation*}
In this case, relations \eqref{eqrelchain} and \eqref{eqrelanti} exactly say that the objects $\O(\Pi_{\{1,\ldots, n\}}, \G_{\min})$ ($n\geq 2)$ and the above morphisms form a well-known algebraic structure called a shuffle operad (see \cite{DK_2010}). To summarize, every $\LBS$-operad contains a shuffle operad obtained by restricting to partition lattices with minimal building sets, and forgetting the symmetries. An alternative way to see this would be to recall from Example \ref{exnestedsetshuffletree} that nested sets in $(\Pi_{\{1,\ldots, n\}}, \G_{\min})$ can be identified with shuffle trees with leaves indexed by $\{1,\ldots, n\}$, and that via this identification the composition of nested sets is the grafting of shuffle trees. 
\end{ex}
\begin{rmq}
One might wonder what kind of object arises if we keep the symmetries in the previous example. For $\n \geq 3$, the automorphism group of $(\Pi_{\{1,\ldots,n\}}, \G_{\min})$ is $\Sym_n$, whereas for $n=2$ the automorphism group of $(\Pi_{\{1,2\}}, \G_{\min})$ is trivial (and not $\Sym_2$). Based on this, one might expect that we obtain a symmetric operad with trivial $\Sym_2$ action, but this is not the case. Indeed, although there exists a functor from symmetric operads with trivial $\Sym_2$ action to $\mathfrak{P}\mathsf{art}$-operads, where $\mathfrak{P}\mathsf{art}$ is the full subcategory of $\LBS$ generated by partition lattices with minimal building sets, there is no natural functor in the reverse direction. Thus, one might describe a $\mathfrak{P}\mathsf{art}$-operad as an intermediate structure between a symmetric operad (with trivial $\Sym_2$ action) and a shuffle operad (both assumed trivial in arity $0$ and $1$). 
\end{rmq}
\subsection{$\LBS$ is a graded Feynman category}\label{secgradedfeycat}
An important feature of $\LBS$ is that morphisms of $\LBScat$ can be graded in the following sense.
\begin{madef}[\red{Kaufmann and Ward \cite{kaufmann_feynman_2017}}]
A degree function on a Feynman category $(\nucat, \fcat, \imath)$ is a map $$\deg: \Mor(\fcat) \rightarrow \N$$ such that
\begin{itemize}
\item $\deg (\phi \circ \psi) = \deg(\phi) + \deg(\psi)$
\item $\deg (\phi \otimes \psi) = \deg(\phi) + \deg(\psi)$
\item Morphisms of degree $0$ and $1$ generate $\Mor(\fcat)$ by compositions and tensor products.
\end{itemize}
Furthermore the degree function is said to be proper if the degree $0$ morphisms are exactly the isomorphisms. A graded Feynman category is a Feynman category equipped with a degree function.
\end{madef}

\begin{ex}
The Feynman category $\Op$ encoding classical operads has a proper grading given by defining the degree of a tree $t$ as the number of inner vertices of $t$ minus one. \end{ex}

For $\LBS$ we define a proper degree function by setting
\begin{align*}
\deg (f) &= 0 \textrm{ for all isomorphisms } f, \\
\deg (\red{[G]}) &= 1 \textrm{ for all } G \textrm{ different from } \un.
\end{align*}
One can see that the relations introduced in Proposition \ref{presfey} are homogeneous with respect to this grading and therefore we can extend this grading to every morphism in $\LBS$, which makes $\LBS$ a properly graded Feynman category.

\section{(Co)operads over $\LBS$}\label{secoperads}

In this section we show that the algebraic invariants introduced in Section \ref{secprelim} (Feichtner--Yuzvinsky algebras, Orlik--Solomon algebras) can be bundled up to form various (co)operads over $\LBS$ (see Definition \ref{defopfeycat}).

\subsection{The Feichtner--Yuzvinsky cooperad}
In this subsection we define and study an $\LBS$-cooperadic structure on the family of Feichtner--Yuzvinsky rings.
\subsubsection{Definition of the cooperad}
\begin{lemma}
The map $(\L, \G) \rightarrow \FYa(\L, \G)$ can be upgraded to a (strong) monoidal functor from $LBS$ to $\gCR^{op}$ where $\gCR$ is the symmetric monoidal category of graded commutative rings.
\end{lemma}
\begin{proof}
Let $(\L, \G)$ and $(\L', \G')$ be two built lattices and let $(\L'', \G'') = (\L, \G) \otimes (\L', \G')$. We have an isomorphism of algebras
\begin{equation*}
\begin{array}{rclcc}
\FYa(\L, \G) &\otimes &\FYa(\L', \G') &\xrightarrow{\sim} &\FYa(\L'', \G'') \\
x_G &\otimes &1 & \rightarrow &x_G \\
1 &\otimes &x_{G'} &\rightarrow &x_{G'} \\
\end{array}
\end{equation*}
with inverse
\begin{align*}
\FYa(\L'', \G'') & \rightarrow \FYa(\L, \G) \otimes \FYa(\L', \G') \\
x_G & \rightarrow
\left\{
\begin{array}{ll}
x_G \otimes 1 & \textrm{ if } G \in \G \\
1\otimes x_G & \textrm{ otherwise.}
\end{array}
\right.
\end{align*}
If $\phi$ is an isomorphism $(\L', \G') \overset{\sim}{\rightarrow} (\L, \G)$ in $LBS$, it induces the isomorphism of algebras
\begin{align*}
\FYa(\L, \G) &\rightarrow \FYa(\L', \G') \\
x_G &\rightarrow x_{\phi(G)}
\end{align*}
which is compatible with composition.
\end{proof}
Next we upgrade $\FYa$ into a monoidal functor from $\LBScat$ to $\gCR^{op}$. Thanks to the presentation of $\LBScat$ given in Subsection \ref{presentation} we only need to specify the action of $\FYa$ on the generators $[G]$ and then check that those morphisms satisfy the right relations. For every $G \in \G\setminus \{\un\}$ we set (using this time the wonderful presentation)
\begin{equation*}
\begin{array}{rccl}
\FYa(\red{[G]}):&\FYa(\L , \G) &\longrightarrow& \FYa([G,\un], \Ind(\G)) \otimes \FYa([\zero, G], \Ind(\G)) \\
&h_{G'} & \longrightarrow &
\left\{
\begin{array}{ll}
1 \otimes h_{G'} &\textrm{ if } G' \leq G \\
h_{G\vee G'} \otimes 1  & \textrm{ otherwise.}
\end{array}
\right.
\end{array}
\end{equation*}
In the realizable case this morphism of algebra is induced by the inclusion of the stratum $\overline{Y}_{\{G, \un\}}$ in the wonderful compactification $\overline{Y}_{\L, \G}$. \red{See \cite{Braden_2022} Section 2.6 for a toric-theoretic perspective on those morphisms.}
{\color{black}
\begin{lemma}
The morphism $\FYa([G])$ is well-defined.
\end{lemma}}
\begin{proof}
We must prove that $\FYa([G])$ sends $\I_{\wond}$ to zero. If $H$ is an atom of $\L$ which is below $G$ then $h_H$ is sent to $h_H\otimes 1$ which is zero in the target algebra. On the contrary if $H$ is an atom of $\L$ which is not below $G$ then this means by sub-modularity of $\L$ that $H\vee G$ is an atom of $[G, \un]$ and thus $h_H$ is sent to zero again. Notice that this is the first time we have actually used the geometricity of our lattices. Now let $X = \{G_1, ..., G_n \}$ be some antichain in $\G$ having join $G' \in \G$. Let us assume that the first $k$ $G_i$'s are the elements of $X$ below $G$. If $G' \leq G$ then $\prod_i(h_G - h_{G_i})$ is sent to $$1 \otimes \prod_i(h_G - h_{G_i}),$$ which is zero in the target algebra. Otherwise if $G' > G$ then $\prod_{i >k}(h_{G'} - h_{G_i})$  is sent to $$\prod_{i > k}(h_{G \vee G'} - h_{G \vee G_i} ) \otimes 1,$$ which is also zero in the target algebra.
\end{proof}

\begin{prop}
The morphisms $\FYa(\red{[G]})$ satisfy the relations given in Proposition \ref{presfey}.
\end{prop}
\begin{proof}
Let $G_1 < G_2$ be two comparable elements in $\G\setminus\{\un\}$. We have to check the equality of algebra morphisms
\begin{equation*}
(\FYa(\red{[G_2]})\otimes \Id) \circ \FYa(\red{[G_1]}) = (\Id \otimes \FYa(\red{[G_1]}) \circ \FYa(\red{[G_2]}),
\end{equation*}
and it is enough to check it on generators. Let $G$ be any element in $\G$. If $G \leq G_1$ one can check that both morphisms send $h_{G}$ to $1\otimes 1 \otimes h_{G}$. If $G \leq G_2$ and $G \nleq G_1$ one can see that both morphisms send $h_{G}$ to $1 \otimes h_{G_1\vee G} \otimes 1$. Lastly, if $G \nleq G_2$ one can check that both morphisms send $h_{G}$ to $h_{G_2 \vee G} \otimes 1 \otimes 1$.\\

Let $G_1, G_2$ be two non-comparable elements of $\G \setminus \{\un\}$ forming a nested set. We need to check the equality of algebra morphisms
\begin{equation*}
(\FYa(\red{[G_1 \vee G_2]}) \otimes \Id) \circ \FYa(\red{[G_1]}) = \sigma^{2,3} \circ (\FYa(\red{[G_1 \vee G_2]}) \otimes \Id) \circ \FYa(\red{[G_2]}).
\end{equation*}
It amounts again to a simple verification on generators with a small \red{trichotomy}. Let $G$ be any element in $\G$. If $G \leq G_1$ this means that $G$ cannot be below $G_2$ and we see that both morphisms send $h_G$ to $1\otimes 1 \otimes h_G$, if $G \leq G_2$ by a similar argument both morphisms send $h_G$ to $1\otimes h_G \otimes 1$. Finally if $G$ is neither below $G_1$ nor below $G_2$ then both morphisms send $h_G$ to $h_{G\vee G_1 \vee G_2}\otimes 1 \otimes 1$. \\

We can give an explicit formula for general morphisms $\FYa(\S)$. If $G$ is any element in $\G$ let $G'$ be the unique minimal element of $\S_{>G}$ (Lemma \ref{lemmaforest}). We then have
\begin{equation*}
\FYa(\S)(h_G) = 1^{\otimes} \otimes h_{\tau_{\S}(G') \vee G } \otimes 1^{\otimes},
\end{equation*}
where $1^{\otimes}$ means that we put a $1$ in every interval which is not $[\tau_{\S}(\red{G'}), \red{G'}]$. \\

Finally, we need to check that the morphisms $\FYa(\red{[G]})$ satisfy relation \eqref{eqreliso}. Let $(\L, \G)$ and $(\L', \G')$ be two built lattices and $f: (\L',\G') \xrightarrow{\sim} (\L, \G)$ an isomorphism in $\LBS$, i.e. an isomorphism of poset\red{s} $f:\L \xrightarrow{\sim}\L'$ such that $f(\G)$ is equal to~$\G'$. Let $G$ be an element in $\G\setminus\{\un\}$. We have to check the equality between algebra morphisms
\begin{equation*}
\FYa(\red{[f(G)]}) \circ \FYa(f) = (\FYa(f_{|[G, \un]})\otimes \FYa(f_{|[\zero,G]})) \circ \FYa(\red{[G]}).
\end{equation*}
Let $h_{K}$ be some generator in $\FYa(\L, \G)$. If $K\leq G$ one can check that both morphisms send $h_K$ to $1\otimes h_{f(K)}$. Otherwise if $K \nleq G$ one can check that both morphisms send $h_K$ to $h_{f(G)\vee f(K)} \otimes 1$.
\end{proof}

In the sequel we will write $\FY$ when referring to the $\LBS$-cooperad and not just the algebras. We also write $\FY^{\vee}$ for the (linear) dual operad (apply the duality functor to all objects and morphisms). This is an $\LBS$-operad in the category of graded coalgebras.

\subsubsection{A quadratic presentation for $\FY^{\vee}$} \label{secfreeop}
In this section we exhibit a quadratic presentation for the operad $\FY^{\vee}$. Let us first quickly recall what \textcolor{black}{that} means in the context of operads over general Feynman categories.\textcolor{black}{This} is all part of the theory developed by Kaufmann and Ward in \cite{kaufmann_feynman_2017}. \red{We refer to \cite{Ward_2020} Appendix A for more details.} \\

Let $\F = (\nucat, \fcat, i)$ be a Feynman category and $M$ an $\F$-module in some monoidal category $\C$ (see Definition \ref{defopfeycat}). If $\C$ is cocomplete there exists a ``free'' $\F$-operad in $\C$ generated by $M$ denoted by $\F(M)$. The $\F$-operad $\F(M)$ satisfies the universal property that for any morphism of $\F$-module between $M$ and some $\F$-operad $\P$ (viewed as an $\F$-module), there exists a unique morphism of $\F$-operad between $\F(M)$ and $\P$ which extends the first morphism. Concretely, $\F(M)$ is given by the left Kan extension of $M: \Symcat(\nucat) \rightarrow \C$ along $\Symcat(\imath)$. The left Kan extension universal property is exactly the freeness universal property. If furthermore $\F$ is assumed to be a graded Feynman category and $\C$ is a category of modules for instance, then free operads are naturally graded (i.e. components in all arity are graded and structural morphisms preserve \textcolor{black}{that} grading). \\

In addition, if we are given an $\F$-operad $\P$ and $M$ a sub $\F$-module of $\P$ then one can define the ideal $\langle M \rangle$ generated by $M$ in $\P$ by considering all possible elements in $\P$ which can be obtained as the composition (along some structural morphism in $\F$) of an element of $M$ with elements of $\P$. Finally, we can define the quotient of an operad by an ideal in an obvious way (just take the quotient in each arity and notice that the morphisms pass to the quotient). With all those notions at hand we can define a quadratic operad over a graded Feynman category to be the quotient of a free operad by an ideal generated by degree $1$ elements. \red{In most examples, degree $1$ morphisms in Feynman categories have a domain with two irreducible monoidal summands (as in $\LBS$ for instance), and therefore degree $1$ elements in free operads can be expressed as sums of structural products of two generators, hence the terminology ``quadratic''.} Since we have proved that $\LBS$ is a graded Feynman category in Section \ref{secgradedfeycat}, the above vocabulary applies to $\LBS$-operads. \\

Let $\Gen$ be the $\LBS$-module in the category of graded $\Z$-modules defined by
\begin{equation*}
\Gen(\L, \G) = \Z[\red{-}2(\rk (\L) - 1)]
\end{equation*}
and
\begin{equation*}
\Gen(f) = \Id
\end{equation*}
for every irreducible built lattice $(\L, \G)$ and \red{every} isomorphism of built lattice\red{s} $f$. If $(\L, \G)$ is an irreducible built lattice we denote by $\Psi_{(\L, \G)}$ the canonical generator of~$\Gen(\L, \G)$. We also denote for all spanning nested sets $\S$ of $(\L, \G)$:
\begin{equation*}
\Psi_{\S}\coloneqq  \LBS(\Gen)(\S)((\Psi_{([\tau_{\S}(G), G], \Ind(\G))})_{G \in \S}),
\end{equation*}
which is an element of $\LBS(\Gen)(\L, \G)$. \red{For instance, for every $G\in \G\setminus\{\un\}$ we have 
\begin{multline*}
    \Psi_{\{G,\un\}} = \Psi_{([\zero, G], \Ind(\G))}\otimes \Psi_{([G,\un], \Ind(\G))} \\ \in \Gen([\zero, G], \Ind(\G))\otimes\Gen([G,\un], \Ind(\G))\subset \LBS(\Gen)(\L, \G).
\end{multline*}}
We are now able to state the main result of this section.

\begin{prop}\label{theopresopFY}
The $\LBS$-operad $\FY^{\vee}$ is isomorphic to the quotient of $\LBS(\Gen)$ by the ideal generated by the elements

\begin{equation}\label{eqrelFY}
\sum_{\un > G\geq H_1} \Psi_{\{G, \un\}} - \sum_{\un > G \geq H_2} \Psi_{\{ G, \un \}},
\end{equation}
for all atoms $H_1, H_2$ in some irreducible built lattice $(\L, \G)$.
\end{prop}
{\color{black}\begin{rmq}\label{rmkhypercom}
Notice that when restricted to partition lattices with minimal building sets, relations \eqref{eqrelFY} are exactly the relations of $\mathsf{Hypercom}$ discovered in \cite{Getzler_1994}. In other words the shuffle operad contained in $\FY^{\vee}$ (see Example \ref{exshuffle}) is the shuffle operad associated to $\mathsf{Hypercom}$.  
\end{rmq}}
\begin{proof}
We have a map of $\LBS$-modules $\Gen \xrightarrow{\pi} \FY^{\vee}$ sending $\Psi_{(\L, \G)}$ to the linear form which is zero in all degrees except degree $2(\rk(\L) -1)$ where it takes value $1$ on $x_{\un}^{\rk(\L)-1}$ (which linearly generates $\FYa^{2(\rk(\L) - 1)}(\L, \G)$). \\

\textcolor{black}{That} map is a natural transformation since for any isomorphism $f: (\L, \G) \rightarrow (\L', \G')$ between built lattices, $f$ preserves the top element and therefore $\FY(f)$ preserves $h_{\un}^{\rk(\L) -1} = x_{\un}^{\rk(\L)- 1}$ which implies that $\FY(f)^{\vee}$ sends $\pi(\Psi_{(\L, \G)})$ to $\pi(\Psi_{(\L', \G')})$.\\

\textcolor{black}{That} map extends to an $\LBS$-operadic map $\LBS(\Gen) \xrightarrow{\hat{\pi}} \FY^{\vee}$ (by \red{the} universal property of free operads). Our goal is to prove that this map passes to the quotient by the elements \eqref{eqrelFY} and that the induced morphism is an isomorphism. The proof splits into three steps.\\

\textbf{Step 1}: The map $\hat{\pi}$ is surjective. \\

Going back to explicit formulas for general left Kan extensions in cocomplete categories yields
\begin{equation*}
\LBS(\Gen)(\L, \G) = \bigoplus_{\otimes_i (\L_i, \G_i) \rightarrow (\L, \G)} \bigotimes \Gen((\L_i, \G_i)) / \sim\,\,\, ,
\end{equation*}
where the sum is over all possible maps $\otimes_i (\L_i, \G_i) \rightarrow (\L, \G)$ in $\LBS$ and the equivalence relation $\sim$ identifies components corresponding to equivalent maps (two maps that can be obtained from one another by precomposing with isomorphisms). \\

More explicitly if we have two maps $\otimes_i (\L_i , \G_i) \xrightarrow{\psi} (\L, \G)$ and $\otimes_i(\L'_i, \G'_i) \xrightarrow{\phi} (\L, \G)$ such that there exist isomorphisms $f_i: (\L_i, \G_i) \xrightarrow{\sim} (\L'_i, \G'_i)$ satisfying
\begin{equation*}
\psi = \phi \circ (\otimes_i f_i),
\end{equation*}
then we have
\begin{equation*}
\otimes_i \alpha_i \sim \otimes_i \Gen(f_i)(\alpha_i)
\end{equation*}
for every element $\otimes \alpha_i$ in $\otimes_i \Gen((\L_i, \G_i))$. Likewise, if we have an equality of the form
\begin{equation*}
\psi = \phi \circ \sigma
\end{equation*}
with sigma some permutation of the $(\L_i, \G_i)$'s then we have
\begin{equation*}
\otimes_i \alpha_i \sim \otimes_i \alpha_{\sigma(i)}.
\end{equation*}
Finally, replacing $\Gen(\L, \G)$ by its definition we get (for every irreducible built lattice~$(\L, \G)$):
\begin{equation*}
\LBS(\Gen)((\L, \G)) = \Z\langle\{ \otimes_i (\L_i, \G_i) \rightarrow (\L, \G)\}/ \sim \rangle
\end{equation*}
(generators are equivalence classes of maps in $\LBS$ having target $(\L, \G)$, with the equivalence relation being the precomposition with isomorphisms).
With this identification the map $\hat{\pi}$ is given by
\begin{align*}
\LBS(\Gen)(\L, \G) &\longrightarrow \FY^{\vee}(\L, \G)\\
[\mu: \otimes_i (\L_i, \G_i) \rightarrow (\L, \G)] & \longrightarrow (\alpha \rightarrow \otimes_i \pi(\Psi_{(\L_i, \G_i)})(\FY(\mu)(\alpha))).
\end{align*}
Let us fix an irreducible built lattice $(\L, \G)$ and some linear order extending the order on $\L$. Amongst all maps of the form $\otimes_i (\L_i, \G_i) \rightarrow (\L, \G)$ we have the maps given by linearly ordered nested sets
\begin{equation*}
(\L_{\S}, \G_{\S}) \xrightarrow{\S} (\L, \G)
\end{equation*}
whose linear order is given by our chosen global linear order on $\L$. It is enough to prove the surjectivity of $\hat{\pi}$ restricted to equivalence classes of those morphims. Passing to the dual we must prove the injectivity of the map
\begin{equation}\label{eqPDop}
\begin{array}{clc}
\FY(\L, \G) &\longrightarrow &\Z\langle\{\textrm{spanning nested sets of } (\L, \G) \}\rangle \\
\alpha & \longrightarrow &(\pi(\Psi_{\S})(\FY(\S)(\alpha)))_{\S}
\end{array}
\end{equation}
where $\pi(\Psi_{\S})$ denotes the tensor product of maps $\otimes_{G \in \S} \pi(\Psi_{([\tau_{\S}(G), G], \Ind(\G))})$.\\

Let $\alpha$ be an element of $\FY(\L, \G)$ which is sent to zero by the above map, meaning that for every nested set $\S$ in $(\L, \G)$ we have $\pi(\Psi_{\S})(\FY(\S)(\alpha)) = 0$. We can assume that $\alpha$ is homogeneous and by Corollary \ref{CoroGrobnerWonderful} we can write it uniquely as a sum of normal monomials

\begin{equation*}
\alpha = \sum_{\substack{\S \textrm{ nested} \\ \mu \, \S-\textrm{admissible}}} \lambda_{S, \mu} h_{\S}^{\mu}
\end{equation*}
where $h_{\S}^{\mu}$ denotes the monomial $\prod_{G\in \S}h_G^{\mu(G)}$. We have to prove that all the $\lambda_{\S, \mu}$'s are zero. \red{If such is not the case, let $G_0$ be some element of $\G$ such that there exists some nested set $\S_0$ and some $\S_0$-admissible index $\mu_0$ such that $\lambda_{\S_0, \mu_0} \neq 0$, and which is minimal for that property.} \\

For any irreducible built lattice $(\L, \G)$ and any $k < \rk \, \L $, one can construct a nested set $\S(\L, \G, k)$ having only local intervals with rank $1$ except the top interval having rank $k$, by picking any maximal chain $\zero  < X_1 < \ldots < X_n < \un$ in $\L$ and putting
\begin{equation*}
\S(\L, \G, k) = \red{[X_1]} \circ \red{[X_2]} \circ ... \circ \red{[X_{n-k+1}]} .
\end{equation*}
The formula makes sense because each $X_i$ is an atom in $[X_{i-1}, \un]$ and must therefore belong to $\Ind_{[X_{i-1}, \un]}(\G)$. Notice that for any $\S$ and $k\geq 1$ we have
\begin{equation*}
\hat{\pi}(\Psi_{\S(\L, \G, k)})(h_{\un}^{k-1}) = 1.
\end{equation*}
Let us denote by $\S_{G_0}^{\mu_0}$ the nested set $\S([\zero, G], \Ind(\G), \mu_0(G_0) +1)$ in $([\zero, G], \Ind(\G))$. For any nested set $\S'$ in $([G_0, \un], \Ind(\G))$ and any nested set $\S$ in $\G$ with admissible index $\mu$ we have
\begin{equation*}
\hat{\pi}(\Psi_{[G_0] \circ(\S', \S_{G_0}^{\mu_0})}) (h_{\S}^{\mu}) = 0
\end{equation*}
if $\S$ does not contain some element with strictly positive index and which is below $G_0$. By minimality of $G_0$ \textcolor{black}{this} means that we have
\begin{align*}
\hat{\pi}(\Psi_{[G] \circ(\S', \S_{G_0}^{\mu_0})}) (\alpha) &= \sum_{\substack{ \S \textrm{ nested} \\ \mu \, \S-\textrm{admissible} \\ G_0 \in \S \\ \mu(G_0) = \mu_0(G_0)}} \lambda_{\S, \mu} \hat{\pi}(\Psi_{[G] \circ(\S', \S_{G_0}^{\mu_0})})(h_{\S}^{\mu}) \\
&= \sum_{\substack{ \S \textrm{ nested} \\ \mu \, \S-\textrm{admissible} \\ G_0 \in \S \\ \mu(G_0) = \mu_0(G_0)}} \lambda_{\S, \mu} \hat{\pi}(\Psi_{\S'})(h_{G_0 \vee \S}^{\mu}/h_{G_0}^{\mu(G_0)}) \\
&= 0.
\end{align*}
By Lemma \ref{lemmarecons} the monomials $h_{G_0 \vee \S}^{\mu}/h_{G_0}^{\mu(G_0)}$ are normal monomials in the irreducible built lattice $([G, \un], \Ind(\G))$. By induction we get that $\lambda_{\S_0, \mu_0}$ is equal to zero which is a contradiction. \\

\textbf{Step 2}: The map $\hat{\pi}$ sends the elements \eqref{eqrelFY} to zero.\\

We must check that the linear forms
\begin{equation*}
\phi_{H} = \hat{\pi}\left(\sum_{\un > G\geq H} \LBS(\Gen)(\red{[G]})\left(\Psi_{([G, \un], \Ind(\G))}, \Psi_{([\zero, G],\Ind(\G))}\right) \right)
\end{equation*}
indexed by atoms $H$ are all equal. \red{To this end, we will prove that those linear forms send every normal monomial of degree $2(\rk(\L) - 2)$ to $1$.} Fortunately, normal monomials of degree $2(\rk(\L) - 2)$ are rather simple.
\begin{lemma}
For any irreducible built lattice $(\L, \G)$, the only degree $2(\rk(\L) - 2)$ normal monomials are the monomials $h_G^{\rk(G) - 1}h_{\un}^{\rk(\un) - \rk(G) - 1}$ with $G$ any element of $\G$ (when $G$ is an atom we get the monomial $h_{\un}^{\rk(\L) - 2}$ ).
\end{lemma}
\begin{proof}
The proof hinges on the following result.
\begin{lemma}
Let $\S$ be a spanning nested set in some irreducible built lattice $(\L, \G)$. We have the equality
\begin{equation*}
\sum_{G\in \S}\rk([\tau_{\S}(G), G]) = \rk(\L).
\end{equation*}
\end{lemma}
\begin{proof}
The proof goes by induction on the rank of $\L$. Let $G_0$ be a minimal element in $\S$. $G_0\vee (\S\setminus\{G_0\})$ is a nested set in $([G_0, \un], \Ind(\G))$ so by induction hypothesis the sum of the rank of its intervals is the rank of $[G_0,\un]$ but by nestedness of $\S$ taking the join with $G_0$ establishes a bijection between intervals of $\S$ which are not the interval $[\zero, G_0]$ and intervals of $G_0 \vee (\S \setminus \{G_0\})$, and \textcolor{black}{that} bijection preserves the rank of the intervals. Consequently, by induction we have
\begin{align*}
\sum_{G\in \S}\rk([\tau_{\S}(G), G]) &= \rk([\zero, G_0]) + \sum_{G \in \S\setminus \{G_0\}} \rk([\tau_{\S}(G), G]) \\
&= \rk([\zero, G_0]) + \rk([G_0, \un]) \\
&= \rk (\L).
\end{align*}
\end{proof}
Now if $\S$ is any spanning nested set then the degree of any normal monomial of the form $h_{\S}^{\mu}$ for some $\S$-admissible index $\mu$ is at most $$2\sum_{G \in \S} (\rk([\tau_{\S}(G), G]) - 1)$$ which is equal by the previous lemma to $2(\rk(\L) - |\S|)$. This means that to have degree $2(\rk(\L) - 2)$ the cardinality of $\S$ must be at most two (counting $\un$) which proves the result.
\end{proof}
Let $h_{G}^{\rk(G) - 1} h_{\un}^{\rk(\L) - \rk(G) - 1}$ be any degree $2(\rk(\L) - 2)$ normal monomial and $H$ some atom. Going back to the definition of $\phi_H$ we get
\begin{multline*}
\phi_{H}(h_{G}^{\rk(G) - 1} h_{\un}^{\rk(\L) - \rk(G) - 1}) = \\
\sum_{G'\geq H} (\pi(\Psi_{[G', \un]})\otimes \pi(\Psi_{[\zero, G']}))\left(\FY(\red{[G]})(h_{G}^{\rk(G) - 1} h_{\un}^{\rk(\L) - \rk(G) - 1})\right).
\end{multline*}
If $H \leq G$ then the only term which is not zero in \textcolor{black}{that} sum is the term $G' = G$ which gives
\begin{align*}
\phi_{H}(h_{G}^{\rk(G) - 1} h_{\un}^{\rk(\un) - \rk(G) - 1}) &= (\pi(\Psi_{[G, \un]})\otimes \pi(\Psi_{[\zero, G]}))\left(\FY(\red{[G]})(h_{G}^{\rk(G) - 1} h_{\un}^{\rk(\un) - \rk(G) - 1})\right) \\
&= (\pi(\Psi_{[G, \un]})\otimes \pi(\Psi_{[\zero, G]}))\left(h_{\un}^{\rk(\un) - \rk(G) - 1}\otimes h_{G}^{\rk(G) - 1}\right) \\
&= 1.
\end{align*}

If $H \nleq G$ the only term which is not zero in \textcolor{black}{that} sum is the term $G' = H$ which gives
\begin{align*}
\phi_{H}(h_{G}^{\rk(G) - 1} h_{\un}^{\rk(\un) - \rk(G) - 1}) &= (\pi(\Psi_{[H, \un]})\otimes \pi(\Psi_{[\zero, H]}))\left(\FY([H])(h_{G}^{\rk(G) - 1} h_{\un}^{\rk(\un) - \rk(G) - 1}\right) \\
&= (\pi(\Psi_{[H, \un]})\otimes \pi(\Psi_{[\zero, H]}))\left( h_{G\vee H}^{\rk(G) - 1}h_{\un}^{\rk{\L}- \rk(G) - 1}  \otimes 1 \right) \numberthis \label{eqHnleqG}\\
&= \pi(\Psi_{[H, \un]})\left( h_{G\vee H}^{\rk(G) - 1}h_{\un}^{\rk{\L}- \rk(G) - 1}\right).
\end{align*}
The monomial in the last equation is not a normal monomial in the irreducible built lattice $([H, \un], \Ind(\G))$ and therefore we must rewrite it. By geometricity of $\L$ there exist atoms $H_1, ..., H_{rk(\L) - \rk(G) - 1}$ such that $\un$ is equal to the join of $G\vee H $ with those atoms. This means that we have the relation
\begin{equation*}
(h_{\un}- h_{G\vee H})(h_{\un} - h_{H_1})...(h_{\un} - h_{H_{\rk(\L) - \rk(G) - 1}}) = 0
\end{equation*}
in the algebra $\FYa([H, \un], \Ind(\G))$, and replacing the $h_{H_i}$'s by zero we get
\begin{equation*}
h_{\un}^{\rk(\L) - \rk(G)}     =   h_{\un}^{\rk(\L)- \rk(G) - 1}h_{G\vee H}
\end{equation*}
which implies
\begin{equation*}
h_{G\vee H}^{\rk(G) - 1}h_{\un}^{\rk(\L)- \rk(G) - 1} = h_{\un}^{\rk(\L) - 2}
\end{equation*}
\textcolor{black}{That} equality together with equation \eqref{eqHnleqG} leads directly to
\begin{equation*}
\phi_{H}(h_{G}^{\rk(G) - 1} h_{\un}^{\rk(\un) - \rk(G) - 1}) = 1
\end{equation*}
as in the case $H \leq G$. \\

\textbf{Step 3}: The kernel of $\hat{\pi}$ is generated by the relations \eqref{eqrelFY}. \\

We postpone the proof of this last step to Subsection \ref{secgrobnerFY}, where it will be an application of our theory of Gröbner bases for $\LBS$-operads.
\end{proof}
\red{\begin{rmq}
Recent articles (\cite{Huh_2018}, \cite{Pagaria_2021}) have demonstrated that even beyond the realizable case the rings $\FYa_{\mathbb{Q}}(\L, \G) \coloneqq \FYa(\L, \G)\otimes \Q$ still satisfy a (purely algebraic) ``Kähler package'' (a set of properties typically enjoyed by rational cohomology rings of projective smooth complex varieties). The simplest (by far) of those properties is Poincaré duality. This states that there exists an integer $n$ such that for all $k >n$ we have $\FYa^k_{\Q}(\L, \G) \simeq \{0\}$, along with an isomorphism $\FYa^n_{\Q}(\L, \G)\xrightarrow[\sim]{\deg}\Q$, such that for all $k \leq n$ the pairing 
\begin{equation*}
    \FYa^k_{\Q}(\L, \G) \otimes \FYa^{n-k}_{\Q}(\L, \G)\xrightarrow{\mathrm{mult}}\FYa^n_{\Q}(\L, \G) \xrightarrow[\sim]{\deg} \Q
\end{equation*}
is non-degenerate. In our case the number $n$ is $\rk(\L) - 1$. Let us now remark that for any spanning nested set $\S$ the linear form $$(\hat{\pi}(\Psi_{\S})\circ \FY(\S))\otimes \Q: \FYa_{\Q}(\L, \G) \rightarrow \mathbb{Q}$$ is zero on $\FYa^{k}_{\Q}(\L, \G)$ for $k$ different from $2(\rk(\L) - \# \S)$, and equal to the multiplication by $x_{\S \setminus \{\un\}}$ (via the identification $\FYa_{\Q}^{2(\rk(\L) -1)}(\L, \G) \overset{\deg}{\simeq} \mathbb{\Q}$) for $k = 2(\rk(\L) - \# \S)$. This means that the injectivity of the map \eqref{eqPDop} can be restated as follows: for all elements $\alpha \in  \FY^{2k}_{\Q}(\L, \G)$, if $\alpha x_\S = 0 $ for all nested sets $\S \subset \G\setminus \{\un\}$ of \textcolor{black}{cardinality} $\rk(\L) - 1 - k$, then $\alpha = 0$. In other words, we have reproved that $\FYa_{\Q}(\L, \G)$ satisfies Poincaré duality. 
\end{rmq}}
\subsubsection{The Feichtner--Yuzvinsky operad}
One can define another operad out of the Feichtner--Yuzvinsky algebras as follow\red{s}. Let $(\L, \G)$ be any irreducible built lattice. We put
\begin{equation*}
\FY^{\PD}(\L, \G) \coloneqq \FYa(\L, \G).
\end{equation*}
For any isomorphism of built lattice\red{s} $f:(\L, \G) \xrightarrow{\sim} (\L', \G')$ we define
\begin{equation*}
\begin{array}{rccc}
\FY^{\PD}(f):& \FYa(\L, \G) &\longrightarrow &\FYa(\L', \G') \\
&x_G& \longmapsto & x_{f^{-1}(G)}
\end{array}
\end{equation*}
and for any $G \in \G\setminus\{\un\}$ we put
\begin{equation*}
\begin{array}{rccl}
\FY^{\PD}(\red{[G]}):& \FY^{\PD}([G, \un], \Ind(\G)) \otimes \FY^{\PD}([\zero, G], \Ind(\G)) &\rightarrow &\FY^{\PD}(\L, \G) \\
& \prod_i x_{G'_i}^{\alpha_i} \otimes \prod_j x_{G_j}^{\alpha_j} &\mapsto& x_G \prod_{i}x_{\Comp_G(G'_i)}^{\alpha_i} \prod_j x_{G_j}^{\alpha_j}
\end{array}
\end{equation*}
if all the $G_j$'s are different from $G$. 
{\color{black} \begin{lemma}
The morphism above is well-defined. 
\end{lemma}}
\begin{proof}
For any antichain $\{G_i\}$ below $G$ and such that $\bigvee_i G_i$ belongs to $\G$ we have
\begin{equation*}
\FY^{\PD}([G])(1 \otimes \prod_i x_{G_i}) = x_G \prod_i x_{G_i} = 0.
\end{equation*}
For any antichain $\{G\vee G_i\}$ in $\Ind_{[G, \un]}(\G)$ having join $G'$ in $\Ind_{[G, \un]}(\G)$ we have
\begin{equation*}
\FY^{\PD}([G])(\prod_{i}x_{G \vee G_i} \otimes 1) = x_G \prod_{i}x_{\Comp_{G}(G\vee G_i)} = 0
\end{equation*}
because either $G'$ belongs to $\G$ and the elements $\{G, \Comp_{G}(G\vee G_i) \}$ have join $G'$ in $\G$, either $G'$ does not belong to $\G$ and the elements $\Comp_{G}(G\vee G_i)$ have join $\Comp_{G}(G')$. \\

For any atom $H$ below $G$ we have
\begin{align*}
\FY^{\PD}([G])(1 \otimes \sum_{\substack{G' \geq H \\ G' < G}}x_{G'} ) &= \sum_{\substack{G' \geq H \\ G' < G}} x_G x_{G'} \\
&= x_G(h_{H} - \sum_{G' \geq G} x_{G'} - \sum_{G',G \textrm{ incomparable }}x_{G'}) \\
&= x_G(-\sum_{G' \geq G}x_{G'} - \sum_{G', G \textrm{ incomparable}}x_{G'}),
\end{align*}
which does not depend on $H$. Finally, if $G\vee H$ is an atom in $[G, \un]$, either $G\vee H$ belongs to $\G$ in which case we have
\begin{equation*}
\FY^{\PD}(\red{[G]})(h_H\otimes \un) = x_{G}h_{G\vee H} = 0,
\end{equation*}
either $G \vee H$ does not belong to $\G$ in which case we have
\begin{equation*}
\FY^{\PD}(\red{[G]})(h_H\otimes \un) = x_{G}h_{H} = 0.
\end{equation*}
\end{proof}
{\color{black}\begin{lemma}
The morphisms above satisfy the relations in Proposition \ref{presfey}. 
\end{lemma}}
\begin{proof}
Let $(\L, \G)$ be some irreducible built lattice. If $G_1$ and $G_2$ are two non comparable elements of $\G\setminus \{\un\}$ forming a nested set, one can check that both $$\FY^{\PD}(\red{[G_1]})\circ(\FY^{\PD}(\red{[G_1\vee G_2]})\otimes \Id)$$ and $$\FY^{\PD}(\red{[G_2]})\circ(\FY^{\PD}(\red{[G_1\vee G_2 ]})\otimes \Id)\circ \sigma_{2,3}$$ send $\prod_{i}x_{G_i}\otimes\prod_{j}x_{G'_j}\otimes \prod_{k}x_{G''_k}$ to $x_{G_1}x_{G_2}\prod_{i}x_{\Comp_{G_1\vee G_2}(G_i)}\prod_{j}x_{G'_j}\prod_{k}x_{G''_k}$. \\

If $G_1 < G_2<\un$ one can check that both $$\FY^{\PD}(\red{[G_1]})\circ(\FY^{\PD}(\red{[G_2]})\otimes \Id)$$ and $$\FY^{\PD}(\red{[G_2]})\circ(\Id \otimes \FY^{\PD}(\red{[G_1]}))$$ send $\prod_{i}x_{G_i}\otimes\prod_{j}x_{G'_j}\otimes \prod_{k}x_{G''_k}$ to $x_{G_1}x_{G_2}\prod_{i}x_{\Comp_{G_2}(G_i)} \prod_j x_{\Comp_{G_1}(G'_j)}\prod_k x_{G''_k}$. \\

Finally, if $f: (\L, \G) \xrightarrow{\sim} (\L', \G')$ is an isomorphism of irreducible built lattice\red{s} and $G'$ is some element in $\G'$, one can check that both $$\FY^{\PD}(f) \circ \FY^{\PD}(\red{[f(G')]})$$ and $$\FY^{\PD}(\red{[G']})\circ(\FY^{\PD}(f_{|[G', \un]})\otimes \FY^{\PD}(f_{[\zero, G']}))$$ send $\prod_i x_{G_i}\otimes \prod_j x_{G'_j}$ to $x_{G'}\prod_i x_{f^{-1}(\Comp_{f(G')}(G_i))} \prod_j x_{f^{-1}(G'_j)}$. 
\end{proof}

The operads $\FY^{\vee}$ and $\FY^{\PD}$ are related via Poincaré duality. If we denote by $\PD$ the isomorphism \red{$\FYa^{\bullet}(\L, \G) \xrightarrow{\sim} \FYa^{\rk \L - 1 - \bullet}(\L, \G)^{\vee}$} given by Poincaré duality we have the equality between morphisms
\begin{equation*}
\FY^{\PD}(\red{[G]}) = \PD^{-1} \circ \FY^{\vee}(\red{[G]})\circ (\PD \otimes \PD).
\end{equation*}

\subsection{The affine Orlik--Solomon cooperad}\label{sectionafforlik}
In this section we introduce an $\LBS$-cooperadic structure on the Orlik--Solomon algebras. \red{When restricted to the collection of partition lattices with minimal building sets (see Example \ref{exshuffle}) this gives the shuffle operad associated to the Gerstenhaber operad (see \cite{Cohen_1973}).}
\subsubsection{Definition}
Let $(\L, \G)$ be an irreducible built lattice. For all $G \in \G\setminus \{ \un \}$ we define a morphism of algebras $\OS(\red{[G]})~:~\OSa(\L)~\rightarrow~\OSa([G, \un])\otimes \OSa([\zero,G])$ by
\begin{equation*}
 \OS(\red{[G]})(e_H) =
 \left\{
 \begin{array}{ll}
 1 \otimes e_H & \textrm{ if } H \leq G \\
 e_{G\vee H} \otimes 1 & \textrm{ otherwise,}
 \end{array}
 \right.
\end{equation*}
for every generator $e_H$.
\begin{lemma}
The morphism $\OS(\red{[G]})$ is well defined.
\end{lemma}
\begin{proof}
We have to check that $\OS(\red{[G]})$ vanishes on elements of the form $\delta \textrm{circuit}$. Let $C = \{H_i\} \sqcup \{H'_j\}$ be a circuit with the $H_i$'s below $G$ and the $H'_j$'s not below $G$. We have
\begin{align*}
\OS(\red{[G]})(\delta (\prod C)) &= \OS(\red{[G]})(\delta(\prod e_{H_i} \wedge \prod e_{H'_j})) \\
&= \OS(\red{[G]})(\delta(\prod e_{H_i}) \wedge \prod e_{H'_j} \pm \prod e_{H_j} \wedge \delta(\prod e_{H'_j})) \numberthis \label{eqdelta} \\
&= \delta(\prod e_{H_i}) \otimes \prod e_{G\vee H'_j} \pm \prod e_{H_j} \otimes \delta(\prod e_{G\vee H'_j}).
\end{align*}
If the $H_i$'s form a set of dependent atoms of $\L$ then we have $\delta(\prod e_{H_i}) = \prod e_{H_i} = 0$ (these identities holding in both $\OSa(\L)$ and $\OSa([\zero,G])$). If on the contrary the $H_i$'s form a set of independent atoms, by the fact that $C$ is dependent there exists an atom $H'_{j_0}$ in $C$ which is below $\bigvee H_i \vee \bigvee_{j<j_0} H'_j$. By taking the join with $G$ in \textcolor{black}{that} relation we obtain $G\vee H'_{j_0} \leq G\vee \bigvee_{j<j_0}H'_j$ which implies that the $G\vee H'_j$'s form a set of dependent atoms of $[G,\un]$ which shows that we have $\delta(\prod e_{G\vee H'_j}) = \prod e_{G\vee H'_j} = 0 $ in the algebra $\OSa([G,\un])$.
 \end{proof}
\begin{lemma}
$\OS$ extends to an $\LBS$-cooperad of graded commutative algebras.
\end{lemma}
\begin{proof}
On objects we set $\OS(\L, \G) = \OSa(\L)$ for every built lattice $(\L, \G)$. For structural morphisms we use the morphisms $\OS(\red{[G]})$ introduced above on each generator of $\LBS$ (one-element nested sets). On isomorphisms the action is
\begin{align*}
\OS(\L, \G) &\rightarrow \OS(\L', \G') \\
e_H & \rightarrow e_{f(H)}
\end{align*}
for any isomorphism $f: (\L', \G') \xrightarrow{\sim} (\L, \G)$ in $\LBS$ (an isomorphism of posets must preserve the atoms). \\

We have to check that those morphisms satisfy the relations given in Proposition \ref{presfey}. Let $G_1 \leq G_2 < \un $ be two comparable elements in the building set $\G$ of a lattice $\L$. Let us prove the equality of morphisms
\begin{equation*}
(\Id \otimes \OS(\red{[G_1]})) \circ \OS(\red{[G_2]}) = (\OS(\red{[G_2]}) \otimes \Id ) \circ \OS(\red{[G_1]}).
\end{equation*}
Since we are dealing with morphisms of algebras it is enough to prove the equality on generators which amounts to a simple verification. For any atom $H$ in $\L$, if $H \leq G_1$ then both morphisms send $e_H$ to $1\otimes 1 \otimes e_H$. If $H \leq G_2$ and $H \nleq G_1$ then both morphisms send $e_H$ to $1\otimes e_{G_1\vee H}\otimes 1$. Lastly, if $H \nleq G_2$ then both morphisms send $e_H$ to $e_{G_2 \vee H} \otimes 1 \otimes 1$. \\

Let $G_1$ and $G_2$ be two non-comparable elements of $\G\setminus \{\un\}$ forming a nested set. We have to show the equality of morphisms
 \begin{equation*}
 (\OS(\red{[G_1 \vee G_2]}) \otimes \Id ) \circ \OS(\red{[G_2]}) = \sigma^{2,3}\circ (\OS(\red{[G_1\vee G_2 ]}) \otimes \Id) \circ \OS(\red{[G_1]}).
 \end{equation*}
For any atom $H$ in $\L$, if $H \leq G_1$ then by nestedness $H\nleq G_2$ and consequently both morphisms send $e_H$ to $1\otimes e_H \otimes 1$. If $H \leq G_2$ then by the same argument both morphisms send $e_H$ to $1\otimes 1 \otimes e_H$. Finally if $H \nleq G_1$ and $H \nleq G_2$ then both morphisms send $e_H$ to $e_{G_1\vee G_2 \vee H} \otimes 1 \otimes 1$ which finishes the proof.\\

Lastly, we need to prove the equality
\begin{equation*}
\OS (\red{[f(G)]}) \circ \OS(f) = (\OS(f_{|[G, \un]})\otimes \OS(f_{|[\zero,G]})) \circ \OS(\red{[G ]})
\end{equation*}
for every isomorphism $f: (\L', \G') \xrightarrow{\sim} (\L, \G)$ in $\LBScat$ and $G \in \G \setminus \{\un \}$. For any atom $H$ of $\L$, if $H \leq G$ then both morphisms send $e_H$ to $1\otimes e_{f(H)}$. If on the contrary $H \nleq G$ then both morphisms send $e_H$ to $e_{f(G) \vee f(H)}\otimes 1$. \\

Finally, one can check that $\OS$ is strong monoidal, which finishes the proof.
\end{proof}

\section{Gröbner bases for operads over $\LBS$}\label{secgrobner}
In this section we develop a theory of Gröbner bases for operads over $\LBS$.  \\

Classical Gröbner bases \cite{BW_1993} are a computational tool which is used to work out quotients of free associative algebras. The general idea is to start by introducing an order on generators of the free algebra. \textcolor{black}{That} order is then used to derive an order on all monomials, which is compatible in some sense with the multiplication of monomials (we call such orders ``admissible''). We then use \textcolor{black}{that} order to rewrite monomials in the quotient algebra:
\begin{equation*}
\textrm{greatest term } \longrightarrow \sum \textrm{lower terms},
\end{equation*}
for every relation $R = \textrm{greatest term} - \sum \textrm{lower terms}$ in some subset $\B$ of the quotient ideal (usually the greatest term is called the ``leading term'' and we will use \textcolor{black}{that} denomination). The subset $\B$ is called a Gröbner basis when it contains ``enough'' elements. To be precise we want that every leading term of some relation in the ideal is divisible by the leading term of some element of $\B$. The goal is to find a Gröbner basis as \red{small} as possible so that the rewriting is as easy as possible. At the end of the rewriting process (which stops if the monomials are well-ordered) we are left with all the monomials which are not rewritable i.e. which are not divisible by a leading term of some element of $\B$. Those monomials are called ``normal'' and they form a linear basis of our algebra exactly when $\B$ is a Gröbner basis. \textcolor{black}{That} basis comes with multiplication tables given by the rewriting process. \\

It turns out that this general strategy can be applied to structures which are much more general and complex than associative algebras, such as operads for instance. Loosely speaking, all we need in order to implement this strategy is to be able to make (reasonable) sense of the key words used above, such as ``monomials'', ``admissible orders'' and ``divisibility between monomials''. For operads over a Feynman category, the only non-trivial part is to construct admissible orders on monomials out of orders on generators. The main issue comes from the symmetries, because usually the compatibility with symmetries is too strong and prevents us from finding any admissible order. In order to circumvent this problem, drawing inspiration from the case of classical operads which was sorted out by Dotsenko and Khoroshkin in \cite{DK_2010}, we introduce a notion of a ``shuffle'' operad over $\LBS$. \red{As we have seen in Example \ref{exshuffle}, the sub-Feynman category of $\LBS$ obtained by restricting to partition lattices and minimal building sets, and forgetting the symmetries, already encodes ``shuffle'' operads, so one could be tempted to define a shuffle $\LBS$-operad simply as an $\LBS$-operad without symmetries. However, this definition would not be satisfactory when trying to define admissible order on monomials. Essentially, we would like to a have a bit more data than just built lattices to construct those orders in a functorial way. In the case of classical shuffle operads this additional data is the linear ordering of the entries.}

\subsection{Shuffle operads over $\LBS$}
In view of what has just been said we make the following definition.
\begin{madef}
A \emph{directed built lattice} is a triple $(\L, \G, \vartriangleleft)$ where $(\L, \G)$ is a built lattice and $\vartriangleleft$ is a linear order on the set of atoms of $\L$. A directed built lattice $(\L, \G, \vartriangleleft)$ is said to be \emph{irreducible} if $(\L, \G)$ is irreducible.
\end{madef}{\color{black} \begin{nota}
For all geometric lattices $\L$ and $X$ any element of $\L$, we will denote by $\At(\L)$ the set of atoms of $\L$ and $\At_{\leq}(X)$ the set of atoms of $\L$ below $X$. 
\end{nota}}
\red{We are going to redo the same constructions as in Sections \ref{secprelim} and \ref{secfeycat}  (induced building set, monoidal structure on built lattices, Feynman category)} but with directed built lattices instead of built lattices.
\begin{madef}
Let $(\L, \G, \vartriangleleft)$ be a directed built lattice and $[G_1, G_2]$ be an interval of~$\L$. The interval $[G_1, G_2]$ admits an induced directed built lattice structure given by the building set $\Ind_{[G_1, G_2]}(\G)$ and the linear order $\vartriangleleft_{\ind}$ defined by
\begin{equation}\label{indor}
K_1 \vartriangleleft_{\ind} K_2 \Leftrightarrow \min \, \{H \,  | \, G_1\vee H = K_1, H \mathrm{ atom }\} \vartriangleleft \min \, \{H \, | \, G_1 \vee H = K_2, H \, \mathrm{ atom }\}
\end{equation}
for any pair of elements $K_1$, $K_2$ covering $G_1$ (i.e. atoms of $[G_1, G_2]$).
\end{madef}
Both minima in \eqref{indor} are well defined by geometricity of $\L$. As in the case of built lattices, doing a double induction is the same as doing a single induction directly on the smallest interval (Lemma \ref{lemmadoubleinduction}).
\begin{madef}
The monoidal product $\otimes$ on directed built lattices is defined by
\begin{equation*}
(\L_1, \G_1, \vartriangleleft_1)\otimes (\L_2, \G_2, \vartriangleleft_2) = (\L_1\times\L_2, \G_1\times \{0\} \sqcup \{0\}\times \G_2, \vartriangleleft)
\end{equation*}
where $\vartriangleleft$ is defined by putting all the atoms of $\L_1$ before the atoms of $\L_2$.
\end{madef}
\red{Recall that in this article, every nested set is assumed to be spanning, i.e. supposed to contain the maximal elements of its ambient building set.} For any nested set $\S$ and $G$ an element of $\S$ we have defined in Subsection \ref{secFeycons} the notation $\tau_{\S}(G) \coloneqq \bigvee \S_{<G}$. If $\S$ is an (ordered) nested set in a directed built lattice $(\L, \G, \vartriangleleft)$ we denote
\begin{equation*}
(\L_{\S}, \G_{\S}, \vartriangleleft_{\S}) \coloneqq \bigotimes_{G \in \S} ([\tau_{\S}(G), G], \Ind(\G), \vartriangleleft_{\ind}).
\end{equation*}
We are now able to make the following definition.
\begin{madef}
The Feynman category $\LBS_{\sh}$ is the triple $(\nucat, \fcat, \imath)$ where
\begin{itemize}
\item The objects of the groupoid $\nucat$ are the directed irreducible lattices and its morphisms are defined by
\begin{multline*}
\Mor_{\nucat}((\L, \G, \vartriangleleft), (\L', \G', \vartriangleleft')) = \{ f:\L' \xrightarrow{\sim} \L \, | \, f \textrm{ poset isomorphism,} \\ \textrm{increasing with respect to} \vartriangleleft \textrm{and} \vartriangleleft' \textrm{ and such that } f(\G') = \G\}.
\end{multline*}
\item \red{The objects of the category $\fcat$ are monoidal products of directed built lattices and its morphisms are generated by structural morphisms $(\L_{\S}, \G_{\S}, \vartriangleleft_{\S}) \xrightarrow{\S} (\L, \G, \vartriangleleft)$ together with the tensored isomorphisms of $\nucat$ and the permutations of monoidal summands. Those morphisms are subject to the relations
\begin{multline}
(\L'_{f(\S)}, \G'_{f(\S)}, \vartriangleleft'_{f(\S)}) \overset{f(\S)}{\rightarrow}(\L', \G', \vartriangleleft') \overset{f}{\rightarrow}(\L, \G, \vartriangleleft)  \sim \\ (\L'_{f(\S)}, \G'_{f(\S)}, \vartriangleleft'_{f(\S)}) \overset{\otimes f}{\rightarrow} (\L_{\S}, \G_{\S}, \vartriangleleft_{\S}) \overset{\S}{\rightarrow} (\L, \G, \vartriangleleft)
\end{multline}
for every isomorphism of irreducible directed built lattice\red{s} $f: (\L', \G', \vartriangleleft'~) \rightarrow (\L, \G, \vartriangleleft)$ with $\S$ any spanning nested set in $(\L, \G)$, and relations 
\begin{equation}
(\L_\S, \G_\S, \vartriangleleft_{\S}) \overset{\sigma}{\rightarrow} (\L_{\S_{\sigma}}, \G_{\S_{\sigma}}, \vartriangleleft_{\S_{\sigma}}) \overset{\S_{\sigma}}{\rightarrow} (\L, \G) \sim (\L_\S, \G_\S, \vartriangleleft) \overset{\S}{\rightarrow} (\L, \G, \vartriangleleft)
\end{equation}
for any permutation $\sigma$ of $\S$. We also impose $\Id_{(\L, \G, \vartriangleleft)} \sim \{\un\}$ for every directed built lattice $(\L, \G, \vartriangleleft)$.}

\item The functor $\imath$ is the obvious inclusion.
\end{itemize}
\end{madef}
We define the composition of nested sets in this context exactly as it was defined in Section \ref{secFeycons}. In subsequent sections a shuffle $\LBS$-operad will mean an operad over the Feynman category $\LBS_{\sh}$.

\subsection{Relating $\LBS$-operads and shuffle $\LBS$-operads}
There is an obvious functor between Feynman categories:
\begin{align*}
\sh: \LBS_{\sh} &\rightarrow \LBS \\
(\L, \G, \vartriangleleft) &\rightarrow (\L, \G) \\
\S &\rightarrow \S.
\end{align*}
This allows us to define a forgetful functor from $\LBS$-operads/modules to shuffle $\LBS$-operads/modules by precomposition.
\begin{madef}
For any $\LBS$-operad (resp. module) $\mathbb{P}$ we define the shuffle $\LBS$-operad (resp. module) $\mathbb{P}_{\sh}$  by
\begin{equation*}
\mathbb{P}_{\sh} \coloneqq \mathbb{P} \circ \sh.
\end{equation*}
\end{madef}
As in the classical operadic case \textcolor{black}{that} functor enjoys very nice properties which are listed and proved below.
\begin{prop}\label{propmainshuffle}
\begin{enumerate}
\item For any $\LBS$-module $M$ we have an isomorphism of shuffle operads:
\begin{equation}\label{eqisoopshuffle}
\LBS_{\sh}(M_{\sh}) \simeq \LBS(M)_{\sh}.
\end{equation}
\item Let $R$ be a sub-$\LBS$-module in some free $\LBS$-operad $\LBS(M)$ with $M$ some $\LBS$-module. The $\LBS$-module $R_{\sh}$ can be identified with a sub-$\LBS_{\sh}$-module of the free $\LBS_{\sh}$-operad $\LBS_{\sh}(M_{\sh})$. The $\LBS_{\sh}$-module $\langle R \, \rangle_{\sh}$ can be identified with an ideal of $\LBS(\Gen)_{\sh}$. The isomorphism \eqref{eqisoopshuffle} sends (via identifications) the (shuffle) ideal $\langle R_{\sh}\rangle$ to the shuffle ideal $\langle R \,\rangle_{\sh}$. As a consequence we have the isomorphism of $\LBS_{\sh}$-operads
\begin{equation*}
\LBS_{\sh}(M_{\sh})/ \langle R_{\sh}\rangle \xrightarrow{\sim} \LBS(M)_{\sh}/\langle R \,\rangle_{\sh}.
\end{equation*}
\item Lastly, we have an isomorphism of shuffle operads
\begin{equation*}
\LBS(M)_{\sh}/\langle R \,\rangle_{\sh} \simeq (\LBS(M)/\langle R\, \rangle)_{\sh}\,\,.
\end{equation*}
\end{enumerate}
\end{prop}
\begin{proof}
\begin{enumerate}
\item By the definition of left Kan extensions we have a morphism of $\LBS$-modules
\begin{equation*}
M \rightarrow \LBS(M).
\end{equation*}
Applying the forgetful functor we get a morphism of $\LBS_{\sh}$-modules
\begin{equation*}
M_{\sh} \rightarrow \LBS(M)_{\sh}.
\end{equation*}
By universal property of free $\LBS_{\sh}$-operads (which comes from the universal property of left Kan extensions) we get a morphism of $\LBS_{\sh}$-operads
 \begin{equation*}
 \LBS_{\sh}(M_{\sh}) \rightarrow \LBS(M)_{\sh}.
 \end{equation*}
 Let us take a closer look at \textcolor{black}{that} morphism. \\

By unpacking the construction of free operads (left Kan extensions) we get the following formula for every irreducible directed built lattice $(\L, \G, \vartriangleleft)$:
\begin{equation*}
\LBS_{\sh}(M_{\sh})(\L, \G, \vartriangleleft) \simeq \bigoplus_{\substack{\S \subset \G \\ \S \textrm{ spanning } \\ \textrm{ nested set }}} M_{\S}.
\end{equation*}
We also have
\begin{equation*}
\LBS(M)_{\sh}(\L, \G, \vartriangleleft) = \LBS(M)(\L, \G) = \left( \bigoplus_{\otimes_i (\L_i, \G_i) \rightarrow (\L, \G)} \bigotimes_i M(\L_i, \G_i)\right) /\sim
\end{equation*}
where $\sim$ identifies components corresponding to equivalent maps via isomorphisms, as explained at the beginning of the proof of Theorem \ref{theopresopFY}. The above morphism sends the component $M(\S)$ to the equivalence class of the component $M(\S)$. \\

However by Lemma \ref{lemmadecmor} the nested sets in $\L$ with linear order $\indtotor$ form a system of representatives for the equivalence classes of morphisms which means that the above morphism is indeed a linear isomorphism in each arity.
\item For the first identification if we start with the injective morphism $R \hookrightarrow \LBS(M)$, then apply the forgetful functor (which preserves injective morphisms since it is the right adjoint to left Kan extension), and then compose with isomorphism \eqref{eqisoopshuffle} we get an injective morphism $R_{\sh} \hookrightarrow \LBS_{\sh}(M_{\sh})$. For the second identification we have an injective morphism $\langle R\, \rangle \hookrightarrow \LBS(M)$ and applying the forgetful functor gives us an injective morphism $\langle R\, \rangle_{\sh} \hookrightarrow \LBS(M)_{\sh}$. By unraveling again the explicit description of isomorphism \eqref{eqisoopshuffle} we see that it sends the shuffle ideal $\langle R_{\sh}\rangle $ to the shuffle ideal $\langle R \, \rangle_{\sh}$.
\item We have an obvious identification of components in each arity between the two shuffle operads and one can check that this identification is operadic.
\end{enumerate}
\end{proof}

\subsection{Monomials and divisibility}
In order to be able to talk about Gröbner bases one needs a suitable notion of ``monomials'' in a free $\LBS_{\sh}$-operad as well as a notion of divisibility between those monomials. Let $M$ be a module over $\LBS_{\sh}$ in the category of vector spaces over some fixed arbitrary field $\K$. We have defined in Section \ref{secfreeop} the free $\LBS_{\sh}$-operad $\LBS_{\sh}(M)$ generated by the module $M$. By the same analysis conducted in the latter section we have the explicit formula for each directed irreducible built lattice $(\L, \G, \vartriangleleft)$:
\begin{equation}\label{eqfreeshop}
\LBS_{\sh}(M)(\L, \G, \vartriangleleft)  = \bigoplus_{\substack{\S \subset  \G \\ \S \textrm{ nested set }}}\bigotimes_{G \in \S}M([\tau_{\S}(G), G], \Ind(\G), \vartriangleleft_{\ind}).
\end{equation}
If furthermore we are given a basis $\Basis(\L, \G, \vartriangleleft)$ of every vector space $M(\L, \G, \vartriangleleft~)$ then we can make the following definition.
\begin{madef}[Monomial]
A \emph{monomial} in $\LBS_{\sh}(M)$ is an element which is a tensor of elements of the basis $\bigsqcup_{(\L', \G', \vartriangleleft')} \Basis(\L', \G', \vartriangleleft')$.
\end{madef}
In other words a monomial corresponds to the datum $(\S, (e_{G})_{G \in \S})$ of a nested set $\S$ and elements of the basis $e_{G} \in \Basis([\tau_{\S}(G), G], \Ind(\G), \vartriangleleft_{\ind})$ for each $G$ in $\S$. Monomials are stable under composition in $\LBS_{\sh}(M)$ and by \eqref{eqfreeshop} they form a basis of every vector space $\LBS_{\sh}(M)(\L, \G, \vartriangleleft)$. Additionally, we have a notion of divisibility between monomials.
\begin{madef}[Division between monomials]
Let $m_1$ and $m_2$ be two monomials in $\LBS_{\sh}(M)(\L, \G, \vartriangleleft)$. We say that $m_1$ \emph{divides} $m_2$ if $m_2$ can be expressed as a composition:
\begin{equation*}
m_2 = \LBS_{\sh}(M)(\S)((\alpha_G)_{G \in \S})
\end{equation*}
for some nested set $\S$, where one of the $\alpha_{G}$'s is $m_1$ and the rest are elements of the basis $\bigsqcup_{(\L, \G, \vartriangleleft)} B(\L, \G, \vartriangleleft)$.
\end{madef}

\subsection{An admissible order on monomials}
Let $M$ be an $\LBS_{\sh}$-module, with a basis $\Basis(\L, \G, \vartriangleleft)$ of $M(\L, \G, \vartriangleleft)$ for each $(\L, \G, \vartriangleleft)$. Assume that we have a total order $\dashv$ of those bases in each arity. In this section we will construct an order on monomials induced by $\dashv$, which is compatible with the composition of monomials in the sense of Proposition \ref{propadmissible}. We start by defining a total order $\indtotor$ on $\L$, induced by the direction $\vartriangleleft$. For any element $G$ in $\L$ we denote by $w(G)$ the word in the alphabet $\At(\L)$ given by the list of atoms below $G$ in increasing order.
\begin{madef}
Given two elements $G_1$, $G_2$ in $\L$, we say $G_1 \indtotor G_2$ if $w(G_2)$ is an initial subword of $w(G_1)$ or if $w(G_1)$ is less than $w(G_2)$ for the lexicographic order.
\end{madef}

There are a few important lemmas/remarks to be made about \textcolor{black}{that} order. First we prove that $\indtotor$ is compatible with the already existing order on $\L$.
\begin{lemma}\label{lemmaext}
The total order $\indtotor$ extends the reversed lattice order on $\L$.
\end{lemma}
\begin{proof}
Let $G_1 < G_2$ be two comparable elements in $\G$, such that $w(G_1)$ is not an initial subword of $w(G_2)$. One can write
\begin{align*}
w(G_1) &= u H_1 w_1 \\
w(G_2) &= u H_2 w_2
\end{align*}
with $u, w_1, w_2$ some words in the alphabet $\At(\L)$ and $H_1, H_2$ two different atoms. Since we have $\At_{\leq}(G_1) \subset \At_{\leq}(G_2)$ we immediately get $H_2 \vartriangleleft H_1$ which shows that $w(G_2)$ is smaller than $w(G_1)$ for the lexicographic order.
\end{proof}
Next we prove that $\indtotor$ behaves well with respect to restriction to intervals.
\begin{lemma}\label{lemmaorderres}
Assume we are given an interval $[G_1, G_2]$ in some irreducible directed built lattice $(\L, \G, \vartriangleleft)$. The two total orders $\indtotor_{\ind}$ and $\indtotor_{|[G_1, G_2]}$ on $[G_1, G_2]$ are the same.
\end{lemma}
\begin{proof}
Since we are comparing two total orders we only need to prove one implication, for instance $ K \indtotor_{|[G_1, G_2]}  K' \Rightarrow K \indtotor_{\ind}K'$. If $w(K')$ is included in $w(K)$ then this means that we have $K' \leq K$ which proves that we have $K \indtotor_{\ind}K'$ by the previous lemma. If not then one can write $w(K) = u H w$ and $w(K') = u H' w'$ with $u$, $w$ and $w'$ some words and $H \vartriangleleft H'$ two atoms of $\L$. The atom $H$ cannot be below $G_1$ otherwise $H$ would also be a letter in $w(K')$. If $H'$ is below $G_1$ then let $H''$ be the first letter bigger than $H'$ in $w(K')$ which is not below $G_1$ and such that $G_1 \vee H''$ does not belong to $G_1 \vee u$ (such a letter exists because we have assumed that $w(K')$ is not included in $w(K)$). In \textcolor{black}{that} case one can write
\begin{align*}
&w_{[G_1, G_2]}(K) = v (G\vee H) t \\
&w_{[G_1, G_2]}(K') = v (G\vee H'') t'
\end{align*}
with $v, t, t'$ some words in $\At([G_1, G_2])$. This implies that we have $K \indtotor_{\ind} K'$.

\end{proof}

Finally, we prove that $\indtotor$ is compatible with the join in $\L$.
\begin{lemma}\label{lemmajoin}
Let $G$, $G_1$ and $G_2$ be three elements in $\L$ such that $\Fact_{\G}(G_1)$ and $\Fact_{\G}(G_2)$ are both disjoint from $\Fact_{\G}(G)$ and such that $\Fact_{\G}(G_1) \cup \Fact_{\G}(G)$ and $\Fact_{\G}(G_2) \cup \Fact_{\G}(G)$ are both nested antichains. Then we have the equivalence
\begin{equation*}
G_1 \indtotor G_2 \Leftrightarrow G \vee G_1 \, \indtotor \, G \vee G_2.
\end{equation*}
\end{lemma}
\begin{proof}
Let us start by proving the direct implication. If $G_1 > G_2$ then we have $G \vee G_1 > G \vee G_2$ (the strictness coming from the nestedness condition) which proves the result by Lemma \ref{lemmaext}. Otherwise write $w(G_1) = u H_1 w_1$ and $w(G_2) = uH_2 w_2$ where $u$, $w_1$ and $w_2$ are some words and $H_1$ is strictly smaller than $H_2$. By the nestedness condition in the proposition we have
\begin{align*}
w(G\vee G_1) &= \textrm{sh}(w(G), w(G_1)), \\
w(G\vee G_2) &= \textrm{sh}(w(G), w(G_2)),
\end{align*}
where $\textrm{sh}(\cdot, \cdot)$ is the operation which merges two given words with increasing letters into a word with increasing letters. From this we see that one can write
\begin{align*}
&w(G\vee G_1) = u' H_1 w'_1, \\
&w(G\vee G_2) = u' H'_2 w'_2,
\end{align*}
with $u'$, $w'_1$, $w'_2$ some words and $H'_2$ an atom of $\L$. If $H'_2$ is in $w(G)$ then $H_1$ is strictly smaller than $H'_2$ (otherwise $H'_2$ would belong to $u'$). If on the contrary $H'_2$ belongs to $w(G_2)$ then $H'_2$ is equal to $H_2$ and $H_1$ is strictly smaller than $H'_2$. \\

For the converse, assume we have $G\vee G_1 \indtotor G\vee G_2$. If $G\vee G_2 < G\vee G_1$ then $G_2 < G_1$ by nestedness which implies $G_1 \indtotor G_2$ by Lemma \ref{lemmaext}. Otherwise, write $w(G\vee G_1) = u H_1 w_1$ and $w(G\vee G_2) = uH_2 w_2$ for $u, w_1, w_2$ some words and $H_1$ strictly smaller than $H_2$. One can check that $H_1$ necessarily belongs to $w(G_1)$ which means that we can write $w(G_1) = u' H_1 w'_1$ and $w(G_2) = u' H'_2 w'_2$ where $H'_2$ is the first letter of $w(G\vee G_2)$ which belongs to $w(G_2)$ and which comes after $H_2$. We immediately get $H_1 \vartriangleleft H_2 \trianglelefteq H'_2$ which finishes the proof.\end{proof}
\smallskip
For any monomial $m = (\S, (e_G)_G)$ and $G_0$ some element in $\min_{\leq}\S$, we denote $$G_0\vee m \coloneqq (G_0 \vee (\S\setminus\{G_0\}), (e_{\Comp_{G_0}(G)})_G),$$ which is a well-defined monomial in $\LBS_{\sh}(M)([G_0, \un], \Ind(\G), \vartriangleleft_{\ind})$, by Lemma \ref{lemmarecons}. For any nested set $\S$ we denote $\mathsf{MM}(\S) \coloneqq \min_{\indtotor_{\dashv}} \min_{\leq} \S$.
\begin{madef}
We define a total order $\indtotor_{\dashv}$ on monomials in the following inductive way. For $m_1 = (\S_1, (e^1_{G})_{G \in \S_1})$ and $m_2 = (\S_2, (e^2_{G})_{G \in \S_2})$ we put $m_1 \indtotor_{\dashv} m_2$ if there exists some $G$ in $\min_{\leq} S_1 \cap \min_{\leq}S_2$ such that $e^1_G = e^2_G$ and $G\vee m_1 \indtotor_{\dashv} G\vee m_2$, or if there is no such $G$ and $\MM(\S_1) \indtotor \MM(\S_2)$ or $\MM(\S_1) = \MM(\S_2)$ and $e^1_{\MM(\S_1)} \dashv e^2_{\MM(\S_2)}$.
\end{madef}
One can check that this definition does not depend on the choice of the element $G$ because if we have two different elements $G$ and $G'$ in $\min_{\leq}\S_1 \cap \min_{\leq}\S_2$ such that $e^1_G = e^2_G$ and $e^1_{G'} = e^2_{G'}$ then we have the equality of monomials $$(G \vee G')\vee(G \vee m_i) = (G \vee G')\vee(G' \vee m_i)$$ for $i = 1, 2$. If there is no ambiguity on the order $\dashv$ we write $\indtotor$ instead of $\indtotor_{\dashv}$.
\begin{prop}\label{propadmissible}
The order on monomials $\indtotor$ is compatible with the composition of monomials. More precisely, if $\S$ is any nested set in some directed irreducible built lattice $(\L, \G, \vartriangleleft)$ and we have some generators $e_{G} \in B([\tau_{\S}(G), G], \Ind(\G), \vartriangleleft_{\ind})$ for all $G$ in $\S$ except for $G_0$ in $\S$ where we have monomials $m_1 = (\S_1, (e^1_G)_G),m_2 = (\S_2, (e^2_G)_G) \in \LBS_{\sh}([\tau_{\S}(G_0), G_0], \Ind(\G), \vartriangleleft_{\ind}~)$, with $\#\S_1 = \#\S_2$, then we have
\begin{equation*}
m_1 \indtotor m_2 \Rightarrow \LBS_{\sh}(M)(\S)((e_G)_G, m_1) \indtotor \LBS_{\sh}(M)(\S)((e_G)_G, m_2).
\end{equation*}
\end{prop}
\begin{proof}
The proof goes by induction on $\#\S + \#\S_1$. The initialization at $\# \S = 0$ or $\#\S_1 = 0$ is obvious. The induction step is a consequence of Lemma \ref{lemmajoin}.
\end{proof}
Since every element in a free $\LBS_{\sh}$-operad can be uniquely written as a sum of monomials, we can make the following definition.

\begin{madef}[Leading term]
If $f$ is an element in some free $\LBS_{\sh}$-operad with total order on generators~$\dashv$ then the \textit{leading term} of $f$, denoted by $\lt(f)$, is the biggest monomial with respect to~$\indtotor_{\dashv}$ which has a non-zero coefficient in $f$.
\end{madef}

Finally, everything has been leading to the following definition.

\begin{madef}[Gröbner basis]
Let $\mathcal{I}$ be an operadic ideal in some free $\LBS_{\sh}$-operad $\LBS_{\sh}(M)$, where $M$ is some $\LBS_{\sh}$-module in some category of vector spaces which is endowed with a basis and a well-order $\dashv$ of \textcolor{black}{that} basis in each arity. A \emph{Gröbner basis} of $\mathcal{I}$ relative to ~$\dashv$ is a subset $\mathcal{B}$ of $\I$ such that every leading term relative to $ \indtotor_{\dashv}$ of some element of $\I$ is divisible by the leading term of some element of $\mathcal{B}$.
\end{madef}
A Gröbner basis is said to be quadratic if it contains only degree $1$ elements.
\begin{madef}[Normal monomial]
A \emph{normal monomial} with respect to some set of elements $\B$ is a monomial which is not divisible by the leading term of some element in $\B$.
\end{madef}
\begin{prop}\label{propgbasisnm}
The set of normal monomials with respect to some set of elements $\B$ in an ideal $\I \subset \LBS_{\sh}(M)$ linearly generates $\LBS_{\sh}(M)/ \I$ in every arity. \textcolor{black}{That} set of monomials is linearly \red{independent} if and only if $\B$ is a Gröbner basis of $\I$.
\end{prop}{\color{black}\begin{proof}
The proof is the same as in any other context involving Gröbner bases. Let us summarize the main arguments. For the first assertion, any element of $\LBS_{\sh}(M)/ \I$ in some arity can be expressed as a sum of monomials. One can repeatedly rewrite the greatest non normal monomial of that sum. Since our linear basis in each arity is well-ordered, the monomials are also well-ordered, which ensures that this process will eventually terminate. The final result will be a sum of normal monomials which is equal in $\LBS_{\sh}(M)/ \I$ to our initial element.  For the second assertion, let $\B$ be a Gröbner basis of $\I$. If there is a nonzero linear relation among the normal monomials of $\B$, then that relation must belong to $\I$, which implies that the leading term of that relation cannot be normal with respect to $\B$, which is a contradiction. Finally, consider $\B$ such that its associated normal monomials are linearly independent in $\LBS_{\sh}(M)/ \I$. If there exists a nonzero element of $\I$ whose leading term is normal, then as previously one can repeatedly rewrite the greatest non normal monomial of this element until we get a sum of normal monomials in $\I$. Since the leading term has not changed, the sum is nonzero, which is a contradiction. \end{proof}}

\subsection{Application: the example of $\FY^{\vee}$}\label{secgrobnerFY}
We have proved in Subsection \ref{secfreeop} that if $\Gen$ is the $\LBS$-module with one generator $\Psi_{(\L, \G)}$ of degree $2(\rk(\L) - 1)$ in each arity $(\L, \G)$, and $\I$ is the ideal generated by the elements
\begin{multline}\label{eqgens}
\sum_{G. \geq H_1} \LBS(\Gen)(\red{[G]})(\Psi_{([G, \un ], \Ind(\G))}, \Psi_{([\zero, G], \Ind(\G))}) - \\
 \sum_{G \geq H_2} \LBS(\Gen)(\red{[G]})(\Psi_{([G, \un ], \Ind(\G))}, \Psi_{([\zero, G], \Ind(\G))})
\end{multline}
for each pair of atoms $H_1$ and $H_2$, then we have a surjective morphism of $\LBS$-operads:
\begin{equation*}
\LBS(\Gen)/ \I \xrightarrow{\pi} \FY^{\vee}.
\end{equation*}
Let us denote by $R$ the linear span of the elements \eqref{eqgens}, which is a sub $\LBS$-module of $\LBS(\Gen)$. By Proposition \ref{propmainshuffle} we have an isomorphism of shuffle $\LBS$-operads
\begin{equation*}
\LBS_{\sh}(\Gen_{\sh})/\langle R_{\sh} \rangle  \xrightarrow{\sim} \LBS(\Gen)_{\sh}/\langle R\, \rangle_{\sh}.\\
\end{equation*}

Let us use our theory of Gröbner bases for shuffle $\LBS$-operads to study the operad $\LBS_{\sh}(\Gen_{\sh})/\langle R_{\sh}\rangle$. Notice that $R_{\sh}$ is just the linear span of elements of the form
\begin{multline*}
\sum_{G \geq H_1} \LBS_{\sh}(\Gen_{\sh})(\red{[G]})(\Psi_{([G, \un ], \Ind(\G))}, \Psi_{([\zero, G], \Ind(\G))}) -  \\
\sum_{G \geq H_2} \LBS_{\sh}(\Gen_{\sh})(\red{[G]})(\Psi_{([G, \un ], \Ind(\G))}, \Psi_{([\zero, G], \Ind(\G))}).
\end{multline*}
We denote by $\B$ the set of those elements, and \red{for every spanning nested set $\S$} we put
\begin{equation*}
\Psi_{\S}\coloneqq  \LBS_{\sh}(\Gen_{\sh})(\S)((\Psi_{([\tau_{\S}(G), G], \Ind(\G))})_{G \in \S}).
\end{equation*}
\red{By Proposition \ref{propmainshuffle}, the surjectivity of $\pi$ implies that the corresponding morphism of $\LBS_{\sh}$-operads
\begin{equation}\label{eqsurjmor}
\LBS_{\sh}(\Gen_{\sh})/ \langle R_{\sh} \rangle  \longrightarrow  \FY^{\vee}_{\sh}. \\
\end{equation}
is also surjective.} We will compute the normal monomials associated to $\B$ and find that they have the desired cardinality, which will prove that $\B$ forms a Gröbner basis of $\langle R_{\sh}\rangle$, and that morphism \eqref{eqsurjmor} is an isomorphism. \red{By Proposition \ref{propmainshuffle} that will imply that morphism $\pi$ is also an isomorphism.}\\

To describe those normal monomials in a natural way we introduce a classical tool in poset combinatorics called ``EL-labeling''.
\begin{madef}[EL-labeling]
Let $P$ be a finite poset with set of covering relations $\Cov(P)$.  An \emph{EL-labeling} of $P$ is a map $\lambda:\Cov(P) \rightarrow \N$ such that for any two comparable elements $X < Y$ in $P$ there exists a unique maximal chain going from $X$ to $Y$ which has increasing $\lambda$ labels (when reading the covering relations from bottom to top) and this unique maximal chain is minimal for the lexicographic order on maximal chains (comparing the words given by the successive $\lambda$ labels from bottom to top).
\end{madef}
We refer the reader to \cite{wachs_poset_2006} for more details on \textcolor{black}{that} notion. The main result we will use about EL-labelings is the following.
\begin{prop}\label{propelgeo}
Let $\L$ be a geometric lattice. Any linear ordering $H_1 \vartriangleleft ... \vartriangleleft H_n$ of the atoms of $\L$ induces an $\EL$-labeling $\lambda_{\vartriangleleft}$ of $\L$ defined by
\begin{equation*}
\lambda_{\vartriangleleft}(X \prec Y) = \min \{i \, | \,  X \vee H_i = Y\}
\end{equation*}
for any covering relation $X \prec Y$ in $\L$.
\end{prop}
\begin{proof}
The proof of this result can be found in \cite{wachs_poset_2006}.
\end{proof}
If $\lambda$ is an EL-labelling of some poset $P$ and $X < Y$ are two comparable elements we denote by $\omega_{X, Y, \lambda}$ the unique maximal chain from $X$ to $Y$ which has increasing $\lambda$ labels. If the EL-labelling can be deduced from the context we will drop it from the notation. We also define $\omega^{k}_{X, Y, \lambda}$ to be the the chain $\omega_{X, Y, \lambda}$ truncated at height $k$ for any positive integer $k$ which is less than the length of $\omega_{X,Y, \lambda}$. More precisely if $\omega_{X, Y, \lambda} = \{X_0 = X \prec X_1 \prec ... \prec X_n = Y\}$ then $\omega^k_{X,Y, \lambda} \coloneqq \{X_0 \prec ... \prec X_k\}$. We will also need a new definition in nested set combinatorics.
\begin{madef}[Cluster]
A spanning nested set $\S$ is called a \textit{cluster} if all its intervals except the top one have rank $1$. A cluster is said to be \emph{proper} if its top interval has rank strictly greater than $1$.
\end{madef}
Clusters can be constructed out of truncated maximal chains in the following way. If $\omega$ is a chain $\omega = \{ X_0 = \zero \prec X_1 \prec ... \prec X_n \}$ in some built geometric lattice $\L$ then we put:
\begin{equation*}
\S(\omega) \coloneqq \{X_1\}\circ \{X_2\} \circ ... \circ \{X_n\},
\end{equation*}
which is a cluster. This formula makes sense even if the $X_i$'s do not belong to $\G$ because each $X_i$ covers $X_{i-1}$ and therefore is an atom in $[X_{i-1}, \un]$ which must belong to the induced building set. We can finally state the main result of this section.
\begin{prop}\label{propnormalmono}
The normal monomials with respect to $\B$ in arity $(\L, \G, \vartriangleleft)$ are the monomials of the form
\begin{equation*}
\LBS_{\sh}(\Gen_{\sh})(\S)\left(\Psi_{\S(\omega^{k_{G}}_{\tau_{\S}(G), G, \lambda_{\vartriangleleft}})}\right)
\end{equation*}
where $\S$ is some nested set without any rank $1$ intervals and the $k_G$'s  are integers strictly less than $\rk( [\tau_{\S}(G), G]) - 1$, except $k_{\un}$ which can be equal to $\rk( [\tau_{\S}(\un), \un]) - 1$.  Furthermore \textcolor{black}{that} decomposition is unique.
\end{prop}
\begin{proof}
We start with the following lemma.
\begin{lemma}\label{lemmadeccluster}
Any spanning nested set $\S$ in some irreducible built lattice can be written as
\begin{equation*}
\S = \S'\circ(\S'_{G'})_{G' \in \S'}
\end{equation*}
where $S'$ is a spanning nested set with no rank $1$ intervals and the $\S'_{G'}s$ are proper clusters except $\S_{\un}$ which is any cluster.
\end{lemma}
The nested set $\S'$ will be called the frame of $\S$ and denoted by $\fr(\S)$.
\begin{proof}
We put $\fr(\S) = \{G \in \S \textrm{ s.t. } \rk([\tau_{\S}(G), G]) > 1\} \cup \{\un\}$ and conclude by Lemma \ref{lemmarecons}.
\end{proof}

Now let us proceed with the proof of the statement. We denote by $\M(\L, \G, \vartriangleleft~)$ the normal monomials with respect to $\B$ in arity $(\L, \G, \vartriangleleft)$ and $\M'(\L, \G, \vartriangleleft)$ the monomials of the form
\begin{equation*}
\LBS_{\sh}(\Gen)(\S)\left(\Psi_{\omega^{k_{G}}_{\tau_{\S}(G), G, \lambda_{\vartriangleleft}}}\right).
\end{equation*}
Our goal is to show the equality $\M(\L, \G, \vartriangleleft) = \M'(\L, \G, \vartriangleleft)$ for all irreducible directed built lattice $(\L, \G, \vartriangleleft)$. However we have a  bijection between $\M'(\L, \G, \vartriangleleft)$ and normal monomials of $\FYa(\L, \G)$ with respect to the Gröbner basis introduced in Theorem \ref{theogrobnerFY}, given by
\begin{equation*}
\LBS_{\sh}(\Gen)(\S)\left(\Psi_{\omega^{k_{G}}_{\tau_{\S}(G), G, \lambda_{\vartriangleleft}}}\right) \rightarrow \prod_{G \in \S} x_G^{\rk([\tau_{\S}(G), G]) - k_G - 1},
\end{equation*}
and we have the surjective morphism of $\LBS_{\sh}$-operads \eqref{eqsurjmor}. This means that it is enough to prove the inclusion $\M(\L, \G, \vartriangleleft) \subset \M'(\L, \G, \vartriangleleft)$ (see Proposition \ref{propgbasisnm}). By Lemma \ref{lemmadeccluster} it is enough to prove that for any cluster $\S$, if $\Psi_{\S}$ is a normal monomial then $\S$ is of the form $\omega^k_{\zero, \un, \lambda_{\vartriangleleft}}$.\\

Leading terms of elements of $\B$ are monomials of the form $\Psi_{\{H, \un\}}$ where $H$ is not the minimal atom. Let $\S$ be any cluster such that $\Psi_{\S}$ is a normal monomial, i.e. is not divisible by any $\Psi_{\{H, \un\}}$ with $H$ not minimal. For any $G \in \S\setminus{\un}$, let us denote by $H_G$ the smallest atom such that $\tau_{\S}(G)\vee H_G = G$.
\begin{lemma}
The map $G \rightarrow H_G$ is increasing (with respect to the order $\leq$ on the domain and the order $\vartriangleleft$ on the codomain).
\end{lemma}
\begin{proof}
Let $G_1 < G_2$ be two elements in $\S$ such that $\rk([\tau_{\S}(G_i), G_i] = 1$ for $i=1, 2$ and such that there is no element in $\S$ strictly between $G_1$ and $G_2$. By Lemma \ref{lemmarecons} we can write $\S = (\S\setminus\{G_1\})\circ \{\tau_{\S\setminus\{G_1\}}(G_2)\vee H_{G_1}\}$. Since $\Psi_{\S}$ is not divisible by any monomial $\Psi_{H}$ where $H$ is not minimal this means that $\tau_{\S\setminus\{G_1\}}(G_2)\vee H_{G_1}$ is the minimal atom in $[\tau_{\S\setminus\{G_1\}}(G_2), G_2]$, but \textcolor{black}{that} interval contains the atom $\tau_{\S\setminus\{G_1\}}(G_2)\vee H_{G_2}$ so we have $\tau_{\S\setminus\{G_1\}}(G_2)\vee H_{G_1} \vartriangleleft_{\ind}\tau_{\S\setminus\{G_1\}}(G_2)\vee H_{G_2}$. \\

Technically this means
\begin{multline}\label{eqminat}
\min \{H \, | \, \tau_{\S\setminus\{G_1\}}(G_2)\vee H = \tau_{\S\setminus\{G_1\}}(G_2) \vee H_{G_1}\} < \\ \min \{H \, | \, \tau_{\S\setminus\{G_1\}}(G_2) \vee H = \tau_{\S\setminus\{G_1\}}(G_2) \vee H_{G_2} \}.
\end{multline}
By nestedness of $\S$ we have
\begin{equation*}
\{H \, | \, \tau_{\S\setminus\{G_1\}}(G_2)\vee H = \tau_{\S\setminus\{G_1\}}(G_2) \vee H_{G_1}\} = \{H \, | \, \tau_{\S}(G_1)\vee H = G_1\}
\end{equation*}
which has minimum $H_{G_1}$, and
\begin{equation*}
\{H \, | \, \tau_{\S\setminus\{G_1\}}(G_2)\vee H = \tau_{\S\setminus\{G_1\}}(G_2) \vee H_{G_2}\} = \{H \, | \, \tau_{\S}(G_2) \vee H = G_2\}
\end{equation*}
which has minimum $H_{G_2}$. Inequality \eqref{eqminat} concludes the proof.
\end{proof}
Let us denote $\S = \{G_1, ..., G_n \}\cup \{\un\}$ with $H_{G_1}\vartriangleleft ...  \vartriangleleft H_{G_n}$. By the previous lemma and successive applications of Lemma \ref{lemmarecons} we get
\begin{equation*}
\S = \{H_{G_1}\}\circ\{H_{G_1}\vee H_{G_2}  \} \circ ... \circ \{H_{G_1} \vee ... \vee H_{G_n}\}.
\end{equation*}
What is left to prove is that the chain $H_{G_1} \prec ... \prec H_{G_1} \vee ... \vee H_{G_n}$ is exactly the chain $\omega^n_{\zero, \un, \lambda_{\vartriangleleft}}$. \\

We consider the concatenation of chains $$H_{G_1} \prec ... \prec H_{G_1} \vee ... \vee H_{G_n} \prec \omega_{H_{G_1} \vee ... \vee H_{G_n}\un\lambda_{\vartriangleleft_{\ind}}}.$$ This chain has increasing labels everywhere except possibly at $H_{G_1} \vee ... \vee H_{G_n}$. \\

By Lemma \ref{lemmarecons} we have $\S = (\S\setminus\{G_n\}) \circ \{\tau_{\S\setminus\{G_n\}}(\un) \vee H_{G_n}\}$ so if $\Psi_{\S}$ is a normal monomial this means that $\tau_{\S\setminus\{G_n\}}(\un) \vee H_{G_n}$ is the minimal atom in $[\tau_{\S\setminus\{G_n\}}(\un), \un]$. This means that $H_{G_n}$ is smaller than all the atoms which are not below $\tau_{\S\setminus\{G_n\}}(\un)$ and consequently the maximal chain introduced previously also has increasing labels at $H_{G_1}\vee ... \vee H_{G_n}$. \\

By Proposition \ref{propelgeo} the chain $H_{G_1} \prec ... \prec H_{G_1} \vee ... \vee H_{G_n} \prec \omega_{H_{G_1} \vee ... \vee H_{G_n}, \un,\lambda_{ \vartriangleleft_{\ind}}}$ must be the chain $\omega_{\zero, \un, \lambda_{\vartriangleleft}}$ and therefore $H_{G_1} \prec ... \prec H_{G_1} \vee ... \vee H_{G_n}$ is the chain $\omega^n_{\zero, \un, \lambda_{\vartriangleleft}}$ which concludes the proof.
\end{proof}

\begin{coro}\label{coroquadgrob}
The morphism
\begin{equation*}
\LBS(\Gen)/\I \xrightarrow{\pi} \FY^{\vee}
\end{equation*}
is an isomorphism and the shuffle $\LBS$-operad $\FY^{\vee}_{\sh}$ admits a quadratic Gröbner basis.
\end{coro}
{\color{black}\begin{proof}
In Subsection \ref{secfreeop} we have proved that the morphism $\pi$ is well-defined and surjective. Proposition \ref{propnormalmono} shows that the normal monomials associated to the subset $\B$ of $\I_{\sh}$ have cardinality the dimension of the Feichtner--Yuzvinsky algebras in each arity. By Proposition \ref{propgbasisnm} this implies that $\B$ is a quadratic Gröbner basis of $\I_{\sh}$ and that the morphism between shuffle $\LBS$-operads: $\LBS_{\sh}(\Gen_{\sh})/\I_{\sh} \rightarrow \FY_{\sh}^{\vee}$ is an isomorphism. By Proposition \ref{propmainshuffle} this implies that $\pi$ is an isomorphism. 
\end{proof}}

\section{Koszulness of $\LBS$-operads}\label{seckoszul}
In \cite{kaufmann_feynman_2017}, Kaufmann and Ward constructed a Koszul duality theory for operads over certain well-behaved Feynman categories called ``cubical''. \textcolor{black}{That} cubicality condition is what allows us to define odd operads and a cobar construction on odd operads, which is central in Koszul duality theory. \\

In the first subsection we prove that $\LBS$ is cubical. Then we unpack Koszulness for operads over $\LBS$, following the definitions in \cite{kaufmann_feynman_2017}. Finally, we prove that having a quadratic Gröbner basis implies being Koszul and we apply \textcolor{black}{that} result to $\FY^{\vee}$.\\

\subsection{$\LBS$ is cubical}
Let us start with some reminders on the notion of ``cubicality'' (we refer to \cite{kaufmann_feynman_2017} for more details). Given $(\nucat, \fcat, \imath)$ a graded Feynman category (see Section \ref{secgradedfeycat}) and $A$, $B$ two objects in $\fcat$ we denote by $C^{+}_n(A, B)$ the set of composable chains of morphisms of degree less \red{or equal} than $1$ having exactly $n$ morphisms of degree $1$, quotiented by relations:
\begin{equation}\label{eqrelisochaine}
A \rightarrow ... \rightarrow X_{i-1} \overset{f}{\rightarrow} X_{i} \overset{g}{\rightarrow} X_{i+1} \rightarrow ... \rightarrow B \sim A \rightarrow ... \rightarrow X_{i-1} \overset{g\circ f}{\rightarrow} X_{i+1} \rightarrow ... \rightarrow B
\end{equation}
provided $f$ or $g$ has degree $0$. There is a composition map going from $C^{+}_n(A, B)$ to $\Hom_{\fcat}(A, B)$ given by composing all the morphisms of the chain (the equivalence relation preserves this composition). This map will be denoted by $c_{A, B}$.
\begin{madef}[Cubical Feynman category]
A graded Feynman category $(\nucat, \fcat, \imath)$ is called \emph{cubical} if the degree function is proper and if for every $A$, $B$ objects of $\fcat$ there is a free $\Sym_n$ action on $C^{+}_n(A,B)$ such that
\begin{itemize}
\item The composition map is invariant over the action of $\Sym_n$.
\item The composition map defines a bijection $c_{A, B}: C^{+}_n(A,B)_{\Sym_n} \overset{\simeq}{\rightarrow} \Hom_{\fcat}(A, B)$.
\item The $\Sym_n$ action is compatible with concatenation of sequences (considering the inclusion $\Sym_p \times \Sym_q \subset~\Sym_{p+q}$).
\end{itemize}
\end{madef}
The cubicality of a Feynman category $\mathfrak{F}$ is a technical condition which ensures that $\mathfrak{F}$ is Koszul in the sense defined in \cite{kaufmann_feynman_2017} (see \cite{KW2021} Theorem 4.1), which is morally what we need in order to be able to define Koszulness of $\mathfrak{F}$-operads.
\begin{ex}
The Feynman category $\Op$, which encodes classical operads as described in Example \ref{exfeycatoperad}, is cubical. Indeed, a morphism in this Feynman category is a graph operation that can be represented by a tree. The degree $1$ generators of this Feynman category are exactly those morphisms which are represented by trees with one inner edge. Expressing a morphism of $\Op$ as a composition of degree $1$ generators is equivalent to choosing an order on the inner edges of the corresponding tree, which has an obvious symmetric action. Recall from example \ref{exshuffle} that shuffle trees can be viewed as nested sets in partition lattices with minimal building sets. There is an obvious bijection between the inner edges of the shuffle tree and the elements of the corresponding nested set. This suggests that $\LBS$ is also cubical, with symmetric action on chains of generators given by changing the order of the underlying nested set. This is the content of the next proposition. 
\end{ex}
\begin{prop}
The Feynman categories $\LBS$ and $\LBS_{\sh}$ are cubical.
\end{prop}
\begin{proof}
We only prove the result for $\LBS$, as the same arguments also work for $\LBS_{\sh}$. By Section \ref{presentation}, degree $0$ and degree $1$ morphisms generate every morphism in $\LBS$. Let us now define an explicit faithful symmetric action on $C_n^+(A,B)$ for every $A,B$ in $\LBS$. It is enough to define it for $B$ an irreducible built lattice. \\

Using relation \eqref{eqrelisoall} and relation \eqref{eqrelisochaine} we can see that every chain $\psi$ in $C_n^+(A,B)$ has a representative of the form
\red{\begin{equation}\label{eqdecisonested}
A \xrightarrow{\varphi } A' \xrightarrow{\phi} B,
\end{equation}
where $\phi$ is an element of $C_{n-1}^{+}(A', B)$ which contains only degree $1$ morphisms which are nested sets of cardinality one, and $\varphi$ is a morphism of degree $1$ of the form 
\begin{equation*}
    \varphi = \left(\bigotimes_i f_i \otimes \Bigl(\red{[G]} \circ (g_1 \otimes g_2)\Bigl) \otimes \bigotimes_j f_j \right) \circ\nu
\end{equation*}
for some element $G$, some isomorphisms $f_i$'s, $f_j's$, $g_1, g_2$ and $\nu$ some permutation of the summands of $A$.} By Lemma \ref{lemmadecmor} \textcolor{black}{that} representative of $\psi$ is in fact unique. We denote $\psi' = [ A'' \xrightarrow{\Id\otimes \red{[G]} \otimes \Id} A' \xrightarrow{\phi} B]$. The composition of the chain $\psi'$ is a nested set $\S$, and we have a linear ordering $G_1, G_2,..., G_n$ of $\S\setminus \{\un\}$ given by reading the chain $\psi'$ say from right to left. Let $\sigma$ be any element of $\Sym_n$. For all $i\leq n$ we put $G^{\sigma}_i \coloneqq G_{\sigma(i)}\vee \bigvee(\{G_{\sigma(1)},\ldots, G_{\sigma(i-1)}\}_{<G_{\sigma(i)}})$. We then define
\begin{equation*}
\sigma\cdot \psi' \coloneqq [A'' \xrightarrow{\sigma'} A''' \xrightarrow{\Id^{\otimes}\otimes [G^{\sigma}_n]\otimes \Id^{\otimes}} ... \xrightarrow{[G^{\sigma}_1]}B ]
\end{equation*}
where
\begin{itemize}
\item The formula $\Id^{\otimes} \otimes [G^{\sigma}_i] \otimes \Id^{\otimes}$ means that we tensor the morphism $[G^{\sigma}_i]$ by the identities of the summands of the codomain not containing $G^{\sigma}_i$. 
\item $\sigma'$ is the only permutation of the summands of $A''$ which gives us $A'''$. 
\end{itemize}
Finally, we put
\begin{equation*}
\sigma \cdot \psi \coloneqq [A \xrightarrow{\bigotimes_i f_i \otimes g_1 \otimes g_2 \otimes \bigotimes_j f_j} A'' \xrightarrow{\sigma\cdot \psi'} B ]. \\
\end{equation*}

First, by iterated use of Lemma \ref{lemmarecons} we see that the composition of the chain $\sigma\cdot \psi'$ is $\S$ and that the linear order on $\S\setminus\{\un\}$ given by reading $\sigma\cdot \psi'$ from right to left is exactly $\sigma$. This implies that $\cdot$ defines an action of $\Sym_n$. Second, let us remark that by relation \eqref{eqrelpermutation}, we have $c_{A'', B}(\sigma\cdot\psi') = c_{A'', B}(\psi')$ which implies $c_{A, B}(\sigma\cdot \psi) = c_{A,B}(\psi)$ i.e. the composition map is invariant by the action of $\Sym_n$. Third, we see that the action is free because of the \red{uniqueness} of decomposition \eqref{eqdecisonested}. Fourth and lastly, the composition map is bijective after passing to the quotient by the action of $\Sym_n$. The surjectivity immediately comes from the fact that morphisms of degree $0$ and degree $1$ generate every morphism in $\LBS$. The injectivity is a consequence of the \red{uniqueness} of decomposition \eqref{eqdecisonested}.
\end{proof}

\subsection{Definition of Koszulness for $\LBS$-operads}

\subsubsection{Odd operads over cubical Feynman categories}
Let $\F = (\fcat, \nucat, \imath)$ be cubical category. An odd operad over $\F$ is an operad over the Feynman category $\F^{\odd} = (\fcat^{\odd}, \nucat, \imath^{\odd})$ where $\fcat^{\odd}$ is the category enriched in abelian groups having the same objects as $\fcat$ and morphisms
\begin{equation*}
\fcat^{\odd}(X ,Y) = \Z<\bigsqcup_n C_n^{+}(X,Y)>/ \sigma.\phi - \epsilon(\sigma)\phi,
\end{equation*}
\red{where $\epsilon(\sigma)$ denotes the signature of $\sigma$.} The composition is given by the concatenation of chains of morphisms. Since $\nucat$ only has isomorphisms it is clear that $\nucat$ is also embedded in $\fcat^{\odd}$ and we call this embedding $\imath^{\odd}$. \\

In the case of $\LBS$, the category $\LBScat^{\odd}$ is generated by isomorphisms and generators $\red{[G]}^{\odd}$ for each element $G$ which is not the maximal element in some building set of some lattice, quotiented by relations
\begin{equation}\label{eqrelantiodd}
\red{[G_1]}^{\odd} \circ ([G_1 \vee G_2]^{\odd}\otimes \Id) = - \red{[G_2]}^{\odd}\circ ([G_1 \vee G_2]^{\odd}\otimes \Id)\circ \sigma_{2,3}
\end{equation}
for every nested antichain $\{G_1 \neq \un, G_2 \neq \un \}$ in some building set, relations
\begin{equation}\label{eqrelchainodd}
\red{[G_1]}^{\odd}\circ (\red{[G_2]}^{\odd}\otimes \Id) = - \red{[G_2]}^{\odd}\circ (\Id \otimes \red{[G_1]}^{\odd})
\end{equation}
for every chain $G_1 < G_2<\un$ in some building set, as well as relations
\begin{equation}\label{eqrelisoodd}
f \circ \red{[f(G)]}^{\odd} = [G]^{\odd} \circ (f_{[G, \un]} \otimes f_{[\zero, G]})
\end{equation}
for every isomorphism $f$ between built lattices. \\

\subsubsection{An example of an odd $\LBS$ cooperad}
The family of projective Orlik--Solomon algebras $\{\OSbara(\L)\}_{(\L, \G)}$ has an odd cooperadic structure over $\LBS$ which extends the dual of the odd operad $\mathsf{Grav}$. It will be denoted by $\OSbar$ and defined as follow\red{s}.
\begin{itemize}
\item For any built lattice $(\L, \G)$ we define
\begin{equation*}
\OSbar(\L, \G) \coloneqq \OSbara(\L)[-1],
\end{equation*}
where $[-1]$ is the shift in grading such that the unit of $\OSbara(\L)$ has grading $1$.
\item For any element $G \in \G \setminus \{\un\}$ in some irreducible built lattice $(\L, \G)$, we define
\begin{equation*}
\begin{array}{rccc}
\OSbar(\red{[G]}^{\odd}):& \OSbara(\L)[-1]& \longrightarrow &\OSbara([G, \un])[-1]\otimes \OSbara([\zero, G])[-1] \\
&\prod_i e_{H_i}\prod_{j} e_{H'_j} &\longrightarrow & \delta(\prod_i e_{G\vee H_i}) \otimes \prod_{j} e_{H'_j}
\end{array}
\end{equation*}
where the $H_i$'s are atoms not below $G$ and the $H'_j$'s are atoms below $G$.
\end{itemize}
Thanks to the shift in grading, those morphisms are morphisms of graded vector spaces. Let us check that the image of $\OSbar(\red{[G]}^{\odd})$ lands in the tensor of projective Orlik--Solomon algebras. By Lemma \ref{lemmaosbar} it is enough to show that $\OSbar(\red{[G]}^{\odd})(\delta(\prod_{\H}e_{H}))$ belongs to $\OSbara([G,\un])\otimes \OSbara([\zero, G])$ for all sets of atoms $\H$. We partition $\H$ into $\{H_i \} \sqcup \{H'_j\}$ where the $H_i$'s are atoms not below $G$ and the $H'_j$'s are atoms below $G$. We then have
\begin{align*}
\OSbar(\red{[G]}^{\odd})(\delta(\prod_{\mathcal{H}}e_H)) &= \OSbar(\red{[G]}^{\odd})(\delta(\prod_i e_{H_i}) \prod_j e_{H'_j} \pm \prod_i e_{H_i} \delta(\prod_j e_{H'_j})) \\
&= \delta(\delta(\prod_i e_{G\vee H_i}))\otimes \prod_j e_{H'_j} \pm \delta(\prod_i e_{G \vee H_i}) \otimes \delta(\prod_j e_{H'_j}) \\
&= \pm \delta(\prod_i e_{G \vee H_i}) \otimes \delta(\prod_j e_{H'_j}) \in \OSbara([G, \un])\otimes \OSbara([\zero, G]).
\end{align*}
\begin{itemize}
\item For any isomorphism of built lattice\red{s} $f:(\L', \G') \xrightarrow{\sim} (\L, \G)$ we define $\OSbar(f)$ as the restriction of $\OS(f)$ to the projective subalgebra.
\end{itemize}

We must check that the morphisms $\OSbar(\red{[G]}^{\odd})$ and $\OSbar(f)$ satisfy relations \eqref{eqrelantiodd}, \eqref{eqrelchainodd} and \eqref{eqrelisoodd} above. Let $\{ G_1 \neq \un, G_2 \neq \un \}$ be a nested antichain in some irreducible built lattice $(\L, \G)$ and let $\alpha = \delta(\prod_{i \leq n}e_{H_i}\prod_{j\leq n'} e_{H'_j}\prod_{k\leq n''} e_{H''_k})$ be an element in $\OSbara(\L)$ where the $H_i$'s are atoms below neither $G_1$ nor $G_2$, the $H'_j$'s atoms below $G_2$ and the $H''_k$'s below $G_1$ (by nestedness of $\{G_1, G_2\}$ there can be no atom below both $G_1$ and $G_2$). In \textcolor{black}{that} case one can check that the morphism $$(\OSbar(\red{[G_1\vee G_2]}^{\odd})\otimes \Id) \circ \OSbar(\red{[G_1]}^{\odd})$$ sends $\alpha$ to $\delta(\prod_i e_{G_1 \vee G_2 \vee H_i})\otimes \delta(\prod_j e_{H'_j})\otimes \delta(\prod_k e_{H''_k})$ whereas the morphism $$\sigma_{2,3}\circ (\OSbar(\red{[G_1\vee G_2]}^{\odd})\otimes \Id) \circ \OSbar(\red{[G_2]}^{\odd})$$ sends $\alpha$ to the opposite. Let $G_1 < G_2 < \un$ be a chain in some irreducible built lattice $(\L, \G)$ and let $\alpha = \delta(\prod_{i \leq n}e_{H_i}\prod_{j\leq n'} e_{H'_j}\prod_{k\leq n''} e_{H''_k})$ be an element in $\OSbara(\L)$ where the $H_i$'s are atoms not below $G_2$, the $H'_j$'s are atoms below $G_2$ and not below $G_1$ and the $H''_k$'s are atoms below~$G_1$. In \textcolor{black}{that} case one can check that the morphism $$(\OSbar(\red{[G_2]}^{\odd})\otimes \Id) \circ \OSbar(\red{[G_1]}^{\odd})$$ sends $\alpha$ to $\delta(\prod_i e_{G_2\vee H_i})\otimes \delta(\prod_j e_{G_1\vee H'_j})\otimes \delta(\prod_k e_{H''_k})$ whereas the morphism $$(\Id \otimes \OSbar(\red{[G_2]}^{\odd})) \circ \OSbar(\red{[G_2]}^{\odd})$$ sends $\alpha$ to the opposite. \red{Finally}, equation \eqref{eqrelisoodd} is also easily verified. \\

To conlude, we have shown that $\OSbar$ is an odd $\LBS$-cooperad (in graded abelian groups).

\subsubsection{The bar/cobar construction}

Let $\C$ be some complete cocomplete symmetric monoidal abelian category. We denote by $\Ch \C$ the category of chain complexes over $\C$. Let $\F$ be a cubical Feynman category. \\

In \cite{kaufmann_feynman_2017}, Kaufmann and Ward define a bar operator
\begin{equation*}
\Bar: \F-Ops_{\Ch \C} \rightarrow \F^{\odd}-Ops_{\Ch \C^{op}}
\end{equation*}
and a cobar operator
\begin{equation*}
\Omega: \F^{\odd}-Ops_{\Ch \C^{op}} \rightarrow \F-Ops_{\Ch \C}
\end{equation*}
which form an adjunction $\Omega \rightleftharpoons \Bar$ and such that the counit
\begin{equation*}
\Omega \Bar \Longrightarrow \Id
\end{equation*}
is a level-wise quasi-isomorphism. Informally those functors are defined by taking free constructions together with a differential constructed using the degree $1$ generators. Let us describe $\Omega$ explicitly in our case. \\

Let $\P$ be an $\LBS^{\odd}$ cooperad in $\Ch \C^{op}$. We have
\begin{equation*}
\Omega (\P) = (\LBS(\P), \d_{\Omega} + \d_{\P}),
\end{equation*}
where $\d_{\P}$ is the obvious differential coming from $\P$ and $\d_{\Omega}$ is defined as follow\red{s}. Recall the explicit formula
\begin{equation*}
\LBS(\P)(\L, \G) = \bigoplus_{\otimes (\L_i, \G_i) \xrightarrow{f} (\L, \G)} \bigotimes_i \P(\L_i, \G_i) / \sim
\end{equation*}
where $\sim$ identifies components corresponding to equivalent maps (maps that can be obtained from each other by precomposition of isomorphisms). For any $\otimes_i p_i \in \bigotimes_{i\leq n} \P(\L_i, \G_i)$ we put
\begin{equation*}
\d_{\Omega} ([\otimes p_i, f:\bigotimes_i (\L_i, \G_i) \rightarrow (\L, \G)]) = \sum_{\substack{j \leq n \\ G \in \G_j\setminus\{\un\}}}[(\Id \otimes \P(\red{[G]}) \otimes \Id)(\otimes p_i), f\circ (\Id \otimes \red{[G]} \otimes \Id))].
\end{equation*}
The cubicality condition and the fact that $\P$ is an odd cooperad ensures that $\d_{\Omega} + \d_{\P}$ is indeed a differential. Let us now describe $\Bar$ explicitely for $\LBS$-operads. Let $\P$ be an $\LBS$-operad in $\Ch \C$. We have
\begin{equation*}
\Bar(\P) = (\LBS^{\odd}(\P), \d_{\Bar} + \d_{\P}),
\end{equation*}
where $\d_{\P}$ is the obvious differential coming from $\P$ and $\d_{\Bar}$ is defined as follow\red{s}. We have the formula
\begin{equation*}
\LBS^{\odd}(\P)(\L,\G) = \bigoplus_{\substack{n \in \N \\ \psi\in C^+_{n}( (\L', \G'), (\L, \G))}} \P(\L',\G') /\sim\, ,
\end{equation*}
where the equivalence relation $\sim$ is given by 
\begin{equation*}
(\P(f) (\alpha) , \psi) \sim (\alpha,\psi \circ f) 
\end{equation*} 
for every isomorphism $f$ and 
\begin{equation*}
(\alpha , \psi) \sim \epsilon(\sigma) (\alpha, \sigma.\psi)
\end{equation*}
for every permutation $\sigma$.  Let $\psi$ be an element in $C^{+}_n((\L', \G'), (\L, \G))$, and let $\alpha$ be an element of  $\P(\L',\G')$. We have
\begin{equation*}
\d_{\Bar}((\alpha, [\psi])) \coloneqq \sum_{\substack{\phi: (\L'', \G'') \rightarrow (\L, \G) \\ G \textrm{ s.t. } c(\phi \circ \red{[G]}) = c(\psi)}}\epsilon(\phi)((\Id\otimes \P(\red{[G]}) \otimes \Id)(\alpha) , \phi),
\end{equation*}
where $c$ is the composition map and $\epsilon(\phi)$ is the signature of the permutation sending $\psi$ to $\phi\circ \red{[G]}$. One can check that this descends to the quotient by the equivalence relation $\sim$.

\subsubsection{Quadratic duality and Koszul duality}

In this subsection we assume that $\C$ is a category of vector spaces over some field. For any graded Feynman category $\F$, an $\F$-quadratic data is a pair $(M, R)$ with $M$ an $\F$-module and $R$ a submodule of $\F_{1}(M)$, where $\F_{1}(M)$ denotes the part of grading $1$ in the free $\F$-operad $\LBS(M)$. Notice that an $\F$-quadratic data can also be seen as an $\F^{\odd}$-quadratic data since $\F$ and $\F^{\odd}$ have the same modules and we have
$$ \F^{\odd}_1(M) = \F_1 (M). $$

If $\P$ is an $\F$-operad which is a quotient $\F(M)/\langle R\, \rangle $ for some $\F$-quadratic data $(M,R)$, we define $$\P^{!} \coloneqq \F^{\odd}(M^{\vee})/\langle R^{\bot}\rangle,$$ which is an $\F^{\odd}$-operad.  We have a morphism of differential graded $\F^{\odd}$-operads
\begin{equation}\label{eqkoszulmor}
(\P^{!})^{\vee} \rightarrow \Bar\P,
\end{equation}
which is induced by the morphism of $\F$-modules given by the composition
\begin{equation*}
(\P^{!})^{\vee} \twoheadrightarrow M^{\vee} \hookrightarrow \P.
\end{equation*}
We say that $\P$ is Koszul with Koszul dual $(\P^{!})^{\vee}$ if morphism \eqref{eqkoszulmor} is a quasi-isomorphism. We refer to \cite{kaufmann_feynman_2017} and \cite{Ward_2020} for more details. This coincides with the classical Koszul duality theories (Koszul duality for operads for instance). \\

In addition to having a homological degree (given by the grading of $\LBS_{\sh}$), the odd cooperad $\Bar \P$ has a weight grading coming from the grading of $\P$. The differential preserves \textcolor{black}{that} weight grading. One can check that the map
\begin{equation*}
(\P^{!})^{\vee} \rightarrow \Bar\P
\end{equation*}
is injective and its image is exactly the kernel of $\d_{\Bar}$ in the diagonal $$\{\textrm{weight grading = degree}\},$$ which is also the homology of the diagonal since the elements on the diagonal are the highest degree elements in their respective weight component. As a consequence $\P$ is Koszul if and only if the homology of $\Bar \P$ is concentrated on the diagonal.

\subsection{Koszulness of $\FY^{\PD}$ using the projective combinatorial Leray model}
In \cite{Bibby_2021} the authors define a differential bigraded algebra $B(\L, \G)$ as follow\red{s}.
\begin{madef}[Projective combinatorial Leray model \cite{Bibby_2021}]
Let $(\L, \G)$ be an irreducible built lattice. The differential bigraded algebra $B(\L, \G)$ is defined as the quotient of the free commutative algebra $\Q[e_G, x_G, G\in \G]$ by the ideal $\I$ generated by 
\begin{enumerate}
\item The elements $e_{\S}x_{\Tree}$ with $\S \cup \Tree$ not nested. 
\item The elements $\sum_{G\geq H}x_G$ for all atoms $H$ of $\L$.
\item The element $e_{\un}$. 
\end{enumerate}
The generators $e_G$ have bidegree $(0, 1)$ and the generators $x_G$ have bidegree $(2,0)$. The differential $\d$ of \textcolor{black}{that} algebra has bidegree $(2,-1)$ and is defined by 
\begin{align*}
&\d(e_G) = x_G \\
&\d(x_G) = 0.
\end{align*}
\end{madef}
The authors of \cite{Bibby_2021} have shown that we have isomorphisms of graded vector spaces 
$$B^{\bullet, d}(\L, \G) \simeq \bigoplus_{\substack{\S \textrm{ spanning } \\ \textrm{nested set of } (\L, \G)  \\ \#\S = d +1 }}\FY^{\PD}(\S)$$
for every integer $d$ (\cite{Bibby_2021} Proposition 5.1.4). Those isomorphisms give an isomorphism of complexes between $B(\L, \G)$ and $\Bar \FY^{\PD}(\L, \G)$.  For each irreducible built lattice $(\L, \G)$ we have a morphism of differential graded algebras
$$\OSbara(\L) \rightarrow B^{\bullet}(\L, \G),$$ 
induced by the map $e_H \rightarrow \sum_{G\geq H}e_G$. One can check that this is a morphism of $\LBS^{\odd}$-cooperads. The main result of \cite{Bibby_2021} is the following. 
\begin{theo}[\cite{Bibby_2021}, Theorem 5.5.1]
The morphism $\OSbara(\L) \xrightarrow{\sim} B^{\bullet}(\L, \G)$ is a quasi-isomorphism for every pair $(\L, \G)$. 
\end{theo}
This immediately implies: 
\begin{coro}
The operad $\FY^{\PD}$ is Koszul with Koszul (co)dual $\OSbar$.
\end{coro}
\begin{rmq}
The algebra structure on $\Bar\FY^{\PD}(\L, \G)$ coming from the isomorphism $$\Bar\FY^{\PD}(\L, \G) \simeq B(\L, \G)$$ can be defined purely operadically as follow\red{s}. Let $\alpha$ be some element in $\FY^{\PD}(\S)$ for some nested set $\S$ in some built lattice $(\L, \G)$  and $\beta$ some element of $\FY^{\PD}(\S')$ for some nested set $\S'$ in the same built lattice. The product of $\alpha$ and $\beta$ in $\Bar\FY^{\PD}(\L, \G)$ is given by 
\begin{equation*}
\alpha\cdot \beta = \left \{
\begin{array}{cl}
\FY(\S')(\alpha)\FY(\S)(\beta)& \textrm { if } \S \cap \S' = \{\un\} \textrm{ and } \S \cup \S' \textrm{ is a nested set}. \\
0  & \textrm{ otherwise.}
\end{array}
\right. 
\end{equation*}
In the first row, $\S'$ is viewed as a nested set of $(\L_{\S}, \G_{\S})$ and $\S$ is viewed as a nested set of $(\L_{\S'}, \G_{\S'})$ via Lemma \ref{lemmarecons}. The product takes place in the algebra $\FYa(\L_{\S\cup \S'}, \G_{\S \cup \S'}).$ It is interesting to note that we have used the operadic structure of $\FY$ and not that of $\FY^{\PD}$. This shows that Poincaré duality plays an important role when trying to relate the properties of $\FY$ and the properties of the Feichtner--Yuzvinsky algebras. 
\end{rmq}
\subsection{Koszulness and the affine Leray model}
In \cite{Bibby_2021} the authors also define a Leray model $\hat{B}(\L, \G)$ for the (affine) \textcolor{black}{Orlik--Solomon} algebras, just by taking out the relation $e_{\un} = 0$ in $B(\L, \G)$. One can also interpret this affine Leray model as a bar construction in a larger Feynman category $\LBSm$ defined as follow\red{s}. The set of objects of the underlying groupoid of $\LBSm$ is 
\begin{equation*}
\mathrm{Ob}(\LBSi) \sqcup \mathrm{Ob}(\LBSi).
\end{equation*}
For each irreducible built lattice $(\L, \G)$ we will denote the two copies of $(\L, \G)$ in $\LBSm$ by $(\L, \G)^{\proj}$ and $(\L, \G)^{\aff}$, for reasons which will be clear later. If $(\L, \G)$ is an irreducible built lattice, the structural morphisms of $\LBSm$ with target $(\L, \G)^{\proj}$ are labelled by spanning nested sets 
\begin{equation*}
\bigotimes_{G \in \S}([\tau_{\S}(G), G], \Ind(\G))^{\proj} \xrightarrow{\S} (\L, \G)^{\proj},
\end{equation*}
with composition as in $\LBS$. In other words when restricting $\LBSm$ to the ``projective'' arities we get the Feynman category $\LBS$. The structural morphisms of $\LBSm$ with target $(\L, \G)^{\aff}$ are labelled by nested sets which can either contain $\un$ or not. If $\S$ contains $\un$ then we have the morphism
\begin{equation*}
\bigotimes_{G \in \S}([\tau_{\S}(G), G], \Ind(\G))^{\proj} \xrightarrow{\S^{\aff}} (\L, \G)^{\aff}
\end{equation*}
and if $\S$ does not contain $\un$ we have the morphism
\begin{equation*}
\bigotimes_{G \in \S}([\tau_{\S}(G), G], \Ind(\G))^{\proj}\otimes ([\tau_{\S}(\un), \un], \Ind(\G))^{\aff} \xrightarrow{\S^{\aff}} (\L, \G)^{\aff}.
\end{equation*}
The composition of those morphisms is defined as in $\LBS$. \textcolor{black}{That} Feynman category encodes pairs $(\P, \Mop)$ with $\P$ an $\LBS$-operad and $\Mop$ a ``$\P$-module'' ($\P$ is the restriction to the projective part and $\Mop$ the restriction to the affine part). A set of generating morphisms of $\LBSm$ is given by 
\begin{equation*}
\red{[G]}^{\proj} \, \, (G \neq \un) \textrm{ and }  \red{[G]}^{\aff}.
\end{equation*}
One can define an odd $\LBSm$-cooperad $\OS_{\tot}$ by setting
\begin{equation*}
\OS_{\tot}((\L, \G)^{\proj}) = \OSbara(\L),\quad 
\OS_{\tot}((\L, \G)^{\aff}) = \OSa(\L),
\end{equation*}
for each irreducible built lattice $(\L, \G)$, and
\begin{equation*}
\OS_{\tot}(\red{[G]}^{\proj}) = \OSbar(\red{[G]})
\end{equation*}
together with
\begin{equation*}
\OS_{\tot}(\red{[G]}^{\aff}) = (\delta \otimes \Id)\circ \OS(\red{[G]}).
\end{equation*}
In \textcolor{black}{that} case the odd $\OSbar$-comodule structure on $\OS$ comes from the morphism of $\LBS$-cooperads $$\OS \xrightarrow{\delta} \OSbar.$$On the other hand one can define an $\LBSm$-operad $\FY^{\PD}_{\tot}$ by setting $$\FY^{\PD}_{\tot}((\L, \G)^{\proj}) = \FYa(\L, \G), \quad \FY^{\PD}_{\tot}((\L, \G)^{\aff}) = \FYa(\L, \G),$$ and 
$$\FY_{\tot}(\red{[G]}^{\proj}) = \FY^{\PD}(\red{[G]}^{\proj}),$$ together with
$$\FY_{\tot}(\red{[G]}^{\aff}) = \FY(\red{[G]})$$
for $G \neq \un$, and finally $\FY_{\tot}(\red{[\un]}^{\aff})$ is set to be the multiplication by $x_{\un}$. Exactly as for the projective part, one can use the results of \cite{Bibby_2021} to see that $\hat{B}(\L, \G)$ is isomorphic to $\Bar \FY^{\PD}_{\tot}$ and the morphism  $e_H \rightarrow \sum_{G\geq H}e_G$ induces a quasi-isomorphism of odd $\LBSm$-cooperads
\begin{equation*}
\OS_{\tot} \xrightarrow{\sim} \Bar\FY^{\PD}_{\tot},
\end{equation*}
which implies the Koszulness of the $\LBSm$-operad $\FY^{\PD}_{\tot}$.
\subsection{Koszulness via shuffle operads}
As in the case of classical operads and their shuffle counterpart, we have the key proposition.
\begin{prop}
Let $\P$ be an $\LBS$-operad. $\P$ is Koszul if and only if $\P_{\sh}$ is Koszul.
\end{prop}
\begin{proof}
By Proposition \ref{propmainshuffle} we have isomorphisms of shuffle $\LBS^{\odd}$-operads
\begin{equation*}
(\P_{\sh})^{!} \simeq (\P^{!})_{\sh}
\end{equation*}
which gives an isomorphism of shuffle $\LBS^{\odd}$-cooperad\red{s} (in $\Ch \C$):
\begin{equation*}
((\P_{\sh})^{!})^{\vee} \simeq ((\P^{!})^{\vee})_{\sh}.
\end{equation*}
On the other hand, we have
\begin{equation*}
(\LBS^{\odd}(\P))_{\sh} \simeq \LBS^{\odd}_{\sh}(\P_{\sh}).
\end{equation*}
By going back to explicit formulas one can check that those isomorphisms are compatible with the Bar differential and that we have a commutative diagram
\begin{equation*}
\begin{tikzcd}
(\P^{!})^{\vee}_{\sh}  \arrow[r] \arrow[d,"\simeq"]
& (\Bar\P)_{\sh} \arrow[d, "\simeq" ] \\
((\P_{\sh})^{!})^{\vee} \arrow[r]
& \Bar(\P_{\sh}),
\end{tikzcd}
\end{equation*}
but of course in every arity $(\L, \G, \vartriangleleft)$ we have the commutative diagram of complexes
\begin{equation*}
\begin{tikzcd}
(\P^{!})^{\vee}(\L, \G) \arrow[r] \arrow[d, "\simeq"]
& \Bar\P(\L, \G) \arrow[d, "\simeq"] \\
((\P^{!})^{\vee})_{\sh}(\L, \G, \vartriangleleft)  \arrow[r]
& (\Bar \P)_{\sh}(\L, \G, \vartriangleleft).
\end{tikzcd}
\end{equation*}
Combining the two diagrams in every arity finishes the proof.
\end{proof}

\subsection{Koszulness and Gröbner bases}

As in the case of classical shuffle operads, we have the key proposition.
\begin{prop}
Let $\P$ be a shuffle $\LBS$-operad. If $\P$ admits a quadratic Gröbner basis then $\P$ is Koszul.
\end{prop}
\begin{proof}
This is just an adaptation of the proof given in \cite{Hoffbeck_2010} to our setting, and translating Gröbner basis language in ``PBW basis'' language. We denote $\P = \LBS_{\sh}(M)/ \langle R \, \rangle $.\\

Let us use our total well-order on monomials to construct a filtration on $\Bar \P$. We will denote this total order by ``$<$''. Let $m = (\S, (e_G)_{G\in \S})$ be some monomial. We define
\begin{equation*}
\Fil_m \Bar \P = < \{ m_1\otimes ... \otimes m_n \in \P(\S') \, | \,  m_1,... \,, m_n \textrm{ monomials s.t. } \LBS_{\sh}(\S')((m_i)_i) \leq m \} >
\end{equation*}
where the brackets $<,>$ denote the linear span. The bar differential preserves \textcolor{black}{that} filtration. We will now show that the associated spectral sequence collapses at the first page and its homology is concentrated on the diagonal. The complex $E^0_m \Bar \P$ is spanned by elements of the form $$\P(\S_1)((e_G)_{G \in \S_1})\otimes ... \otimes \P(\S_n)((e_G)_{G\in \S_n}),$$ where the nested sets $\S_i$ are such that there exist some nested set $\S'$ satisfying $\S = \S' \circ (\S_i)_i$, and such that the monomials $\P(\S_i)((e_G)_{\G \in \S_i}$ are all normal. \\

For any $G$ in $\S\setminus\{\un\}$ we denote by $\n(G)$ the unique minimum of $\S_{>G}$. We also denote by $\Adm(m)$ the set of elements of $\S\setminus\{\un\}$ such that $\P(\red{[G]})(e_G, e_{\n(G)})$ is a normal monomial. By the fact that our Gröbner basis is quadratic we see that $E^0_m \Bar \P$ is isomorphic to the augmented dual of the combinatorial complex $C_{\bullet}(\Delta_{\Adm(m)})$, which has trivial homology except when $\Adm(m) = \emptyset$, in which case the complex is reduced to $\K$ on the diagonal (with generator given by $\otimes_G e_G$). By a standard spectral sequence argument this concludes the proof.
\end{proof}
As a corollary of \textcolor{black}{that} proposition and \ref{coroquadgrob} we get 
\begin{coro}
The operad $\FY^{\vee}$ is Koszul.
\end{coro}

\section{Further directions}\label{secdirec}
In this section we highlight some possible ways to extend/refine $\LBS$ which seem natural to us and may lead to further applications.
\subsection{Working with matroids instead of geometric lattices}
One possible refinement of $\LBS$ would be to do everything with matroids instead of geometric lattices, which would allow us to take loops and parallel elements into account. For now this refinement is useless because all the operads we know (Feichtner--Yuzvinsky rings, Orlik--Solomon algebras) do not ``see'' the loops and parallel elements (i.e. factor through the lattice of flat construction). However, it may happen that some finer invariants of matroids which detect loops and parallel elements may also have an operadic structure. In order to implement this refinement it will be beneficial to have a purely matroidal axiomatization of building sets. Let us describe a possible way to obtain that. Recall that a matroid can be defined by its rank function as follow\red{s}.
\begin{madef}[Matroid, via rank function]
Let $E$ be a finite set. A \emph{matroid} structure on $E$ is the datum of a map
$$ \rk : \mathcal{P}(E) \rightarrow \N $$
called the rank function, satisfying the following properties.
\begin{enumerate}
\item The rank function takes value $0$ on the empty set.
\item For every $A, B \in \mathcal{P}(E)$ we have $$ \rk(A \cup B) + \rk(A \cap B) \leq \rk(A) + \rk(B).$$
\item For every $A \in \mathcal{P}(E)$ and $x \in E$ we have $$\rk(A \cup \{x\}) \leq \rk(A) +1.$$
\end{enumerate}
\end{madef}
Here is a possible way of axiomatizing a building set in terms of the rank function.
\begin{madef}[Building decomposition]
A \emph{building decomposition} of a matroid $(E, \rk)$ is a function $\nu$ which assigns to every subset $X\subset E$ a partition of $X$ and which satisfies the following axioms.
\begin{enumerate}
\item If $X\subset Y$ are two subsets of $E$, then $\nu(X)$ refines the restriction of $\nu(Y)$ to $X$.
\item If $\nu(X)$ is the partition with blocks $P_1| ... |P_n$ then for all $i\leq n$ the partition $\nu(P_i)$ is the partition with only one block.
\item For all $X\subset E$, if $\nu(X)$ is the partition $P_1| ... | P_n$ then $\rk(X) = \rk(P_1) + ... + \rk(P_n)$.
\end{enumerate}
\end{madef}

On simple loopless matroids the datum of a building decomposition is equivalent to the datum of a building set on the lattice of flats. One can construct a building set out of a building decomposition by considering the flats which have a partition with only one block. On the other hand, one can construct a building decomposition out of a building set by setting $\nu(X)$ to be the partition induced by the factor decomposition of $\sigma(X)$. 

\begin{ex}
Let $\textcolor{black}{\Gamma} = (V ,E)$ be a graph and $M_{\textcolor{black}{\Gamma}}$ its cycle matroid. $M_{\textcolor{black}{\Gamma}}$ admits a building decomposition given by the partitions into connected components for each subset of $E$. Naturally if we look at the induced building set on the lattice of flats this gives the graphical building set introduced in Example \ref{exbs}.
\end{ex}

We also have a natural notion of induced building decomposition on restrictions and contractions of matroids.
\begin{madef}
Let $M = (E,\rk)$ be a matroid with building decomposition $\nu$ and $S$ a subset of the ground set $E$. The contraction $M^{S}$ admits a building decomposition $\Ind^{S}(\nu)$ given by $\Ind^S(\nu)(X) = \nu(X \cup S)_{|X}$ for every $X \subset E\setminus S$. The restriction $M_{S}$ admits a building decomposition $\Ind_S(\nu)$ given by $\Ind_S(\nu)(X) = \nu(X)$ for every $X \subset S$.
\end{madef}

Working with those definitions, we are fairly certain everything should work in the same fashion as in Section \ref{secfeycat}, by just replacing the built lattices $([\zero, G], \Ind(\G)), ([G,\un], \Ind(\G))$ by the matroidal restrictions/contractions $(M_{G}, \Ind(\nu)), (M^{G}, \Ind(\nu))$, for $G$ such that $\nu(G)$ is the trivial partition (a priori we would not even need $G$ to be closed, which would give additional structural morphisms). 

\subsection{The polymatroidal generalization}

One can also naturally consider an exension of $\LBS$ to polymatroids, which form a combinatorial abstraction of subspace arrangements. This is justified by the fact that the wonderful compactification story also works for subspace arrangements and the cohomology algebras give us a natural candidate for an operad over \textcolor{black}{that} bigger Feynman category. It has been shown by Pagaria and Pezzoli \cite{Pagaria_2021} that those cohomology rings also admit natural generalizations to arbitrary polymatroids and that they also have a Hodge theory. Here are some reminders on polymatroids.
\begin{madef}[Polymatroid]
Let $E$ be a finite set. A \emph{polymatroid} structure on $E$ is the datum of a function $$\cd: \mathcal{P}(E) \rightarrow \N $$ satisfying
\begin{enumerate}
\item $\cd(\emptyset) = 0$.
\item For any subsets $A \subset B$ of $E$ we have $\cd(A) \leq \cd(B)$.
\item For any subsets $A, B$ of $E$ we have $$\cd(A\cap B) + \cd(A \cup B) \leq \cd(A) + \cd(B).$$
\end{enumerate}
\end{madef}
The letters ``cd'' stand for codimension. If we ask that $\cd$ take value $1$ on singletons we get a classical matroid. We can define the lattice of flats of a polymatroid by considering the subsets $F$ of $E$ such that $\cd(F \cup \{x\}) > \cd(F)$ for all $x \notin F$. However, for general polymatroids the lattice of flats does not contain enough information and needs to be considered together with $\cd$ to recover the polymatroid (for matroids ``$\cd$'' is just the rank function of the lattice of flats and does not add any information). In \cite{Pagaria_2021} the authors introduced a notion of building set for polymatroids.
\begin{madef}
Let $P = (E, \cd)$ be a polymatroid with lattice of flats $\L$. A \emph{building set} of $P$ is a subset $\G$ of $\L\setminus\{\zero\}$ such that for any $X$ in $\L$ the join gives an  isomorphism of posets $$ \prod_{G \in \Fact_{\G}(X)}[\zero, G] \xrightarrow{\sim} [\zero, X] $$ and we additionally have
$$\cd(X) = \sum_{G \in \Fact_{\G}(X)} \cd(G).$$
\end{madef}
Notice that the last condition is automatically verified for matroids ($\cd = \rk$). The authors also give suitable generalizations of nested sets, and they show that for any $G$ in $\L$, the (polymatroidal) contraction $([G, \un], \Ind(\cd))$ has an induced (polymatroidal) building set given (as in the matroidal case) by  $$\Ind_{[G, \un]}(\G)  = (\G \vee G) \cap (G, \un],$$ and the same goes for (polymatroidal) restrictions. With those definitions we are fairly certain one can readily extend $\LBS$ to polymatroids. \\

In \cite{Pagaria_2021}, the authors also introduce a generalization of the Feichtner--Yuzvinsky algebras to the polymatroidal setting as follows.
\begin{madef}
Let $(\L, \cd)$ be a polymatroid with some building set $\G$. The algebra $\FYa(\L, \G, \cd)$ is defined by $$\FYa(\L, \G, \cd) = \mathbb{\Q}[x_G, \, G \in \G]/\I,$$ with all the generators in degree $2$ and $\I$ the ideal generated by elements $$\prod_{i\leq n}x_{G_i} $$ where $\{G_1, ..., G_n\}$ is not nested and elements $$ \left(\sum_{G \geq H}x_{G}\right)^{\cd(H)}$$ for any atom $H$.
\end{madef}
As in the matroidal case, one can get another presentation by considering the change of variable $h_{G} \coloneqq \sum_{G' \geq G} x_{G'}$. The algebra morphisms
\begin{equation*}
\begin{array}{ccc}
\FYa(\L, \cd, \G) &\longrightarrow &\FYa([G, \un], \Ind(\cd), \Ind_{[G, \un]}(\G)) \otimes \FYa([\zero, G], \Ind(\cd), \Ind_{[\zero, G]}(\G)) \\
h_{G'} &\longrightarrow &
\left\{
\begin{array}{cc}
h_{G\vee G'} \otimes 1 &\textrm{ if } G' \nleq G \\
1 \otimes h_G & \textrm{ otherwise.} 
\end{array}
\right.
\end{array}
\end{equation*}
are well-defined and give us an operadic structure on the family of generalized Feichtner--Yuzvinsky algebras.

\subsection{Adding morphisms to $\LBS$}

One could also consider adding more morphisms of degree $0$ in $\LBS$. This is justified by the fact that the Feichtner--Yuzvinsky algebras have a lot more functoriality than what we have in $\LBS$. More precisely let $(\L, \G)$ and $(\L', \G')$ be two built lattices and let $f:\L \longrightarrow \L'$ be a poset morphism which sends $\G$ to $\G'$, atoms of $\L$ to atoms of $\L'$ and which is compatible with the join on both sides, i.e.
\begin{equation*}
f(G_1\vee G_2) = f(G_1) \vee f(G_2)
\end{equation*}
for all $G_1, G_2$ in $\L$. With those hypotheses the map induced by
\begin{equation*}
\begin{array}{ccc}
\FYa(\L, \G) &\xrightarrow{\FYa(f)} &\FYa(\L', \G') \\
h_G &\longrightarrow & h_{f(G)}
\end{array}
\end{equation*}
is a well-defined map of algebras. This incentivizes us to formally add such morphisms in $\LBS$. Some of those morphisms are very natural to add in their own right. For instance if $\G \subset \G'$ are two building sets of some lattice $\L$ then the identity of $\L$ satifies the above conditions. In the realizable case the corresponding map $\FYa(f)$ is induced by the blow down $$\overline{Y}_{\L, \G'} \rightarrow \overline{Y}_{\L, \G}.$$

If $f$ is the inclusion of some interval $[\zero, G] \hookrightarrow \L$ then $f$ satisfies the above conditions when taking the induced building set on $[\zero, G]$. For instance this includes the various inclusions $\Pi_n \simeq [\zero, \llbracket 1, n+1 \rrbracket\setminus \{i\} | i ]\hookrightarrow \Pi_{n+1}$. The corresponding morphisms $$\FYa^{\vee}(\Pi_{n+1}, \G_{\min}) \rightarrow \FYa^{\vee}(\Pi_{n}, \G_{\min})$$
are induced by forgetting some marked point on the genus 0 curve.

\bibliography{ArticleFeyCat}

\begin{thebibliography}{10}

\bibitem{Huh_2018}
Karim Adiprasito, June Huh, and Eric Katz.
\newblock Hodge theory for combinatorial geometries.
\newblock {\em Ann. of Math. (2)}, 188(2):381--452, 2018.

\bibitem{Backman_Spencer_Eur_2020}
Spencer Backman, Christopher Eur, and Connor Simpson.
\newblock Simplicial generation of {C}how rings of matroids.
\newblock {\em S\'{e}m. Lothar. Combin.}, 84B:Art. 52, 11, 2020.

\bibitem{Batanin_2015}
Michael Batanin and Martin Markl.
\newblock Operadic categories and duoidal {D}eligne's conjecture.
\newblock {\em Advances in Mathematics}, 285:1630--1687, 2015.

\bibitem{BW_1993}
Thomas Becker and Volker Weispfenning.
\newblock {\em Gröbner Bases, A Computational Approach to Commutative Algebra}.
\newblock Springer, 1993.

\bibitem{Bibby_2021}
Christin Bibby, Graham Denham, and Eva~Maria Feichtner.
\newblock A {L}eray model for the {O}rlik{\textendash}{S}olomon algebra.
\newblock {\em International Mathematics Research Notices}, sep 2021.

\bibitem{Borisov2008}
Dennis~V. Borisov and Yuri~I. Manin.
\newblock {\em Generalized Operads and Their Inner Cohomomorphisms}, pages 247--308.
\newblock Birkh{\"a}user Basel, Basel, 2008.

\bibitem{Braden_2022}
Tom Braden, June Huh, Jacob~P. Matherne, Nicholas Proudfoot, and Botong Wang.
\newblock A semi-small decomposition of the {C}how ring of a matroid.
\newblock {\em Advances in Mathematics}, 409:108646, 2022.

\bibitem{Carr_Devadoss_2004}
Michael~P. Carr and Satyan~L. Devadoss.
\newblock Coxeter complexes and graph-associahedra.
\newblock {\em Topology Appl.}, 153(12):2155--2168, 2006.

\bibitem{Cohen_1973}
Fred Cohen.
\newblock Cohomology of braid spaces.
\newblock {\em Bull. Amer. Math. Soc.}, 79:763--766, 1973.

\bibitem{de_concini_wonderful_1995}
Corrado De~Concini and Claudio Procesi.
\newblock Wonderful models of subspace arrangements.
\newblock {\em Selecta Mathematica. New Series}, 1(3):459--494, 1995.

\bibitem{DK_2010}
Vladimir Dotsenko and Anton Khoroshkin.
\newblock Gr\"{o}bner bases for operads.
\newblock {\em Duke Math. J.}, 153(2):363--396, 2010.

\bibitem{Dotsenko_2012}
Vladimir Dotsenko and Anton Khoroshkin.
\newblock Quillen homology for operads via {G}r{\"o}bner bases.
\newblock {\em Documenta Mathematica}, 2012.

\bibitem{EHKR_2010}
Pavel Etingof, Andr\'{e} Henriques, Joel Kamnitzer, and Eric~M. Rains.
\newblock The cohomology ring of the real locus of the moduli space of stable curves of genus 0 with marked points.
\newblock {\em Ann. of Math. (2)}, 171(2):731--777, 2010.

\bibitem{FK_2004}
Eva-Maria Feichtner and Dmitry~N. Kozlov.
\newblock Incidence combinatorics of resolutions.
\newblock {\em Selecta Math. (N.S.)}, 10(1):37--60, 2004.

\bibitem{feichtner_chow_2003}
Eva~Maria Feichtner and Sergey Yuzvinsky.
\newblock Chow rings of toric varieties defined by atomic lattices.
\newblock {\em Invent. Math.}, 155(3):515--536, 2004.

\bibitem{Forcey_Ronco_2022}
Stefan Forcey and Mar\'{\i}a Ronco.
\newblock Algebraic structures on graph associahedra.
\newblock {\em J. Lond. Math. Soc. (2)}, 106(2):1189--1231, 2022.

\bibitem{Getzler_1994}
Ezra Getzler.
\newblock Operads and moduli spaces of genus {$0$} {R}iemann surfaces.
\newblock In {\em The moduli space of curves ({T}exel {I}sland, 1994)}, volume 129 of {\em Progr. Math.}, pages 199--230. Birkh\"{a}user Boston, Boston, MA, 1995.

\bibitem{Getzler2009}
Ezra Getzler.
\newblock {\em Operads Revisited}, pages 675--698.
\newblock Birkh{\"a}user Boston, Boston, 2009.

\bibitem{Hoffbeck_2010}
Eric Hoffbeck.
\newblock A {P}oincar\'{e}-{B}irkhoff-{W}itt criterion for {K}oszul operads.
\newblock {\em Manuscripta Math.}, 131(1-2):87--110, 2010.

\bibitem{KW2021}
Ralph Kaufmann and Benjamin Ward.
\newblock Koszul {F}eynman categories.
\newblock {\em Proceedings of the American Mathematical Society}, 08 2021.

\bibitem{kaufmann_feynman_2017}
Ralph~M. Kaufmann and Benjamin~C. Ward.
\newblock Feynman categories.
\newblock {\em Ast\'{e}risque}, (387):vii+161, 2017.

\bibitem{LV_2012}
Jean-Louis Loday and Bruno Vallette.
\newblock {\em Algebraic Operads}.
\newblock Springer, 2012.

\bibitem{LM_2000}
Andrey Losev and Yuri Manin.
\newblock {New moduli spaces of pointed curves and pencils of flat connections.}
\newblock {\em Michigan Mathematical Journal}, 48(1):443 -- 472, 2000.

\bibitem{Manin_1999}
Yuri~I. Manin.
\newblock {\em Frobenius Manifolds, Quantum Cohomology, and Moduli Spaces}.
\newblock American Mathematical Society, 1999.

\bibitem{OS_1980}
Peter Orlik and Louis Solomon.
\newblock Combinatorics and topology of complements of hyperplanes.
\newblock {\em Invent. Math.}, 56(2):167--189, 1980.

\bibitem{Pagaria_2021}
Roberto Pagaria and Gian~Marco Pezzoli.
\newblock Hodge theory for polymatroids, 2021.

\bibitem{Rains_2010}
Eric~M. Rains.
\newblock The homology of real subspace arrangements.
\newblock {\em J. Topol.}, 3(4):786--818, 2010.

\bibitem{wachs_poset_2006}
Michelle~L. Wachs.
\newblock Poset {Topology}: {Tools} and {Applications}, February 2006.
\newblock arXiv:math/0602226.

\bibitem{Ward_2020}
Benjamin~C. Ward.
\newblock Six operations formalism for generalized operads.
\newblock {\em Theory Appl. Categ.}, 34:Paper No. 6, 121--169, 2019.

\bibitem{welsh_matroid_1976}
Dominic J.~A. Welsh.
\newblock {\em Matroid theory}.
\newblock L. {M}. {S}. {Monographs}, {No}. 8. Academic Press [Harcourt Brace Jovanovich, Publishers], London-New York, 1976.

\bibitem{Yuzvinsky_2001}
Sergey Yuzvinski\u{\i}.
\newblock Orlik-{S}olomon algebras in algebra and topology.
\newblock {\em Uspekhi Mat. Nauk}, 56(2(338)):87--166, 2001.

\end{thebibliography}
\bibliographystyle{plain}
\end{document}